\documentclass[12pt]{article}
\usepackage{amsmath}
\usepackage{latexsym}
\usepackage{amssymb}
\newtheorem{thm}{Theorem}[section]
\newtheorem{la}[thm]{Lemma}
\newtheorem{Defn}[thm]{Definition}
\newtheorem{Remark}[thm]{Remark}
\newtheorem{prop}[thm]{Proposition}
\newtheorem{cor}[thm]{Corollary}
\newtheorem{Example}[thm]{Example}
\newtheorem{Examples}[thm]{Examples}
\newtheorem{Number}[thm]{\!\!}
\newenvironment{defn}{\begin{Defn}\rm}{\end{Defn}}
\newenvironment{example}{\begin{Example}\rm}{\end{Example}}
\newenvironment{examples}{\begin{Examples}\rm}{\end{Examples}}
\newenvironment{rem}{\begin{Remark}\rm}{\end{Remark}}
\newenvironment{numba}{\begin{Number}\rm}{\end{Number}}
\newenvironment{proof}{{\noindent\bf Proof.}}%
                  {\nopagebreak\hspace*{\fill}$\Box$\medskip\medskip\par}   
\newcommand{\Punkt}{\nopagebreak\hspace*{\fill}$\Box$}
\newcommand{\wb}{\overline}
\newcommand{\ve}{\varepsilon}

\newcommand{\wt}{\widetilde}
\newcommand{\tensor}{\otimes}

\newcommand{\isom}{\cong}
\newcommand{\mto}{\mapsto}
\newcommand{\N}{{\mathbb N}}
\newcommand{\A}{{\mathbb A}}
\newcommand{\R}{{\mathbb R}}
\newcommand{\C}{{\mathbb C}}

\newcommand{\K}{{\mathbb K}}

\newcommand{\MFD}{{\mathbb M}{\mathbb F}{\mathbb D}}
\newcommand{\TG}{{\mathbb T}{\mathbb G}}
\newcommand{\TOP}{{\mathbb T}{\mathbb O}{\mathbb P}}
\newcommand{\LIE}{{\mathbb L}{\mathbb I}{\mathbb E}}

\newcommand{\cO}{{\mathcal O}}
\newcommand{\cS}{{\mathcal S}}
\newcommand{\cT}{{\mathcal T}}
\newcommand{\cL}{{\mathcal L}}

\newcommand{\ch}{{\mathfrak h}}
\DeclareMathOperator{\Aut}{Aut}

\newcommand{\one}{{\bf 1}}
\newcommand{\sub}{\subseteq}
\DeclareMathOperator{\GL}{GL}
\DeclareMathOperator{\im}{im}

\DeclareMathOperator{\id}{id}

\newcommand{\sbull}{{\scriptscriptstyle \bullet}}

\DeclareMathOperator{\Diff}{Diff}
\DeclareMathOperator{\GermDiff}{GermDiff}
\DeclareMathOperator{\GermEnd}{GermEnd}
\DeclareMathOperator{\Supp}{supp}
\DeclareMathOperator{\Exp}{Exp}
\DeclareMathOperator{\Hol}{Hol}
\DeclareMathOperator{\Ad}{Ad}
\DeclareMathOperator{\ad}{ad}
\DeclareMathOperator{\ob}{ob}
\newcommand{\dl}{{\displaystyle\lim_{\longrightarrow}}}
\newcommand{\pl}{{\displaystyle\lim_{\longleftarrow}}}

\begin{document}
\begin{center}
{\Large\bf Direct Limits of Infinite-Dimensional\\[0.6mm]
Lie Groups Compared to
Direct Limits\\[2.6mm]
in Related Categories}\\[6mm]
{\bf Helge Gl\"{o}ckner}\vspace{4mm}
\end{center}
\begin{abstract}\vspace{1mm}
\hspace*{-7.2 mm}
Let\hspace*{-.15mm}
$G$\hspace*{-.15mm} be\hspace*{-.15mm} a\hspace*{-.15mm}
Lie\hspace*{-.15mm} group\hspace*{-.15mm} which\hspace*{-.15mm}
is\hspace*{-.15mm} the\hspace*{-.15mm} union\hspace*{-.15mm}
of\hspace*{-.15mm} an\hspace*{-.15mm}
ascending\hspace*{-.15mm} sequence\linebreak
$G_1\sub G_2\sub\cdots$\hspace*{-.15mm}
of\hspace*{-.15mm} Lie\hspace*{-.15mm} groups\hspace*{-.15mm}
(all\hspace*{-.15mm} of\hspace*{-.15mm} which\hspace*{-.15mm}
may\hspace*{-.15mm} be\hspace*{-.15mm} infinite-dimensional).
We study the question when
$G=\dl\; G_n$\vspace{-.3mm} in the category of
Lie groups, topological groups,
smooth manifolds, resp., topological spaces.
Full answers are obtained
for $G$ the group $\Diff_c(M)$
of compactly supported
$C^\infty$-diffeomorphisms of a $\sigma$-compact
smooth manifold~$M$; and
for test function groups
$C^\infty_c(M,H)$
of compactly supported smooth
maps with values in a finite-dimensional
Lie group~$H$.
We also discuss the cases where $G$
is a direct limit of unit groups
of Banach algebras, a Lie group
of germs of Lie group-valued analytic maps,
or a weak direct product of Lie groups.
\vspace{3mm}
\end{abstract}
{\footnotesize {\em Classification}:
22E65 (main); % inf-dim Lie groups
22E67, % loop gps and related
46A13, % FA: spaces defined by inductive limits
46F05, % spaces of test functions
46T20, % cts and diff'ble maps in non-lin FA
54B30, % categorical methods in gen topology
54H11, % relations of gen topology to top gps
58B10, % inf-dim mfds: diff questions
58D05.\\[2.5mm] % groups of diffeos
}\vspace{5mm}
\begin{center}
{\bf\Large Introduction}\vspace{.6mm}
\end{center}
It frequently happens that an infinite-dimensional
Lie group~$G$ of interest is a union
$G=\bigcup_{n\in \N}G_n$
of an ascending sequence $G_1\sub G_2\sub\cdots$
of Lie groups
which are easier to handle.
Typically, each $G_n$ is finite-dimensional,
a Banach-Lie group,
or at least a Fr\'{e}chet-Lie group,
while $G$ is modelled on a more complicated
locally convex space
(e.g., an LB-space, LF-space, or a Silva space).
Then good tools of infinite-dimensional
calculus are available
to establish differentiability properties
of mappings on the groups~$G_n$,
while properties of mappings (or homomorphisms) on~$G$
are more elusive and can be difficult to access.\\[3mm]
In this article, we consider
mappings $f\colon G\to X$,
where $X$ is a topological space or smooth manifold;
and also homomorphisms $f\colon G\to H$,
where $H$ is a topological group or Lie group.
We analyze the question
when continuity (or smoothness)
of $f|_{G_n}$ for each~$n$ implies
continuity (or smoothness) of~$f$.
To rephrase this problem
in category-theoretical terms,
let $\TOP$ be the category
of topological spaces and continuous maps,
$\TG$ the category of topological
groups and continuous homomorphisms,
$\MFD_\infty$ the category of smooth manifolds
(modelled on real locally convex spaces)
and smooth maps, and
$\LIE$ be the category
of Lie groups (modelled on real locally convex spaces)
and smooth homomorphisms.
We assume that $G=\bigcup_{n\in \N}G_n$
for an ascending sequence $G_1\sub G_2\sub\cdots$
of Lie groups,
such that all inclusion maps
$i_{n,m}\colon G_m\to G_n$ (for $m\leq n$)
and $i_n\colon G_n\to G$
are smooth homo\-morphisms.
Thus $\cS:=((G_n)_{n\in\N}, (i_{n,m})_{n\geq m})$
is a direct system in~$\LIE$,
and we can consider $\cS$ also as a direct system
in $\TOP$, $\TG$ and $\MFD_\infty$,
forgetting extraneous structure.
We are asking whether
$G=\dl\, G_n$\vspace{-.5mm}
(more precisely, $(G,(i_n)_{n\in \N})=\dl\,\cS$)
holds in $\TOP$, $\MFD_\infty$, $\TG$, resp., $\LIE$.\\[3mm]
If each $G_n$ is a finite-dimensional Lie group
and $G=\dl\,G_n$\vspace{-.5mm}
the direct limit Lie group
constructed in~\cite{FUN}
(see also \cite{DIR}, \cite{NRW1}
and \cite{NRW2} for special cases),
then $G=\dl\,G_n$\vspace{-.4mm}
holds in each of the preceding
categories \cite[Theorem~4.3]{FUN}.
The goal of this article is to shed light
on the case where the Lie groups~$G_n$
are infinite-dimensional.
Then the situation changes
drastically, and some of the direct
limit properties can get lost.\\[3mm]
We are interested both in general techniques
for the investigation of direct limit
properties, and a detailed analysis
of the properties of concrete groups.
In the following,
we summarize some of the main results.\\[5mm]
{\bf Tools for the identification
of direct limits.}
Direct limit properties of $G=\bigcup_{n\in \N}G_n$
are accessible provided that $G$ admits a
``direct limit chart'' composed of charts
of the groups~$G_n$ (see Definition~\ref{DLchart}
for details).
If this is the case, then
the following information becomes available:
\begin{itemize}
\item
If $G=\dl\,G_n$\vspace{-.3mm}
as a topological group, then $G=\dl\,G_n$
as a Lie group (Theorem~\ref{reduce});
\item
$G=\dl\,G$\vspace{-.3mm}
as a topological space if and only if
$L(G)=\dl\,L(G_n)$\vspace{-.3mm}
as a topological space
(Theorem~\ref{usefcor}\,(a));
\item
If
$L(G)$ is smoothly regular,
then $G=\dl\,G_n$\vspace{-.3mm}
as a smooth manifold if and only if
$L(G)=\dl\,L(G_n)$
as a smooth manifold
(Theorem~\ref{usefcor}\,(b)).\vspace{-.3mm}
\end{itemize}
Here, the existence of a direct limit chart is
a very weak requirement, which is satisfied
by all relevant examples known to the author.
Another result (Proposition~\ref{prodDL})
provides a criterion
(also satisfied by all relevant examples
inspected so far) which ensures
that $G=\dl\,G_n$\vspace{-.3mm}
as a topological group.
This criterion is closely related
to investigations in~\cite{TSH}.
In this paper, a condition was formulated
which facilitates a quite explicit
(``bamboo-shoot'')
description of the group topology
on a direct limit of topological groups.\\[5mm]
{\bf Direct limit properties of the prime examples.}
We now summarize our results concerning
concrete examples of direct limit groups.\\[4mm]
{\bf 1. \,Groups of compactly supported
functions or diffeomorphisms.}
Let $\Diff_c(M)$ be
the group of all $C^\infty$-diffeomorphisms $\gamma$
of a $\sigma$-compact, non-compact finite-dimensional
smooth manifold~$M$
which are compactly supported
(i.e., $\gamma(x)=x$ for all
$x$ outside some compact subset of~$M$).
Then $\Diff_c(M)$ is a Lie group
modelled on the LF-space $C^\infty_c(M,TM)$
of compactly supported smooth vector fields
(see \cite{Mic} or \cite{DIF}).
It is a union $\Diff_c(M)=\bigcup_K \Diff_K(M)$
of the Lie subgroups $\Diff_K(M)$
of diffeomorphisms supported in~$K$
(which are Fr\'{e}chet-Lie groups),
for $K$ ranging through the directed
set of compact subsets of~$M$.
Given a finite-dimensional
Lie group $H$, we are also interested
in the ``test function group''
$C^\infty_c(M,H)$ of compactly supported
smooth $H$-valued maps,
which is a Lie group modelled on the LF-space
$C^\infty_c(M,L(H))$ (see \cite{GCX}).
It is a union $C^\infty_c(M,H)=\bigcup_K C^\infty_K(M,H)$
of groups of mappings supported in a given compact set~$K$.
Now assume that
$M$ and~$H$ are non-discrete.
The following table
(compiled from Propositions~\ref{negdiff},
\ref{posdiff}, \ref{negtestf} and \ref{postestf})
describes in
which categories
$C^\infty_c(M,H)=\dl\, C^\infty_K(M,H)$,\vspace{-1.1mm}
resp., $\Diff_c(M)=\dl\,\Diff_K(M)$
holds:\vspace{2mm}
\begin{center}
\begin{tabular}{||c||c|c||}\hline\hline
category $\backslash$ group &  $C^\infty_c(M,H)$ 
& $\text{Diff}_c(M)$
\\ \hline\hline
Lie groups & yes & yes \\ \hline
topological groups & yes & yes \\ \hline
smooth manifolds & no & no \\ \hline
topological spaces & no & no \\ \hline\hline
\end{tabular}\vspace{5mm}
\end{center}
It was known before that
the direct limit topology does not
make $\Diff_c(M)$ a topological group
(see \cite[Theorem~6.1]{TSH}),
and Yamasaki's Theorem \cite[Theorem~4]{Yam}
implies that the direct limit topology does
not make $C^\infty_c(M,H)$ a topological group.
Hence, the Lie group topologies must
be properly coarser than the direct limit topologies
(as asserted in the last line of the table).
The other results
compiled in the table are new.\\[4mm]
\textbf{2. \,Weak direct products of Lie groups.}
The weak direct product
$\prod_{n\in \N}^*G_n$
of a sequence $(G_n)_{n\in \N}$ of Lie groups
is defined as the group of all
$(g_n)_{n\in \N}\in \prod_{n\in \N}G_n$
such that $g_n=1$ for all but finitely many~$n$;
it carries a natural Lie group structure~\cite[\S7]{MEA}.
We shall see that
$\prod_{n\in \N}^*G_n=\dl\,\prod_{k=1}^nG_k$\vspace{-.3mm}
as a topological group and as a Lie group
(Proposition~\ref{DLwprodLie}).
The direct limit properties
in the categories of topological groups
or smooth mani\-folds depend on the sequence $(G_n)_{n\in\N}$
(see Remark~\ref{remwdpr}).\\[4mm]
\textbf{3. \,Direct limits of unit groups of Banach algebras.}
Consider a\linebreak
sequence
$A_1\sub A_2\sub\cdots$ of unital Banach algebras
(such that all inclusion maps
are continuous homomorphisms of unital algebras)
and equip
$A:=\bigcup_{n\in \N}A_n$
with the locally convex direct limit topology.
We show that $A^\times=\dl\,A^\times_n$\vspace{-.7mm}
as a topological group
and (provided $A$ is Hausdorff)
also as a Lie group
(Proposition~\ref{Edamer}).
Non-unital
Banach algebras are discussed as well (Proposition~\ref{Edamer2}).
In the case where each inclusion map is an isometry,
the topological result has been
obtained earlier by Edamatsu
\cite[Theorem~1]{Eda}.\\[4mm]
\textbf{4. \,Lie groups of germs of analytic mappings.}
We also discuss the direct limit
properties of the groups $\Gamma(K,H)$
of germs of $H$-valued analytic functions
on open neighbourhoods of a non-empty compact set
$K\sub X$, where $H$ is a (real or complex)
Banach-Lie group and $X$ a (real or complex)
metrizable locally convex space.
Such groups are interesting
in this context because they
are prototypical examples of direct limits
$G=\bigcup_{n\in \N}G_n$
of a direct system $(G_n)_{n\in \N}$
which is not strict (nor $(L(G_n))_{n\in \N}$).
For $X$ and~$H$ finite-dimensional,
$\Gamma(K,H)$ is modelled
on a Silva space, whence
$\Gamma(K,H)$ has all desired direct
limit properties by general results
concerning Silva-Lie groups
prepared in Section~\ref{silauto}
(see Proposition~\ref{propsgemsilva}).
For infinite-dimensional~$X$ (or~$H$),
we can still show that $\Gamma(K,H)=\dl\, G_n$\vspace{-.9mm}
as a topological group and as a Lie group
(Corollary~\ref{gemnonsilva}).\\[4mm]
Finally, we
obtain results concerning
\emph{locally convex direct limits}.
Consider a Hausdorff locally
convex space $E=\bigcup_{n\in\N}E_n$ which is the
locally convex direct limit $\dl\,E_n$\vspace{-.6mm}
of Hausdorff locally convex spaces~$E_n$.
Then $E=\dl\,E_n$\vspace{-.4mm}
as a topological group
(see \cite[Proposition~3.1]{Hir}).
Our results imply that
$E=\dl\,E_n$\vspace{-.5mm}
also as a Lie group
(Example~\ref{exlcx}).
It is well known that the direct limit topology on~$E$ can
be properly finer than the locally convex
direct limit topology.
We provide concrete criteria
ensuring that
$E\not=\dl\,E_n$\vspace{-.9mm}
as a topological space and
smooth mani\-fold
(Lemma~\ref{lemKM}).
These criteria are also useful for the study of general
Lie groups because of Theorem~\ref{usefcor}
(already mentioned).
Similarly, weak direct products of
Lie groups (and topological groups)
are useful tools for the study of general
direct limits (see Section~\ref{secprset}).\\[2.7mm]
{\bf Further results.}
In Section~\ref{secconstr},
we construct a Lie group
structure on a union $G=\bigcup_{n\in \N}G_n$
of infinite-dimensional
Lie groups, under suitable hypotheses.
The construction can be used to
turn groups of germs of analytic
diffeomorphisms around a non-empty compact subset
$K$ of $\C^n$ (or $\R^n$) into analytic
Lie groups (Section~\ref{diffgerm}).
Section~\ref{covgp} broaches on
universal covering groups
of direct limit groups.
Section~\ref{secopen}
compiles open questions.\\[2.7mm]
{\bf Structure of the article.}
Besides the introduction
and the list of references,
the article is structured as follows:\\[2.7mm]
{\footnotesize
\S\,1 \,Preliminaries and notation\,\dotfill\,\pageref{secprelim}\\
\S\,2 \,Tools to identify direct limits of
Lie groups\,\dotfill\,\pageref{secreduce}\\
\S\,3 \,Tools to identify direct limits of topological spaces
and manifolds\,\dotfill\,\pageref{secruleout}\\
\S\,4 \,Example: Weak direct products of
Lie groups\,\dotfill\,\pageref{secwdprod}\\
\S\,5 \,Example: Diffeomorphism groups\,\dotfill\,\pageref{secdiff}\\
\S\,6 \,Proof of the Fragmentation Lemma for
diffeomorphism groups\,\dotfill\,\pageref{secfrag}\\
\S\,7 \,Example: Test function
groups\,\dotfill\,\pageref{sectesf}\\
\S\,8 \,Proof of the Fragmentation Lemma for test
function groups\,\dotfill\,\pageref{profrgtf}\\
\S\,9 \,Direct limit properties of Lie
groups modelled on Silva spaces or
$k_\omega$-spaces\,\dotfill\,\pageref{silauto}\\
\S\,10 \,Example: Lie groups of germs of Lie-group
valued mappings\,\dotfill\,\pageref{secgems}\\
\S\,11 \,Tools to identify direct limits of topological
groups\,\dotfill\,\pageref{secprset}\\
\S\,12 \,Example: Direct limits of unit groups of
Banach algebras\,\dotfill\,\pageref{secunit}\\
\S\,13 \,Example: Lie groups of germs beyond
the Silva case\,\dotfill\,\pageref{beysilva}\\
\S\,14\, Construction of Lie group structures on direct
limit groups\,\dotfill\,\pageref{secconstr}\\
\S\,15\, Example: Lie groups of germs of analytic
diffeomorphisms\,\dotfill\,\pageref{diffgerm}\\
\S\,16\, Covering groups of direct limit groups\,\dotfill\,\pageref{covgp}\\
\S\,17\, Open problems\,\dotfill\,\pageref{secopen}\\[.5mm]
Appendix:
\,Smooth regularity of direct sums\,\dotfill\,\pageref{appdsum}}
\section{Preliminaries and notation}\label{secprelim}
%
%\ma{secprelim}
%
%
%
We fix notation
and terminology concerning differential
calculus, direct limits and
properties of locally convex spaces
($C^r$-regularity and related concepts).
\subsection*{{\normalsize Calculus in locally convex spaces and Lie groups}}
We recall some basic definitions of
Keller's $C^r_c$-theory
and the theory of analytic
mappings
(see \cite{RES}, \cite{GaN}, \cite{Mic} and \cite{Mil}
for the proofs and further details).
Throughout the article, $\K\in \{\R,\C\}$.
\begin{defn}\label{defsmth}
Let $E$ and $F$ be Hausdorff locally convex topological
$\K$-vector spaces,
$U\sub E$ be open,
$f\colon U\to F$ be a map,
and $r\in \N_0\cup\{\infty\}$.
We say that $f$ is $C^r_\K$ (or simply $C^r$)
if it is continuous
and, for all $k\in \N$ such that $k\leq r$,
the iterated directional derivatives $d^kf(x,y_1,\ldots,y_k):=
D_{y_1}\cdots D_{y_k}f(x)$ exist
for all $x\in U$ and $y_1,\ldots,y_k\in E$,
and define a continuous map $d^kf\colon U\times E^k\to F$.
The $C^\infty_\R$-maps are also called {\em smooth}.
\end{defn}
Occasionally, we shall also encounter
analytic maps.
\begin{defn}
If $\K=\C$, then $f$ as before is called
\emph{complex analytic}
if it is continuous and given locally by a
pointwise convergent series of continuous
homogeneous polynomials (see \cite[Definition~5.6]{BaS}
for details).
If $\K=\R$,  we call $f$ \emph{real analytic}
if it extends to a complex analytic map between
open subsets of the complexifications of~$E$ and~$F$.
\end{defn}
It is well known that $f$ is $C^\infty_\C$ if and only if it
is complex analytic (see, e.g., \cite{GaN}).
Since compositions of $C^r_\K$-maps
(resp., $\K$-analytic maps)
are $C^r_\K$ (resp., $\K$-analytic),
$C^r_\K$-manifolds and $\K$-analytic
Lie groups (modelled on Hausdorff locally convex spaces)
can be defined as usual
(see \cite{RES}, \cite{GaN}; cf.\ \cite{Mil}).
As in \cite{RES} and \cite{GaN},
we shall not presume that the modelling spaces
are complete.\footnote{Some
readers may prefer
to work with categories of Lie groups
and manifold modelled on complete,
sequentially complete, or Mackey complete
locally convex spaces. This only causes
minor changes of our results.}
Unlike \cite{Mil},
we do not require that manifolds
are regular topological spaces.\\[3mm]
We warn the reader that topological
spaces and locally convex spaces
are not assumed Hausdorff in this article.
However,
the topological
spaces underlying manifolds and Lie groups
are assumed Hausdorff.
We remark that switching to the categories of Hausdorff
topological groups (or Hausdorff topological spaces)
would not affect the validity (or failure)
of direct limit properties of Lie groups
in such categories (as Lie groups
are Hausdorff anyway).\\[3mm]
{\bf General conventions.}
The word ``Lie group'' (without further
specification) refers to a smooth
Lie group modelled on a real locally convex space.
A $C^\infty_\K$-Lie
group is
a smooth real Lie group (if $\K=\R$),
resp., a complex analytic Lie group (if $\K=\C$).
A $C^\omega_\K$-Lie group
group means a $\K$-analytic Lie group.
Likewise for manifolds.
If $(E,\|.\|)$ is a normed
space, $x\in E$ and $r>0$, we\vspace{-.3mm}
set $B^E_r(x):=\{y\in E\colon \|y-x\|<r\}$
and $\wb{B}^E_r(x):=\{y\in E\colon \|y-x\|\leq r\}$.\\[3mm]
See
\cite{GaN} for the following useful
fact (or \cite[Lemma~10.1]{Ber}, if $r\in \N_0\cup\{\infty\}$).
\begin{la}\label{closd}
Let $f\colon U\to F$
be as in Definition~{\rm\ref{defsmth}},
and $r\in \N_0\cup\{\infty,\omega\}$.
If $f(U)\sub F_0$ for a closed
vector subspace $F_0\sub F$,
then $f$ is $C^r_\K$
if and only if its co-restriction
$f|^{F_0}\colon U\to F_0$ is $C^r_\K$.\vspace{-3mm}\Punkt
\end{la}
\subsection*{{\normalsize Conventions and basic facts
concerning direct limits}}
We recall some basic
definitions and facts concerning direct limits.\vspace{-.6mm}
\begin{numba}
A {\em direct system\/}
in a category $\A$ is a pair
$\cS=((X_i)_{i\in I},(\phi_{i,j})_{i\geq j})$,
where $(I,\leq)$ is a directed set,
each $X_i$ an object of~$\A$,
and each $\phi_{i,j}\colon X_j\to X_i$ a morphism
(``bonding map'')
such that $\phi_{i,i}=\id_{X_i}$
and $\phi_{i,j}\circ \phi_{j,k}=\phi_{i,k}$ if
$i\geq j\geq k$. A {\em cone over $\cS$\/}
is a pair $(X,(\phi_i)_{i\in I})$, where $X\in \ob \,\A$
and $\phi_i\colon X_i\to X$ is a morphism
for $i\in I$ such that $\phi_i\circ \phi_{i,j}=\phi_j$
if $i\geq j$. A cone $(X,(\phi_i)_{i\in I})$
is a {\em direct limit cone\/} over $\cS$
in the category~$\A$
if, for every cone $(Y,(\psi_i)_{i\in I})$ over~$\cS$,
there exists a unique morphism $\psi\colon X\to Y$
such that $\psi\circ \phi_i=\psi_i$ for each~$i$.
We then write $(X,(\phi_i)_{i\in I})=\dl\, \cS$.\vspace{-.8mm}
If the bonding maps and ``limit maps''
$\phi_i$ are understood,
we simply call $X$ the {\em direct limit\/} of $\cS$
and write $X=\dl\, X_i$.\vspace{-.8mm}
If also $\cT=((Y_i)_{i\in I},(\psi_{i,j})_{i\leq j})$
is a direct system over $I$ and
$(Y,(\psi_i)_{i\in I})$ a cone over $\cT$,
we call a family $(\eta_i)_{i\in I}$
of morphisms $\eta_i\colon X_i\to Y_i$
{\em compatible\/}
if $\eta_i\circ \phi_{i,j}=\psi_{i,j}\circ \eta_j$
for $i\geq j$.
Then $(Y,(\psi_i\circ \eta_i)_{i\in I})$
is a cone over $\cS$;
we write $\dl\,\eta_i:=\eta$\vspace{-1mm}
for the morphism $\eta\colon X\to Y$
such that $\eta\circ \phi_i=\psi_i\circ \eta_i$.\vspace{-.6mm}
\end{numba}
\begin{numba}
For all direct systems
$\cS=((X_i)_{i\in I},(\phi_{i,j})_{i\geq j})$
encountered in the article,
$I$ will contain a co-final
subsequence. It therefore suffices
to state all results for the case
where the directed set
is $(\N,\leq)$, i.e., for direct
sequences.\vspace{-.6mm}
\end{numba}
\begin{numba}\label{dirset}
\emph{Direct limits of sets, topological spaces, and groups}.
For basic facts concerning direct limits
of topological spaces and topological groups,
the reader is referred to
\cite{DIR}, \cite{Han}, \cite{Hir}
and \cite{TSH}
(where also many of the pitfalls and subtleties
of the topic are described).
In particular, we shall frequently
use that the set underlying
a direct limit of groups is the
corresponding direct limit in the category
of sets, and that the direct limit
of a direct system of topological spaces
in the category of (not necessarily
Hausdorff) topological spaces
is its direct limit in the category
of sets, equipped with the final topology with respect
to the limit maps.
A direct system $((X_i)_{i\in I}, (\phi_{i,j})_{i\geq j})$
of topological spaces is called
\emph{strict} if each bonding map $\phi_{i,j}\colon X_j\to X_i$
is a topological embedding.
\end{numba}
The following simple fact will be used:
%
%
%
%
%
%\ma{baseDL}
\begin{la}\label{baseDL}
Consider a direct sequence
$((X_n)_{n\in\N}, (i_{n,m})_{n\geq m})$
of topological spaces~$X_n$
and continuous maps $i_{n,m}\colon X_m\to X_n$,
with direct limit $(X,(i_n)_{n\in \N})$
in the category of topological spaces.
Let $U_n\sub X_n$ be open subsets
such that $i_{n,m}(U_m)\sub U_n$
whenever $m\leq n$, and $U:=\bigcup_{n\in \N}i_n(U_n)$.
Then $U$ is open in~$X$,
and $(U,(i_n|_{U_n}^U)_{n\in \N})\!=\!\dl\,((U_n)_{n\in \N},
(i_{n,m}|_{U_m}^{U_n})_{n\geq m})$\vspace{-.9mm}
in the category of topological spaces.
\end{la}
\begin{proof}
It is clear that $U=\dl\,U_n$\vspace{-.4mm}
as a set, and that the inclusion map
$\dl\,U_n\to\dl\,X_n$\vspace{-.4mm}
is continuous (being continuous
on each~$U_n$).
If $V\sub \dl\,U_n$\vspace{-.4mm}
is open (e.g., $V=U$), then $V_n:=(i_n|_{U_n})^{-1}(V)$
is open in $U_n$ and hence in~$X_n$,
for each $n\in \N$.
If $x\in i_n^{-1}(V)$, there exists
$k\geq n$ such that $i_n(x)=i_k(y)$
for some $y\in V_k$.
Hence, there is $\ell\geq k$ such that
$i_{\ell,n}(x)=i_{\ell,k}(y)\in V_\ell$.
We deduce that $i_n^{-1}(V)=\bigcup_{\ell\geq n}i_{\ell,n}^{-1}(V_\ell)$
is open in~$X_n$. Hence $V$ open in~$X$.
\end{proof}
All necessary background concerning direct
limits of locally convex spaces
can be found in \cite{BTV}, \cite{Jar},
\cite{Sch}, and \cite{Tre}.\\[3mm]
{\bf General conventions:}
If we write $G=\bigcup_{n\in \N}G_n$
for a topological group (resp., $C^\infty_\K$-Lie group) $G$,
we always presuppose that each $G_n$ is
a topological group (resp., $C^\infty_\K$-Lie group),
$G_n\sub G_{n+1}$ for each $n\in \N$,
and that all of the inclusion maps
$G_n\to G_{n+1}$ and $G_n\to G$
are continuous (resp., $C^\infty_\K$-)
homomorphisms.
Analogous conventions apply
to topological spaces $X=\bigcup_{n\in \N}X_n$
and manifolds $M=\bigcup_{n\in\N}M_n$.
\subsection*{{\normalsize Smooth regularity and related concepts}}
\begin{defn}
Given $r\in \N_0\cup\{\infty\}$,
a Hausdorff real locally convex space~$E$ is called
\emph{$C^r$-regular} if, for each $0$-neighbourhood
$U\sub E$, there exists a $C^r$-function
$f\colon E\to \R$ such that
$f(0)=1$ and $f|_{E\setminus U}=0$.
If $E$ is $C^\infty$-regular,
we also say that~$E$ is \emph{smoothly regular}.
\end{defn}
After composing a suitable smooth
self-map of~$\R$ with~$f$, we may assume that
$f|_V=1$ for some $0$-neighbourhood
$V\sub U$, $f(E)\sub [0,1]$,
and $\Supp(f)\sub U$.
%
%
%\ma{bassr}
\begin{rem}\label{bassr}
Note that every Hausdorff locally convex space is
$C^0$-regular, being a completely regular
topological space
(cf.\ \cite[Theorem~8.4]{HaR}).
It is easy to see that
every Hilbert space is smoothly regular.
Furthermore, vector subspaces and
(finite or infinite)
direct products of
$C^r$-regular locally convex spaces
are $C^r$-regular.
This implies that every nuclear locally
convex space is smoothly regular,
because it can be realized as
a vector subspace of a direct product
of Hilbert spaces
(cf.\ \cite[\S7.3, Corollary~2]{Sch},
also \cite{Pie} and \cite{Tre}).
\end{rem}
To prove that certain locally convex spaces
have pathological properties, at some point
we shall find it useful
to use ideas from Convenient Differential Calculus
(see \cite{FaK}, \cite{KaM}).
We recall that the final topology on a locally
convex space~$E$ with respect to the set
$C^\infty(\R,E)$ of smooth curves is called
the \emph{$c^\infty$-topology}.
We write $c^\infty(E)$ for~$E$, equipped
with the $c^\infty$-topology;
open subsets of $c^\infty(E)$
are called \emph{$c^\infty$-open}.
Given locally convex spaces $E$ and~$F$,
a map $f\colon U\to F$ on a $c^\infty$-open set
$U\sub E$ will be called a \emph{$c^\infty$-map}
if $f\circ \gamma\colon \R\to F$ is smooth
for each smooth curve $\gamma\colon \R\to E$
with $\gamma(\R)\sub U$.
\begin{defn}
We say that a locally convex space~$E$
is \emph{$c^\infty$-regular}
if, for each $0$-neighbourhood $U\sub c^\infty(E)$,
there exists a $c^\infty$-function $f\colon E\to\R$
such that $f(0)=1$ and $f|_{E\setminus U}=0$.
\end{defn}
\section{Tools to identify direct limits of Lie groups}\label{secreduce}
%
%\ma{secreduce}
%
%
Consider a Lie group $G=\bigcup_{n\in \N}G_n$
such that $G=\dl\,G_n$\vspace{-.3mm}
as a topological group.
Then automatically
also $G=\dl\,G_n$\vspace{-.3mm}
as a Lie group,
provided that $G$
admits a ``direct limit chart''
in a sense defined presently.
The existence of a direct limit chart
is a very natural requirement,
which ties together the
direct system of Lie groups
and its associated direct system of
locally convex topological Lie algebras.
The concept can be defined more generally
for direct systems of manifolds
modelled on locally convex spaces.
%
%
%\ma{DLchart}
\begin{defn}\label{DLchart}
Let $r\in \N_0 \cup\{\infty,\omega\}$
and $M$ be a $C^r_\K$-manifold
such that $M=\bigcup_{n\in \N}M_n$
for an ascending sequence
$M_1\leq M_2\leq\cdots$ of $C^r_\K$-manifolds,
such that all inclusion maps $i_{n,m}\colon M_m\to M_n$,
($n\geq m$), and $i_n\colon M_n\to M$
are~$C^r_\K$.
Let $E$ and $E_n$ be the modelling spaces of~$M$
and $M_n$, respectively.
A chart $\phi\colon U\to V\sub E$
of~$M$
is called a \emph{weak direct limit chart}
if (a) and (b) hold for some $n_0\in \N$:
\begin{itemize}
\item[(a)]
There exist continuous linear maps
$j_{n,m}\colon E_m\to E_n$ for $n\geq m\geq n_0$
and $j_n\colon E_n\to E$
such that $\cS:=((E_n)_{n\geq n_0},(j_{n,m})_{n\geq m\geq n_0})$
is a direct system of locally convex spaces,
$(E,(j_n)_{n\geq n_0})$ is a cone of locally convex spaces
over~$\cS$, and $E=\bigcup_{n\geq n_0}j_n(E_n)$.
\item[(b)]
There exist charts
$\phi_n\colon U_n\to V_n\sub E_n$
of~$M_n$ such that
$U_m\sub U_n$ and $j_{n,m}(V_m)\sub V_n$
if $n\geq m\geq n_0$,
$U=\bigcup_{n\geq n_0}U_n$,
$V=\bigcup_{n\geq n_0}j_n(V_n)$
and
%\ma{spllot}
\begin{equation}\label{spllot}
\phi\circ i_n|_{U_n}\;=\; j_n|_{V_n} \circ \phi_n\quad
\mbox{for each $\,n\geq n_0$.}
\end{equation}
\end{itemize}
If, furthermore, $(E,(j_n)_{n\geq n_0})=\dl\,\cS$\vspace{-.4mm}
as a locally convex space, then $\phi$ is called
a \emph{direct limit chart}.
If $\phi$ is a direct limit
chart and $U\cap M_n=U_n$ for each $n\geq n_0$,
we call $\phi$ a \emph{strict direct limit chart}.
We say that a Lie group~$G$
admits a direct limit chart if it has a direct
limit chart around~$1$ (and hence also
a direct limit chart around any $x\in G$).
\end{defn}
%
%
%
%\ma{remdlcha}
\begin{rem}\label{remdlcha}
With notation as in Definition~\ref{DLchart},
we have:
\begin{itemize}
\item[(a)]
Each of the linear maps $j_n\colon
E_m\to E_n$
(and hence also each $j_{n,m}$)
is injective because
$j_n|_{V_n}=\phi\circ i_n\circ \phi_n^{-1}$
by (\ref{spllot}), which is injective.
Identifying
$E_n$ with a subspace of~$E$
via $j_n$,
we may assume that $j_n$ (and each $j_{n,m}$)
simply is the inclusion map.
Then (\ref{spllot}) becomes $\phi|_{U_n}=\phi_n$,
and we have $\phi_n|_{U_m}=\phi_m$
if $n\geq m$.
\item[(b)]
If $M=\bigcup_{n\in\N}M_n$
admits a weak direct limit chart
$\phi\colon U\to V$ (as above)
around $x\in M$, we may assume
that $x\in M_{n_0}$
and $\phi(x)=0$ (after a translation).
We shall usually assume this
in the following.
\item[(c)]
Let $\phi$ be a weak direct limit chart
around $x\in M$. If $r\geq 1$,
we can identify $T_xM_n$ with $E_n$ and
$T_xM$ with $E$, using the
linear automorphisms $d\phi_n(x)$, resp., $d\phi(x)$.
Then
$j_{n,m}=T_x(i_{n,m})$ and $j_n=T_x(i_n)$
for all integers $n\geq m\geq n_0$.
\item[(d)]
If $\phi\colon U\to V$
is a (weak) direct limit chart around~$x$
and $W\sub U$ an open
neighbourhood of~$x$, then also
$\phi|_W\colon W \to \phi(W)$ is a (weak) direct limit chart,
because $W=\bigcup_{n\in \N}(W\cap U_n)$
with $W\cap U_n$ open in~$U_n$.
\end{itemize}
\end{rem}
\begin{examples}
We shall see later that countable weak direct products
of Lie groups,
groups of compactly supported
diffeomorphisms and
test function groups
admit a (strict) direct
limit chart
(see Remark~\ref{DLforsum},
Remark~\ref{strDLchadf} and Remark~\ref{strchatf},
respectively).
The Lie groups
of germs encountered in Section~\ref{secgems}
admit a direct limit chart
(albeit not a strict one).
\end{examples}
In the absence of
additional conditions (like
existence of a direct limit chart),
we cannot hope to establish direct
limit properties.
The following examples
illustrate some of the possible
pathologies.
\begin{example}
Let $G$ be the additive topological group
of the locally convex space $\R^{(\N)}$.
Give $G_n:=\R^n$ and $D:=\R^{(\N)}$
the discrete topology.
Then $G=\bigcup_{n\in \N}G_n$,
but the topologies on the subgroups
$G_n$ are just too fine compared to the topology on~$G$
to be of any use,
and the discontinuity of the
homomorphism $\id\colon G\to D$
(which is smooth on each $G_n$)
shows that $G\not=\dl\,G_n$\vspace{-.4mm}
as a Lie group,
topological group, smooth manifold,
and as a topological space.
\end{example}
The next example shows that the
existence of a weak direct
limit chart on a Lie group $G=\bigcup_{n\in \N}G_n$
is not enough
for the discussion of direct limit properties of~$G$.
For other purposes,
it suffices
(e.g., Proposition~\ref{complf} below).
\begin{example}
If we give $G:=C^\infty_c(\R,\R)=\bigcap_{k\in \N_0}C^k_c(\R,\R)$
the (unusually coarse\,!)
topology of the projective
limit of LB-spaces $\pl_{k\in \N_0}\,C^k_c(\R,\R)$,\vspace{-.4mm}
then $G_n:=C^\infty_{[-n,n]}(\R,\R)$ is a closed
vector subspace (and hence a Lie subgroup)
of~$G$, and $\id_G\colon G\to G$
is a weak direct limit chart
for $G=\bigcup_{n\in \N}G_n$.
Let $H:=C^\infty_c(\R,\R)$,
with the usual LF-topology.
The discontinuity of the homomorphism
$\id\colon G\to H$ (which is smooth on
each~$G_n$)
shows that $G\not=\dl\,G_n$\vspace{-.3mm}
as a Lie group, topological group,
topological space, and smooth manifold.
\end{example}
Recall a well-known fact:
If $f\colon G\to H$ is $C^\infty_\K$-homomorphism,
and $x\in G$, then
$\lambda^H_{f(x)}\circ f=f\circ \lambda_x^G$
holds with left translations
as indicated, whence
\begin{equation}\label{obvi}
T_x(f)\;=\; T_1(\lambda^H_{f(x)})\circ T_1(f)\circ T_x(\lambda^G_{x^{-1}})\,.
\end{equation}
%
%
%\ma{reduce}
\begin{thm}[{\bf Reduction\hspace*{-.6mm} to\hspace*{-.6mm}
Topological\hspace*{-.6mm} Groups}]\label{reduce}
Consider a $C^\infty_\K$-Lie\linebreak
group $G=\bigcup_{n\in \N}G_n$.
If $G=\dl\,G_n$\vspace{-.3mm}
as a topological group
and $G$ admits a direct limit
chart, then $G=\dl\,G_n$\vspace{-.3mm}
as a $C^\infty_\K$-Lie group.
\end{thm}
\begin{proof}
In view of Remark~\ref{remdlcha} (a) and (c),
we can identify $L(G_m)$
with a subalgebra of~$L(G_n)$
(if $n\geq m$) and~$L(G)$.
Let $f\colon G\to H$ be a homomorphism
to a $C^\infty_\K$-Lie group~$H$ such that $f_n:=f|_{G_n}$
is~$C^\infty_\K$, for each $n\in \N$.
Then $f$ is continuous, by the direct
limit property of~$G$ as a topological group.
Pick a chart $\psi \colon P\to Q\sub L(H)$
of~$H$ around~$1$
such that $\psi(1)=0$.
Let $\phi=\dl\,\phi_n\colon U\to V$\vspace{-.5mm}
be a direct limit chart of~$G$ around~$1$
such that $f(U)\sub P$,
where $\phi_n\colon U_n\to V_n$
and $\phi(0)=0$.
To see that $f$ is $C^1_\K$,
we pass to local coordinates:
We define
$h:=\psi\circ f|_U\circ \phi^{-1}\colon V\to Q$
and $h_n:=\psi\circ f|_U\circ \phi_n^{-1}\colon V_n\to Q$.
If $x\in V$, then $x\in V_{n_0}$ for some~$n_0$.
Given $n\geq n_0$ and $y\in L(G_n)$,
the limit
%
%\ma{dllc}
\begin{equation}\label{dllc}
dh(x,y) \; =\; {\textstyle\frac{d}{dt}\big|_{t=0}}h(x+ty)
\;=\; {\textstyle\frac{d}{dt}\big|_{t=0}}h_n(x +ty)
\; =\; dh_n(x,y)
\end{equation}
exists in $L(H)$.
We abbreviate
$\theta:=dh(0,\sbull)\colon L(G)\to L(H)$
and $h_n'(0):=dh_n(0,\sbull) \colon
L(G_n)\to L(H)$.
Since $\theta|_{L(G_n)}=h_n'(0)$
for each $n\geq n_0$
by (\ref{dllc}),
which is a continuous linear map,
and $L(G)=\dl\,L(G_n)$\vspace{-.5mm}
as a locally convex space, we deduce that
$\theta$ is continuous linear.
Since $f_n$ is a $C^\infty_\K$-homomorphism,
(\ref{obvi})
implies that
$dh_n(x,y)=d\lambda_{h_n(x)}^H(0,h_n'(0).d\lambda^{G_n}_{x^{-1}}(x,y))$
for all $x$ and~$y$ as before
(using the respective locally defined left translation
maps in local coordinates).
Hence
\[
dh(x,y)\;=\;d\lambda_{h(x)}^H (0,\theta. d\lambda^G_{x^{-1}}(x,y))\,,
\]
entailing that $dh\colon V\times L(G)\to L(H)$
is continuous.
Thus $h$ is $C^1_\K$ and hence also $f|_U$
is $C^1_\K$.
Since $f$ is a homomorphism,
it readily follows that~$f$ is~$C^1_\K$
(see \cite[Lemma~3.1]{PAD})
and hence~$C^\infty_\K$, by \cite[Lemma~2.1]{HOE}.
\end{proof}
We recall a simple fact
(cf.\ \cite[Prop.\,3.1]{Hir} and \cite[Exercise~14 to
Ch.\,II, \S4]{BTV}).
%
%\ma{sllma}
\begin{la}\label{sllma}
Let $\cS:=((E_n)_{n\in \N}, (f_{n,m})_{n\geq m})$
be a direct sequence of topological
$\K$-vector spaces.
\begin{itemize}
\item[\rm (a)]
Let $(E,(f_n)_{n\in \N})$ be the
direct limit of~$\cS$ in the category
of topological $\K$-vector spaces.
Then $E=\dl\,E_n$\vspace{-.4mm}
as a topological group.
If each $E_n$ is locally convex, then
$E=\dl\,E_n$\vspace{-.3mm}
also as a locally convex space.
\item[\rm (b)]
If each $E_n$ is Hausdorff,
let $F:=\dl\,E_n$\vspace{-.4mm}
in the category of Hausdorff
topological $\K$-vector spaces.
Then $F=\dl\,E_n$\vspace{-.4mm}
also as a Hausdorff topological group.
If each $E_n$ is locally convex,
then furthermore
$F=\dl\,E_n$\vspace{-.4mm}
as a Hausdorff locally convex space.
\end{itemize}
\end{la}
\begin{proof}
It is well known that the box topology on
$S:=\bigoplus_{n\in \N}E_n$ is a
vector topology which makes~$S$ the
direct sum in the category of
topological vector spaces
(cf.\ \cite[\S4.1, Proposition~4]{Jar})
and in the category
of topological abelian
groups (cf.\ also Lemma~\ref{prodmaps} below).
If each $E_n$ is locally convex, then also~$S$,
and the box topology coincides with the locally convex
direct sum topology
(cf.\
\cite[Exercise~14 to
Ch.\,II, \S4]{BTV}).
Let $i_n\colon E_n\to E$
be the canonical embedding.
Then the subgroup
$R:=\langle \bigcup_{n\geq m}\im\,((i_n\circ f_{n,m})-i_m) \rangle \leq S$
is a vector subspace of~$S$.
The direct limit~$E$ in each of the categories
described in~(a) can be realized as $S/R$.
If each $E_n$ is Hausdorff,
then the direct limit~$F$ in the categories
from~(b) can be realized as $S/\wb{R}$,
using the closure of~$R$.
\end{proof}
%
%
%\ma{exlcx}
\begin{example}\label{exlcx}
Let $E=\bigcup_{n\in \N}E_n$
be a Hausdorff locally convex space
over $\K\in \{\R,\C\}$
such that $E=\dl\,E_n$\vspace{-.3mm}
as a locally convex space.
Then the identity map $E\to E$
is a direct limit chart.
Since $E=\dl\,E_n$\vspace{-.5mm}
as a topological group (by Lemma~\ref{sllma}),
Theorem~\ref{reduce} shows that
$E=\dl\,E_n$\vspace{-.5mm} also as a $C^\infty_\K$-Lie group.
\end{example}
\begin{rem}
In more specialized categories of Lie
groups, direct limit properties can
be quite automatic, even in the absence
of a direct limit chart
(cf. also \cite[Proposition~4.24]{GCX}).
For example, consider a Lie group
$G=\bigcup_{n\in \N}G_n$ which is locally
exponential (i.e., $\exp_G\colon L(G)\to G$
exists and is a local $C^\infty$-diffeomorphism
at~$0$). If each $G_n$ has a
smooth exponential map and $L(G)=\dl\,L(G_n)$\vspace{-.4mm}
as a locally convex space,
then $G=\dl\,G_n$\vspace{-.4mm}
in the category of Lie groups
possessing a smooth exponential map.\footnote{At the time of writing,
no Lie group modelled on a complete locally convex space is
known which does not have a smooth exponential map.}
If, furthermore, each $G_n$ is locally
exponential, then $G=\dl\,G_n$\vspace{-.5mm}
also in the category of locally exponential
Lie groups (as an immediate consequence).
For the proof, consider a homomorphism
$f\colon G\to H$ to a Lie group~$H$
having a smooth exponential function,
such that $f_n:=f|_{G_n}$ is smooth
for each~$n$.
Let $\theta\colon L(G)\to L(H)$
be the unique continuous linear
map such that $\theta|_{L(G_n)}=L(f_n)$
for each~$n$.
Then
$f \circ \exp_G|_{L(G_n)}=
f_n \circ\exp_{G_n}=\exp_H\circ L(f_n)=
\exp_H\circ \,\theta|_{L(G_n)}$
for each~$n$ and thus
$f\circ \exp_G=\exp_H\circ \theta$,
entailing that $f$ is smooth on some
identity neighbourhood and hence smooth.
\end{rem}
The following application of
weak direct limit charts
is the main result of~\cite{HOM}.
Information on homotopy groups
(notably on $\pi_1(G)$
and $\pi_2(G)$)
is important for the extension
theory of infinite-dimensional
groups (see \cite{NeC}--\cite{NeN}).
%
%\ma{complf}
%
\begin{prop}\label{complf}
If a Lie group
$G=\bigcup_{n\in\N}G_n$
admits a weak direct
limit chart,
then its connected component of the identity
is $G_0=\bigcup_{n\in \N}(G_n)_0$.
Furthermore,
$\pi_k(G)=\dl\, \pi_k(G_n)$\vspace{-.9mm},
for each $k\in \N$.\Punkt
\end{prop}
Further applications of direct limit charts
can be found in~\cite{SMA}.
\section{Tools to identify direct limits of topological spaces
and manifolds}\label{secruleout}
%
%\ma{secruleout}
%
Given a Lie group $G=\bigcup_{n\in \N}G_n$,
it is natural to hope
that $G=\dl\,G_n$\vspace{-.3mm}
as a topological space if and only if
$L(G)=\dl\,L(G_n)$\vspace{-.3mm}
as a topological space.
In this section, we show that
this hope is justified
if~$G$ has a direct limit chart.
The analogous problem
for the category of smooth manifolds
is also addressed.
At the end of the section, we
consider a locally convex direct
limit $E=\dl\,E_n$\vspace{-.3mm}
and compile conditions ensuring that
$E\not=\dl\,E_n$\vspace{-.3mm}
as a topological space (resp.,
smooth manifold).
%
%
%\ma{retrlem}
\begin{la}[{\bf Localization Lemma}]\label{retrlem}
Let $M$ be a $C^r$-manifold
modelled on a real locally convex space~$E$,
where $r\in \N_0\cup\{\infty\}$,
such that
$M$ is a regular topological space
$($e.g., $M$ is a Lie group$)$.
If $r\geq 1$, we assume that $E$ is $C^r$-regular.
Let $P\sub M$ be open and $x\in P$.
Then there exists a
$C^r$-map $\rho \colon M\to P$ with the following properties:
\begin{itemize}
\item[\rm (a)]
$\rho(y)=y$ for all $y$
in an open neighbourhood $Q\sub P$ of~$x$;
\item[\rm (b)]
The closure of $\{y\in M\colon \rho(y)\not=x\}$
in~$M$ is a subset of~$P$.
\end{itemize}
Furthermore, the following can be achieved:
\begin{itemize}
\item[\rm (c)]
If $M=\bigcup_{n\in \N}M_n$ and $M$ admits a
direct limit chart
around~$x$, then for every $n\in \N$ and $y\in M_n$
there exists $m\geq n$ and a neighbourhood $W$
of~$y$ in~$M_n$ such that $\rho(W)\sub M_m$
and $\rho|_W\colon M_n\supseteq W\to M_m$ is~$C^r$.
\item[\rm (d)]
If $M=\bigcup_{n\in \N}M_n$ and $M$ admits a
strict direct limit chart
around~$x$, then it can be achieved
that $\rho(M_n)\sub M_n$
for all $n\geq n_0$ and that
$\rho|_{M_n}\colon M_n \to M_n$ is a $C^r$-map,
for a suitable $n_0\in \N$.
\end{itemize}
\end{la}
\begin{proof}
If $r=0$, then $E$ is $C^0$-regular
(see Remark~\ref{bassr}).\vspace{1mm}

(a) and (b): After shrinking~$P$ if necessary,
we may assume that
there exists a chart $\phi\colon P\to V\sub E$ of~$M$
such that $\phi(x)=0$
and $[0,1]V=V$.
Since~$M$ is a regular topological space,
we find a $0$-neighbourhood $B\sub V$
such that $A:=\phi^{-1}(B)$ is closed in~$M$.
Since $E$ is a $C^r$-regular,
there exists a $C^r$-function
$\beta\colon E\to \R$
such that $\Supp(\beta)\sub B$,
$\im(\beta)\sub [0,1]$,
and such that $\beta|_R=1$ for some
$0$-neighbourhood $R\sub V$.
Set $Q:=\phi^{-1}(R)$.
Then
\[
\psi\colon P \to V\,,\qquad \psi(x)=\beta(\phi(x))\cdot \phi(x)
\]
is a $C^r$-mapping such that $\psi|_Q=\phi|_Q$
and $\psi(x)=0$ for each $x\in P\setminus A$.
Extending $\psi$ by $0$,
we obtain a $C^r$-map $\tilde{\psi}\colon M\to V$.
Then $\rho :=\phi^{-1}\circ \tilde{\psi}$
is a $C^r$-map $M\to M$ such that $\rho(M)\sub P$,
$\rho|_Q=\id_Q$ and $\rho|_{M\setminus A}=x$.\vspace{1mm}

(c) Because we can always pass to a cofinal subsequence
$(M_{n+n_0})_{n\in \N}$,
we may assume that $x\in M_1$.
After shrinking~$P$ if necessary, we may assume
that $\phi\colon P\to V$ is a direct limit chart,
say $\phi=\bigcup_{n\in\N}\phi_n$,
where
$P=\bigcup_{n\in \N}P_n$ and
$V=\bigcup_{n\in \N}V_n$ for certain
compatible charts
$\phi_n\colon P_n\to V_n\sub E_n$
around~$x$, where $E_n$ is the modelling spaces of~$M_n$.
Let $y\in M$,
say $y\in M_n$.
If $y\not\in A$, then $W:=M_n\setminus A$
is an open neighbourhood of~$y$ in~$M_n$
such that $\rho|_W$ is constant
(with value~$x$) and hence a $C^r$-map
into $M_n$.
If $y\in P$,
let $z:=\phi(y)\in V$.
Then $z\in V_k$ for some $k\geq n$
and $\beta(z)z\in V_m$ for some $m\geq k$.
By continuity of~$\beta$, scalar multiplication
in $E_m$ and the inclusion map $E_k\to E_m$,
there exists an open neighbourhood
$Z$ of $z$ in $V_k$ such that
$\beta(v)v\in V_m$ for all
$v\in Z$.
Then $W:=M_n\cap \phi_k^{-1}(Z)$ is an open
neighbourhood of~$y$ in $M_n$ such that
$\rho|_W$ is a $C^r$-map into $U_m\sub M_m$.\vspace{1mm}

(d) If $M$ admits a strict direct
limit chart at~$x$, then we may assume
that $V\cap E_n=V_n$
for each~$n$. Given $y$ and~$n$
as in the proof of~(c),
we can now take $m:=k:=n$,
from which (d) follows.
\end{proof}
%
%
%
%\ma{mainneg}
\begin{la}\label{mainneg}
Let
$M=\bigcup_{n\in \N}M_n$\vspace{-.3mm}
be a $C^r$-manifold,
$r\in \N \cup\{\infty\}$,
which
is a regular topological space
and admits a direct limit chart around some
$x\in M$.
\begin{itemize}
\item[\rm (a)]
If $T_xM\not=\dl\,T_xM_n$\vspace{-.4mm}
as a topological space,
then also $M\not=\dl\,M_n$\vspace{-.5mm}
as a topological space.
\item[\rm (b)]
If
$T_xM\not=\dl\,T_xM_n$\vspace{-.4mm}
as a $C^r_\R$-manifold
and $M$ is modelled on a $C^r$-regular locally convex space,
then $M\not=\dl\,M_n$\vspace{-.8mm}
as a $C^r_\R$-manifold.
\end{itemize}
\end{la}
\begin{proof}
We may assume that $x\in M_1$.
To prove~(a), set $s:=0$;
for~(b), set $s:=r$. 
Let $\phi\colon P\to V\sub T_xM$ be a direct limit
chart around~$x$ and
$\rho\colon \! M\to P$
be a $C^r$-map as in Lemma~\ref{retrlem}\,(c).
Since $T_xM\not=\dl\,T_xM_n$\vspace{-.3mm}
as a topological space (resp., $C^r$-manifold),
there exists a map $h\colon T_xM\to X$
to a topological space (resp., $C^r$-manifold) $X$
that is not~$C^s$,
although $h|_{M_n}$ is~$C^s$ for each~$n$.
Hence, there is $z\in T_xM$
such that $h$ is not $C^s$ on any open neighbourhood
of~$z$. We may assume that $z\in T_xM_1$;
after replacing $h$ with $h(\sbull-z)$,
we may assume that $z=0$.
Then $f:=h \circ \phi\circ \rho$ is not $C^s$.
We claim that $f|_{M_n}$
is~$C^s$ for each $n\in \N$;
thus $M$ is not the direct limit
topological space (resp., $C^r$-manifold).
To prove the claim, let $A\sub P$
be a closed subset of~$M$ such that $\rho|_{M\setminus A}=x$,
and $y\in M_n$.
If $y\not\in A$, then $W:=M_n\setminus A$ is an open neighbourhood
of~$y$ in~$M_n$ such that $f(w)=h(0)$
for each $w\in W$;
thus $f|_W$ is~$C^s$.
If $y\in P$, let~$m$ and~$W$ be as in the proof
of Lemma~\ref{retrlem}\,(c).
Then $\rho(W)\sub U_n$ and $f|_W
=h|_{V_n}\circ\phi_n\circ \rho|_W^{U_n}$
is $C^s$, as claimed.
\end{proof}
Replacing the topological space~$X$ by a $C^0$-manifold in the
proof of~(a), we see that $M\not=\dl\,M_n$\vspace{-.5mm}
as a $C^0$-manifold if $T_xM\not=\dl\, T_xM_n$\vspace{-.3mm}
as a $C^0$-manifold.
%
%\ma{usefcor}
\begin{thm}[{\bf Reduction to the Lie Algebra Level}]\label{usefcor}
\hspace*{-.4mm}Let
$G\!=\!\bigcup_{n\in \N}G_n$\linebreak
be a real Lie group admitting a direct limit chart.
\begin{itemize}
\item[\rm (a)]
Then $L(G)=\dl\,L(G_n)$\vspace{-.3mm}
as a topological space
if and only if $G=\dl\,G_n$\vspace{-.3mm}
as a topological space.
\item[\rm (b)]
If $L(G)$ is $C^r$-regular
for $r\in \N_0\cup\{\infty\}$, then
$L(G)=\dl\,L(G_n)$\vspace{-.3mm}
as a $C^r$-manifold
if and only if
$G=\dl\,G_n$\vspace{-.3mm}
as a $C^r$-manifold.
\end{itemize}
\end{thm}
\begin{proof}
We set $r:=0$ in the situation of~(a).
If $L(G)\not=\dl\,L(G_n)$\vspace{-.4mm}
as a topological space (resp., as
a $C^r$-manifold, in
the situation of~(b)),
then $G\not=\dl\, G_n$\vspace{-.3mm}
as a topological space (resp., $C^r$-manifold),
by Lemma~\ref{mainneg}.\vspace{1mm}

Conversely, assume that $G\not=\dl\,G_n$\vspace{-.3mm}
as a topological space
(resp., as a $C^r$-manifold,
in the situation of\,(b)).
Then there exists a map
$f\colon G\to X$ to a topological space
(resp., $C^r$-manifold) $X$
which is not $C^r$,
although $f|_{G_n}$ is $C^r$ for each $n\in \N$.
There is $x\in G$ such that $f$ is not $C^r$ on any open
neighbourhood
of~$x$. After replacing $f$ with
$y\mto f(x^{-1}y)$, we may assume that $x=1$.
Let $\phi\colon U\to V\sub L(G)$ be
a direct limit chart of~$G$ around~$1$
such that $\phi(1)=0$,
with $\phi=\bigcup_{n\in \N}\phi_n$,
$U=\bigcup_{n\in \N}U_n$
and $V=\bigcup_{n\in \N}V_n$
for charts $\phi_n\colon U_n\to V_n\sub L(G_n)$
such that $\phi_{n+1}|_{U_n}=\phi_n$.
By Lemma~\ref{retrlem}, there exists
a $C^r$-map
$\rho\colon L(G)\to V$
such that $\rho|_Q=\id_Q$
for an open $0$-neighbourhood $Q\sub L(G)$
such that $\rho|_{V_n}$ locally
is a $C^r$-map into some~$V_m$, $m\geq n$.
Then $f\circ \phi^{-1}\circ \rho$
is not $C^r$,
although $(f\circ \phi^{-1}\circ \rho)|_{L(G_n)}$
is $C^r$ for each $n\in \N$.
\end{proof}
%
%
%\ma{recallyam}
\begin{rem}\label{recallyam}
Theorem~\ref{usefcor} complements
Yamasaki's Theorem \cite[Thm.\,4]{Yam}:\\[2.5mm]
\emph{Consider a group $G=\bigcup_{n\in \N}G_n$,
where each $G_n$ is a metrizable topological group and
each inclusion map $G_n\to G_{n+1}$ a topological embedding.
Assume that neither {\rm (a)} nor {\rm (b)} holds}:
\begin{itemize}
\item[\rm (a)]
\emph{For each $m\in \N$,
there exists $n\geq m$ and an identity neighbourhood
$U\sub G_m$ whose closure
in $G_n$ is compact};
\item[\rm (b)]
\emph{There exists $m\in \N$
such that $G_m$ is open in $G_n$
for each $n\geq m$.}
\end{itemize}
\emph{Then the direct limit topology does
not make $G$ a topological group.}
\end{rem}
Stimulated by Theorem~\ref{mainneg},
we turn to locally convex direct
limits and their direct limit properties
as topological spaces
and manifolds.
%
%
%\ma{lemKM}
\begin{la}\label{lemKM}
Let $E_1\sub E_2\sub\cdots$
be an ascending sequence
of Hausdorff locally convex spaces
which does not become stationary,
such that each $E_n$
is a vector subspace of $E_{n+1}$
and $E_{n+1}$ induces the given topology
on~$E_n$. Let $E:=\bigcup_{n\in\N}E_n$,
equipped with the locally convex direct
limit topology.
\begin{itemize}
\item[\rm (a)]
If each $E_n$ is
infinite-dimensional
and metrizable,
then $E\not=\dl\,E_n$\vspace{-.5mm}
as a topological space.
\item[\rm (b)]
If each $E_n$ is an infinite-dimensional
nuclear Fr\'{e}chet space,
then $E\not=\dl\,E_n$\vspace{-.5mm}
as a $C^r$-manifold,
for each $r\in \N_0\cup\{\infty\}$.
Furthermore, $E$ is smoothly regular.
\end{itemize}
\end{la}
\begin{proof}
(a) (Cf.\ Theorem 4.11\,(3)
and Proposition~4.26\,(ii) in~\cite{KaM}
if each $E_n$ is a Fr\'{e}chet space).
Let $E_1\subset E_2\subset \cdots$
be a strict direct sequence
of infinite-dimensional
metrizable topological vector spaces
which is strictly increasing.
Then $E_m$ is not open in $E_n$
for any integers $m<n$,
and the closure of a
$0$-neighbourhood $U$ of $E_m$ in~$E_n$
cannot be compact because then
$U$ would be pre-compact and
thus $\dim(E_m)<\infty$.
Now Yamasaki's Theorem
(see Remark~\ref{recallyam})
shows that the direct limit topology
does not make $E=\bigcup_{n\in \N}E_n$ a topological
group. The assertion follows.

(b) It is well known that the
$c^\infty$-topology on~$E$
coincides with the direct limit topology 
(cf.\ \cite[Theorem~4.11\,(3)]{KaM}).
By Part\,(a) just established
(or \cite[Proposition~4.26\,(ii)]{KaM}),
the latter is properly finer than
the locally convex direct limit topology.
Therefore, there exists a $0$-neighbourhood
$U\sub c^\infty(E)$
which is not a $0$-neighbourhood
of~$E$. Since $E$
is $c^\infty$-regular (see \cite[Theorem~16.10]{KaM}),
there exists a $c^\infty$-function
$f\colon E\to\R$ such that $f(0)=1$
and $f|_{E\setminus U}=0$.
Then $f|_{E_n}$ is a $c^\infty$-map
and hence smooth (as $E_n$
is metrizable), for each $n\in \N$.
However, $f$ is discontinuous
(and hence not $C^r$ for any
$r\in \N_0\cup\{\infty\}$).
In fact, if $f$ was continuous, then
$f^{-1}(\R^\times)\sub U$
would be a $0$-neighbourhood in~$E$
and hence also~$U$,
contradicting our choice of~$U$.
Like any countable locally convex direct
limit of nuclear spaces, $E$ is nuclear \cite[\S7.4, Corollary]{Sch}
and hence smoothly regular (see Remark~\ref{bassr}).
\end{proof}
We close this section with a variant of
Yamasaki's Theorem for Lie groups.
%
%
%\ma{yamlocexp}
\begin{prop}\label{yamlocexp}
Let $G=\bigcup_{n\in \N}G_n$ be a Lie group,
where $L(G_n)$ is metrizable
for each $n\in \N$.
Assume that
condition {\rm (i)} or {\rm (ii)} is satisfied:
\begin{itemize}
\item[\rm (i)]
$G$
has a direct limit chart, and the direct sequence
$(L(G_n))_{n\in \N}$ is strict;
\item[\rm (ii)]
For each $n\in \N$,
the Lie group $G_n$ has an exponential function which
is a local homeomorphism at~$0$,
and the direct sequence $(G_n)_{n\in\N}$ is strict.
\end{itemize}
If $G=\dl\,G_n$\vspace{-.3mm}
as a topological space,
then {\rm (a)} or {\rm (b)} holds:
\begin{itemize}
\item[\rm (a)]
$G_n$ is a finite-dimensional Lie group,
for each $n\in \N$; or:
\item[\rm (b)]
There exists $m\in \N$ such that
$G_m$ is open in $G_n$ for each $n\geq m$.
\end{itemize}
\end{prop}
\begin{proof}
If (i) holds but neither (a) nor (b),
then we find $n_0\in \N$ such that
$L(G_{n_0})\sub L(G_{n_0+1})\sub\cdots$
is a strict direct sequence
of infinite-dimensional Fr\'{e}chet
spaces which does not become stationary,
whence $G\not=\dl\,G_n$\vspace{-1.2mm}
as a topological space
by Theorem~\ref{usefcor}\,(a) and Lemma~\ref{lemKM}\,(a).\vspace{1mm}

If (ii) holds,
then $(L(G_n))_{n\in \N}$ is strict.
To see this, given $n\geq m$ define $j:=L(i_{n,m})$
and let $V\sub L(G_m)$
and $W\sub L(G_n)$ be open $0$-neighbourhoods
such that
$\phi:=\exp_{G_m}|_V$
is a homeomorphism onto an open identity neighbourhood
$\tilde{V}\sub G_m$,
$\psi:=\exp_{G_n}|_W$
is a homeomorphism onto an open identity neighbourhood
$\tilde{W}\sub G_n$,
and $j(V)\sub W$.
Since $i_{n,m}$ is a topological embedding,
there exists an open identity neighbourhood
$\tilde{X}\sub \tilde{W}$ such that
$\tilde{X}\cap G_m=\tilde{V}$;
define $X:=\psi^{-1}(\tilde{X})\sub W$.
Since $\exp_{G_n}\circ \,j|_V=i_{n,m}\circ \exp_{G_m}|_V$
is a topological embedding, also $j|_V$ is an
embedding, whence the continuous linear map~$j$
is injective. To see that~$j$ is a topological embedding,
let $(x_k)_{k\in \N}$ be a sequence
in $L(G_m)$ such that
$j(x_k)\to 0$ in $L(G_n)$.
After omitting finitely many terms,
we may assume that $j(x_k)\in X$
for each~$k$.
Then $\exp_{G_n}(j(x_k))=i_{n,m}(\exp_{G_m}(x_k))
\in \tilde{X}\cap G_m=\tilde{V}$
and $y_k:=\phi^{-1}(\exp_{G_n}(j(x_k)))\to 0$
in $L(G_m)$.
Since $y_k\in V$, we have that $j(y_k)\in W$.
Now $\psi(j(x_k))=\exp_{G_n}(j(x_k))=i_{n,m}(\exp_{G_m}(y_k))
=\exp_{G_n}(j(y_k))=\psi(j(y_k))$
by the definition of~$y_k$ and
naturality of~$\exp$.
Since $\psi$ is injective, we deduce that
$j(x_k)=j(y_k)$ and hence $x_k=y_k\to 0$.
Thus $j$ is an embedding.\\[2.5mm]
Now assume that (ii) holds but neither (a)
nor (b). After passing to a suitable subsequence
$(G_{n_k})_{k\in \N}$,
we may assume that each $G_m$ is infinite-dimensional
and $G_m$ is not open in
$G_n$ whenever $n>m$.
Then $L(G_m)$ (identified with $L(i_{n,m}).L(G_m)$)
is a proper vector subspace
of $L(G_n)$, because otherwise
$i_{n,m}(G_m)\supseteq i_{n,m}(\exp_{G_m}(L(G_m)))
=\exp_{G_n}(L(i_{n,m}).L(G_m))=\exp_{G_n}(L(G_n))$
would contain an open identity neighbourhood
and hence be an open subgroup
(which we just ruled out).
For any identity neighbourhood
$U\sub G_m$, we now show that
its closure $K:=\wb{U}$ in~$G_n$
cannot be compact.\\[2.5mm]
To see this, suppose to the contrary
that~$K$ was compact.
Let $V$, $W$, $\phi$ and $\psi$
be as earlier in the proof.
Since~$G_n$ is a regular topological
space, it has a closed identity neighbourhood
$A\sub G_n$ such that $A\sub \wt{W}$.
Then $Q:=\phi^{-1}(G_m\cap A\cap K)$
is a $0$-neighbourhood in~$L(G_m)$
such that $\psi(Q)$ has compact closure
$\wb{\psi(Q)}\sub A\cap K$,
whence $Q$ has compact closure
$\psi^{-1}(\wb{\psi(Q)})$ in $L(G_n)$.
The inclusion map $L(G_m)\to L(G_n)$
being an embedding, this entails
that $Q$ is precompact
and thus $\dim(L(G_m))<\infty$
(which is absurd).\\[2.5mm]
Now Yamasaki's Theorem
shows that
the direct limit topology
does not make $G$ a topological
group. It therefore differs
from given topology on~$G$.
\end{proof}
\section{Weak direct products of Lie groups}\label{secwdprod}
%
%\ma{secwdprod}
%
%
%
In this section, we recall the definition of weak direct
products of Lie groups and
analyze their direct limit properties.
\begin{numba}\label{boxtp}
If $(G_i)_{i\in I}$
is a family of topological groups,
we let
$\prod_{i\in I}^*G_i\leq \prod_{i\in I}G_i$
be the subgroup of all families
$(g_i)_{i\in I}$ such that
$g_i=1$ for all but finitely many~$i$.
A \emph{box} is a set of the form $\prod^*_{i\in I}U_i:=
(\prod_{i\in I}^*G_i)\cap \prod_{i\in I}U_i$,
where $U_i\sub G_i$ is open and
$1\in U_i$ for all but finitely many~$i$.
It is well known
that the set of boxes
is a basis for a topology
on
$\prod_{i\in I}^*G_i$
making it a topological group.
\end{numba}
\begin{numba}\label{chawpro}
If $(G_i)_{i\in I}$
%\ma{chawpro}
is a family of $C^r_\K$-Lie groups, where $r\in \{\infty,\omega\}$,
then
$\prod_{i\in I}^*G_i$ can be made a 
$C^r_\K$-Lie group, modelled on the locally
convex direct sum $\bigoplus_{i\in I}L(G_i)$
(see \cite{MEA}).
The Lie group structure is characterized by
the following property (cf.\ \cite[proof of Proposition~7.3]{MEA}):
{\em For each $i\in I$,
let $\phi_i\colon \tilde{U}_i\to \tilde{V}_i\sub L(G_i)$
be a chart of $G_i$ around~$1$
such that $\phi_i(1)=0$.
Let $U_i\sub \tilde{U}_i$
be an open, symmetric identity neighbourhood
such that $U_iU_i\sub \tilde{U}_i$,
$V_i:=\tilde{\phi}(U_i)$,
and $\phi_i:=\tilde{\phi}_i|_{U_i}^{V_i}$.
Then
\[
\kappa:=\bigoplus_{i\in I}\phi_i\colon
{\textstyle \prod_{i\in I}^*}U_i\to
\bigoplus_{i\in I}V_i\,,\quad
(x_i)_{i\in I}\mto (\phi_i(x_i))_{i\in I}
\]
is a chart for $\prod_{i\in I}^*G_i$.}
If $I$ is countable,
then $\bigoplus_{i\in I}L(G_i)$ carries
the box topology, entailing that the
topology underlying the Lie group
$\prod_{i\in I}^*G_i$
is the box topology from
{\bf\ref{boxtp}} (because the sets
of identity neighbourhoods coincide).
\end{numba}
\begin{rem}\label{DLforsum}
Assume that $I=\N$ in the preceding situation.
Then, as is clear,
$\prod_{n\in \N}^*G_n=\dl_{n\in \N}
\,\prod_{k=1}^nG_k$
as an abstract group. 
Since $\kappa$ restricts
to the chart $\prod_{k=1}^n\phi_k \colon
\prod_{k=1}^n U_k\to \prod_{k=1}^n V_k$
of $\prod_{k=1}^nG_k$,
we see that\vspace{.3mm}
$\prod_{n\in \N}^*G_n=\bigcup_{n\in \N}
\,\prod_{k=1}^nG_k$
admits a strict direct limit chart.
\end{rem}
%
%
%\ma{prodmaps}
\begin{la}\label{prodmaps}
Let $(G_n)_{n\in \N}$ be a sequence of
topological groups, $H$ be a topological
group and $f_n\colon G_n\to H$ be a
map which is continuous at~$1$
and such that $f_n(1)=1$.
Then the map
$f\colon \prod^*_{n\in \N}G_n\to H$
taking $x=(x_n)_{n\in \N}$ to
\[
f(x)\, :=\, f_1(x_1)f_2(x_2)\cdots f_N(x_N)\quad
\mbox{if $\, x_n=1$ for all $n>N$}
\]
is continuous at~$1$.
In particular, $\prod^*_{n\in \N}G_n
= \dl_{n\in\N}\, \prod_{k=1}^n G_k$\vspace{-1mm}
in the category of topological groups.
\end{la}
\begin{proof}
Given an identity neighbourhood $V_0\sub H$,
there is a sequence $(V_n)_{n\in \N}$
of identity neighbourhoods of~$H$
such that $V_nV_n\sub V_{n-1}$
for each $n\in \N$.
Then $V_1V_2\cdots V_n\sub V_0$,
for each $n\in \N$.
Since $f_n$ is continuous at~$1$,
the preimage $U_n:=f^{-1}_n(V_n)$
is an identity neighbourhood in~$G_n$.
Then $U:=\prod_{n\in \N}^*U_n$
is an identity neighbourhood
in $G:=\prod^*_{n\in \N}G_n$
such that $f(U)\sub V_0$.
In fact, if $x\in U$ and $x_n=1$
for all $n>N$, then
$f(x)=f_1(x_1)\cdots f_N(x_N)
\in V_1\cdots V_N\sub V_0$. Hence $f$
is continuous at~$1$.\\[3mm]
To prove the final assertion,
let $f\colon G \to H$
be a homomorphism to a topological group~$H$
such that $f_n:=f|_{G_n}\colon G_n\to H$
is continuous for each $n\in\N$.
Since $f$ is a homomorphism,
given $x\in G$ with $x_n=1$
for $n>N$ we have
\[
f(x)=f(x_1\cdots x_N)=f(x_1)\cdots f(x_N)
=f_1(x_1)\cdots f_N(x_N)\,.
\]
Thus $f$ is a mapping of the form
just discussed.
Therefore~$f$ is continuous
at~$1$ and hence continuous,
being a homomorphism.
\end{proof}
\begin{prop}\label{DLwprodLie}
Let $(G_n)_{n\in \N}$ be a sequence of
$C^\infty_\K$-Lie groups modelled on locally
convex spaces.
Then
%
%\ma{weakprcol}
\begin{equation}\label{weakprcol}
{\textstyle \prod^*_{n\in \N}} G_n
\;=\;
\dl_{n\in\N}\, \prod_{k=1}^n G_k
\end{equation}
holds as a topological group,
and as a $C^\infty_\K$-Lie group.
\end{prop}
\begin{proof}
By Lemma~\ref{prodmaps},
$G :=\prod^*_{n\in \N} G_n$ has the desired direct limit property
in the category of topological groups.
Since $G$ admits a direct limit chart
(see Remark~\ref{DLforsum}),
it also is the
desired direct limit in the category
of $C^\infty_\K$-Lie groups, by Theorem~\ref{reduce}.
\end{proof}
\begin{rem}\label{remwdpr}
There is no uniform
answer concerning the validity
of~(\ref{weakprcol})
in the categories of topological spaces
resp., smooth manifolds.
\begin{itemize}
\item[(a)]
If each $G_n$ is modelled
on a Silva space (or on a $k_\omega$-space,
as in {\bf\ref{defkom}}),
then (\ref{weakprcol})
holds in the category of topological
spaces and the category of
$C^r_\K$-manifolds, for each $r\in \N_0\cup\{\infty\}$
(see Proposition~\ref{gdDLSilva}\,(ii)
below). In particular, this is the case
if each~$G_n$ is finite-dimensional.
\item[(b)]
If each $G_n$ is modelled on an infinite-dimensional
Fr\'{e}chet space,
then $\prod^*_{n\in \N}G_n\not=\dl\,\prod_{k=1}^nG_k$\vspace{-.3mm}
as a topological space, by Theorem~\ref{usefcor}\,(a)
and Lemma~\ref{lemKM}\,(a).
If $L(G_n)$ is an infinite-dimensional
nuclear Fr\'{e}chet space
for each $n\in \N$, then
$\prod^*_{n\in \N}G_n\not=\dl\,\prod_{k=1}^nG_k$\vspace{-.3mm}
as a $C^r_\R$-manifold
for each $r\in \N_0\cup\{\infty\}$,
by Theorem~\ref{usefcor}\,(b)
and Lemma~\ref{lemKM}\,(b).
\item[(c)]
Suppose that, for some $r\in \N_0\cup\{\infty\}$,
each $L(G_n)$ is $C^r$-regular,
and $\bigoplus_{n\in \N}L(G_n)\not=\dl\,\prod_{k=1}^nL(G_k)$\vspace{-.4mm}
as a $C^r_\R$-manifold.
Then
$\bigoplus_{n\in \N}L(G_n)$ is $C^r$-regular
(by Proposition~\ref{bumponsum} in the appendix)
and therefore\linebreak
$\prod^*_{n\in \N}G_n\not=\dl\,\prod_{k=1}^nG_k$\vspace{-.3mm}
as a $C^r_\R$-manifold, by Theorem~\ref{usefcor}\,(b).
\end{itemize}
\end{rem}
\section{Diffeomorphism groups}\label{secdiff}
%
%\ma{secdiff}
%
%
In this section, we outline the proofs
of the results concerning diffeomorphism
groups described in the Introduction.
\begin{numba}\label{chdf}
Recall that the Lie group $\Diff_c(M)$ is modelled
on the space $C^\infty_c(M,TM)$
of compactly supported smooth vector fields.
To obtain a chart around $\id_M$,
one chooses a smooth Riemannian metric $g$
on~$M$, with associated exponential
map $\exp_g$. Then there is
an open $0$-neighbourhood
$V\sub C^\infty_c(M,TM)$
with the following properties:
For each $\gamma\in V$,
the composition $\psi(\gamma):=\exp_g\circ \, \gamma$
makes sense and is a $C^\infty$-diffeomorphism of~$M$;
$\psi\colon V\to\Diff_c(M)$ is
injective; and $\phi:=\psi^{-1}\colon U\to V$
(with $U:=\psi(V)$) is a chart
for $\Diff_c(M)$.
Furthermore, it can be achieved that,
for each compact subset $K\sub M$,
$\phi(U\cap \Diff_K(M))=V\cap
C^\infty_K(M,TM)$
and the restriction of $\phi$ to a map
$U\cap \Diff_K(M)\to V\cap
C^\infty_K(M,TM)$ is a chart for $\Diff_K(M)$
(see \cite{DIF}; cf.\ \cite{Mic}).
\end{numba}
\begin{rem}\label{strDLchadf}
Let $K_1\sub K_2\sub\cdots$
be an exhaustion of~$M$ by compact sets,
i.e., $M=\bigcup_{n\in \N}K_n$
and $K_n\sub K_{n+1}^\circ$
for each $n\in \N$,
where $K_{n+1}^\circ$ is the interior
of~$K_{n+1}$.
Then $(K_n)_{n\in \N}$
is a cofinal subsequence of the directed
set of all compact subsets of~$M$.
It is clear that any chart $\phi$ of $\Diff_c(M)$
of the form just described is
a strict direct limit chart
of $\Diff_c(M)=\bigcup_{n\in \N}\Diff_{K_n}(M)$.
\end{rem}
We begin with the negative results.
%
%\ma{negdiff}
\begin{prop}\label{negdiff}
Let $M$ be a
$\sigma$-compact, non-compact,
finite-dimensional
smooth manifold of positive dimension.
Then there exists
a discontinuous function
$f\colon \Diff_c(M)\to \R$
whose restriction to $\Diff_K(M)$
is smooth, for each compact set $K\sub M$.
Hence $\Diff_c(M)\not=\dl\, \Diff_K(M)$\vspace{-.5mm}
as a topological
space and as a $C^r$-manifold,
for each $r\in \N_0\cup\{\infty\}$.
\end{prop}
\begin{proof}
Let $(K_n)_{n\in \N}$ be an exhaustion
of~$M$ by compact sets.
Then
the space $C^\infty_{K_n}(M,TM)$
of smooth vector fields on~$M$ supported
in~$K_n$ is a nuclear Fr\'{e}chet
space (cf.\ \cite{Pie}, \cite{Sch}, \cite{Tre}),
and the direct sequence $(C^\infty_{K_n}(M,TM))_{n\in\N}$
is strict (each topology being induced
by $C^\infty(M,TM)$)
and does not become stationary.
Hence Lemma~\ref{lemKM}\,(b) and its proof
show that the locally convex direct limit
$C^\infty_c(M,TM)=\dl\,C^\infty_{K_n}(M,TM)$\vspace{-.3mm}
is smoothly regular,
and that there is a discontinuous map
$h\colon C^\infty_c(M,TM)\to\R$\vspace{-.3mm}
which is smooth on $C^\infty_{K_n}(M,TM)$
for each $n\in \N$.
After composing with a translation, we may assume
that $h$ is discontinuous at~$0$.
Since
$\Diff_c(M)=\bigcup_{n\in \N}\Diff_{K_n}(M)$
has a direct limit chart (Remark~\ref{strDLchadf}),
the proof of Lemma~\ref{mainneg}\,(b)
enables us to manufacture
a function $f\colon \Diff_c(M)\to \R$\vspace{-.3mm}
which is discontinuous at $\id_M$,
although its restriction to $\Diff_{K_n}(M)$
is smooth for each $n\in \N$.
\end{proof}
%
%
%
%\ma{posdiff}
\begin{prop}\label{posdiff}
Let $M$ be a non-compact,
$\sigma$-compact finite-dimensional
smooth manifold.
Then
$\Diff_c(M)=\dl\, \Diff_K(M)$\vspace{-.3mm}
in the category of topological
groups, and in the category of Lie groups.
\end{prop}
The proof hinges on the technique of
\emph{fragmentation}.
The idea of fragmentation
is to write a compactly supported
diffeomorphism as a composition of
diffeomorphisms supported in given sets
(cf.\ \cite[\S2.1]{Ban}
and the references therein;
cf.\ also \cite{HaT} for fragmentation in the convenient setting
of analysis).
The following lemma (proved in Section~\ref{secfrag})
establishes a link
between fragmentation and weak direct products;
it asserts that, close to $\id_M$,
diffeomorphisms can be decomposed smoothly
into pieces supported in some locally
finite cover of compact sets.
%
%
%\ma{fragment}
\begin{la}[{\bf Fragmentation Lemma for
Diffeomorphism Groups}]\label{fragment}$\;$\linebreak
For any finite-dimensional,
$\sigma$-compact
$C^\infty$-manifold~$M$,
the~following~holds:
\begin{itemize}
\item[\rm (a)]
There exists a locally finite cover $(K_n)_{n\in \N}$
of~$M$ by compact sets, an open identity neighbourhood
$\Omega\sub \Diff_c(M)$ and a smooth
map
\[
\Phi :\; \Omega\to {\textstyle \prod_{n\in \N}^*}\; \Diff_{K_n}(M)\,,
\quad
\gamma\mto \Phi(\gamma)=:(\gamma_n)_{n\in \N}
\]
such that $\Phi(1)=1$ and $\gamma=\gamma_1\circ \cdots\circ
\gamma_n$ for each $\gamma\in \Omega$
and each sufficiently large~$n$.
\item[\rm (b)]
If $(U_n)_{n\in \N}$ is a locally finite cover
of~$M$ by relatively compact, open sets, then
$(K_n)_{n\in \N}$ in {\rm (a)}
can be chosen such that $K_n\sub U_n$
for all $n\in \N$.\,\Punkt
\end{itemize}
\end{la}
{\bf Proof of Proposition~\ref{posdiff}.}
Because the hypotheses of Theorem~\ref{reduce}
are satisfied 
by $\Diff_c(M)$ and a cofinal subsequence
of its Lie subgroups
$\Diff_K(M)$ (see Remark~\ref{strDLchadf}),
we only need to show that
$\Diff_c(M)=\dl\, \Diff_K(M)$\vspace{-.5mm}
as a topological group.
To this end,
let $f\colon \Diff_c(M)\to H$ be a homomorphism
to a topological group~$H$
whose restriction $f_K\colon \Diff_K(M)\to H$
to $\Diff_K(M)$ is continuous,
for each compact subset $K\sub M$.
We have to show that $f$ is continuous.
Let $(K_n)_{n\in \N}$ and $\Phi$ be as in the Fragmentation
Lemma and consider the auxiliary function
$h\colon \prod_{n\in \N}^*\, \Diff_{K_n}(M)\to H$,
\[
h((\gamma_n)_{n\in \N})
\, :=\,
f_{K_1}(\gamma_1)\cdot\ldots \cdot  f_{K_N}(\gamma_N)\quad
\mbox{if $\,\gamma_n=1$ for all $n>N$,}
\]
which is continuous at $1$ by Lemma~\ref{prodmaps}.
Then
$f|_\Omega=h\circ \Phi$,
since $f$ is a homomorphism.
Thus also $f$ is continuous at~$1$
and hence continuous.\Punkt

\section{Smooth fragmentation of diffeomorphisms}\label{secfrag}
%
%\ma{secfrag}
%
%
%
In\hspace*{-.2mm} this\hspace*{-.2mm} section, we\hspace*{-.2mm}
prove\hspace*{-.2mm} the\hspace*{-.2mm} Fragmentation\hspace*{-.2mm}
Lemma\hspace*{-.2mm} for\hspace*{-.2mm} diffeomorphism\hspace*{-.2mm}
groups
(Lemma~\ref{fragment}). We start with a preparatory lemma.
\begin{la}\label{laaufsteig}
Let $r\in \N_0\cup\{\infty\}$
and $\pi\colon E\to M$ be a $C^r$-vector bundle
over a $\sigma$-compact finite-dimensional $C^r$-manifold~$M$.
Let $P \sub C^r_c(M,E)$
be a $0$-neighbourhood
in the space of compactly supported $C^r$-sections,
and
$(s_n)_{n\in \N}$ be a sequence in
$C^r(M,\R)$
such that, for each compact set $K\sub M$,
there exists $N\in \N$ such that
$s_n|_K=s_m|_K$ for all $n,m\geq N$.
Then there is a $0$-neighbourhood
$Q\sub P$ such that $s_n\cdot Q\sub P$
for all~$n\in \N$.
\end{la}
\begin{proof}
Let $(U_n)_{n\in \N}$
be a locally finite cover of~$M$ by relatively compact,
open sets. Then the linear map
\[
p\colon C^r_c(M,E)\to\bigoplus_{n\in \N}
C^r(U_n,E|_{U_n})\,,\quad
\gamma\mto (\gamma|_{U_n})_{n\in \N}
\]
is a topological embedding
(see \cite{SEC} or \cite[Proposition~F.19]{ZOO}).
Hence, there are
open $0$-neighbourhoods $V_n\sub C^r(U_n,E|_{U_n})$
such that $p^{-1}(\bigoplus_{n\in \N}V_n)\sub P$.
The hypothesis
entails that $F_n:=\{s_k|_{U_n}\colon k\in \N\}$
is a finite set, for each $n\in \N$.
Fix $n\in \N$.
Since $C^r(U_n,E|_{U_n})$
is a topological $C^r(U_n,\R)$-module
(see \cite{SEC} or \cite[Corollary~F.13]{ZOO}),
for each $s\in F_n$
the multiplication operator
$C^r(U_n,E|_{U_n})\to
C^r(U_n,E|_{U_n})$,
$\gamma\mto s \cdot \gamma$
is continuous. 
Hence, there exists an open $0$-neighbourhood 
$W_n\sub V_n$ such that
$s\cdot W_n\sub V_n$ for all
$s\in F_n$.
Then $Q:=p^{-1}(\bigoplus_{n\in \N}W_n)\sub P$
is an open $0$-neighbourhood
such that $s_n\cdot Q\sub P$ for each $n\in \N$.
The proof is complete.
\end{proof}
\begin{numba}
To prove Lemma~\ref{fragment},
let $(U_n)_{n\in \N}$ be a locally finite cover of~$M$
be relatively compact, open subsets $U_n\sub M$,
and $(h_n)_{n\in \N}$ be a smooth partition of unity
of~$M$ such that $K_n:=\Supp(h_n)\sub U_n$
for each $n\in \N$.
For each $n\in \N_0$, we set $s_n:=\sum_{i=1}^nh_i$.
We define $U_0:=K_0:=\emptyset$.
For each $n\in \N$, we set
$W_n:=U_n\cup U_{n-1}$
and choose $\xi_n\in C^\infty_c(W_n,\R)$
such that
$\xi_n|_{K_n\cup K_{n-1}}=1$.
We abbreviate $L_n:=\Supp(\xi_n)$.
\end{numba}
\begin{numba}
Let $\tilde{\phi}\colon \tilde{U}\to \tilde{V}\sub C^\infty_c(M,TM)$
be a chart of $\Diff_c(M)$ around $\id_M$
such that $\tilde{\phi}(\id_M)=0$ and
$\tilde{\phi}$ restricts to a chart
$\tilde{U}\cap \Diff_K(M)\to \tilde{V}\cap C^\infty_K(M,TM)$
of $\Diff_K(M)$, for each compact
subset $K\sub M$
(see {\bf\ref{chdf}}).
There exists an open, symmetric
identity neighbourhood $U\sub \tilde{U}$
such that $UU\sub \tilde{U}$.
Set $V:=\tilde{\phi}(U)$,
$\phi:=\tilde{\phi}|_U^V$,
and let $\phi_n\colon U\cap \Diff_{K_n}(M)\to V\cap C^\infty_{K_n}(M,TM)$
be the restriction
of $\phi$ to a chart of $\Diff_{K_n}(M)$.
By {\bf\ref{chawpro}}, the map
\begin{equation}\label{definkapp}
\kappa:=\bigoplus_{n\in \N}\,\phi_n \colon
{\textstyle\prod_{n\in \N}^*}(U\cap \Diff_{K_n}(M)) \to
\bigoplus_{n\in \N}(V\cap C^\infty_{K_n}(M,TM))
\end{equation}
sending
$(\eta_n)_{n\in \N}$ to $(\phi_n(\eta_n))_{n\in \N}$
is a chart of $\prod_{n\in\N}^* \Diff_{K_n}(M)$
around~$1$.
\end{numba}
\begin{numba}
Pick an open,
symmetric identity neighbourhood
$P\sub U$ such that $PP\sub U$.
Let $Q:=\phi(P)$.
By Lemma~\ref{laaufsteig},
there is an open $0$-neighbourhood
$S\sub Q$ such that $s_n \cdot S\sub Q$
for all $n\in \N$.
We set $R:=\phi^{-1}(S)$.
In local coordinates,
the group multiplication of $\Diff_c(M)$
corresponds to the smooth map
$\mu\colon Q\times Q\to V$,
$\mu(\sigma,\tau):=\sigma*\tau:=\phi(\phi^{-1}(\sigma)\circ \phi^{-1}(\tau))$.
The group inversion corresponds to the smooth map
$Q\to Q$, $\sigma\mto \sigma^{-1}:=\phi(\phi^{-1}(\sigma)^{-1})$.
For each compact set $K\sub M$,
the local multiplication restricts
to a smooth map $(Q\cap C^\infty_K(M,TM))^2\to V\cap C^\infty_K(M,TM)$,
and the local inversion to a smooth map
$Q\cap C^\infty_K(M,TM)\to Q\cap C^\infty_K(M,TM)$.
\end{numba}
\begin{numba}
Given $\gamma\in R$
and $n\in \N_0$,
we define $s_n \odot \gamma\in P$
via $s_n \odot \gamma:=
\phi^{-1}\big(s_n \cdot \phi(\gamma)\big)$.
For each $n\in \N$,
we let $\gamma_n:=(s_{n-1}\odot \gamma)^{-1} \circ (s_n\odot \gamma)
\in U$.
Since $(s_n\odot\gamma)(x)=(s_{n-1}\odot\gamma)(x)$
for all $x\in M\setminus K_n$,
we have $\gamma_n(x)=x$ for such~$x$ and thus
$\gamma_n\in C^\infty_{K_n}(M,TM)\cap U$.
Given $\gamma$, there is $N\in \N$
such that $U_n\cap \Supp(\gamma)=\emptyset$
for all $n>N$.
Then $s_n|_{\Supp(\gamma)}=1$ and thus
$s_n\odot \gamma=\gamma$ for all
$n\geq N$, whence $\gamma_1\circ \cdots \circ \gamma_n=s_n\odot\gamma
=\gamma$ for all $n\geq N$.
Thus
\[
\Phi\colon R \to
{\textstyle \prod_{n\in \N}^*}\Diff_{K_n}(M)\,,\quad
\Phi(\gamma):=(\gamma_n)_{n\in \N}
\]
has the desired properties, except for smoothness.
To complete the proof, we show that
$\Phi|_\Omega$ is smooth
for some open identity neighbourhood
$\Omega\sub R$.
\end{numba}
\begin{numba}
Since $(W_n)_{n\in \N}$
is a locally finite cover
of~$M$ by relatively compact, open
sets,
the map
\[
p\colon C^\infty_c(M,TM)\to\bigoplus_{n\in \N}
C^\infty(W_n,TW_n)\,,\quad
\gamma\mto (\gamma|_{W_n})_{n\in \N}
\]
is continuous linear
(and in fact an embedding
onto a closed vector subspace,
see \cite{SEC} or \cite[Proposition~F.19]{ZOO}).
Because $C^\infty(W_n,TW_n)$
is a topological $C^\infty(W_n,\R)$-module
(see \cite{SEC} or \cite[Corollary~F.13]{ZOO}),
the multiplication operators
$\mu_n\colon C^\infty(W_n,TW_n)\to
C^\infty_{L_n}(W_n,TW_n)$,
$\gamma\mto \xi_n\cdot s_n|_{W_n}\cdot \gamma$
and
$\lambda_n
\colon C^\infty(W_n,TW_n)\to C^\infty_{L_n}(W_n,TW_n)$,
$\gamma\mto \xi_n\cdot s_{n-1}|_{W_n} \cdot\gamma$
are continuous linear.
Then
\[
\lambda:=\bigoplus_{n\in \N}(\lambda_n,\mu_n)\colon
\bigoplus_{n\in \N}
C^\infty(W_n,TW_n)\to
\bigoplus_{n\in \N}
C^\infty_{L_n}(W_n,TW_n)
\times
C^\infty_{L_n}(W_n,TW_n)
\]
is continuous linear.
The restriction map
$\rho_n\colon
C^\infty_{L_n}(M,TM)\to C^\infty_{L_n}(W_n,TW_n)$
is an isomorphism of topological
vector spaces for each $n\in \N$
(see \cite{SEC} or \cite[Lemma~F15\,(b)]{ZOO}),
whence so is
\[
\rho:=\bigoplus_{n\in \N}(\rho_n^{-1}\times\rho_n^{-1})\,\colon\;
\bigoplus_{n\in \N}
C^\infty_{L_n}(W_n,TW_n)^2
\to \,
\bigoplus_{n\in \N}
C^\infty_{L_n}(M,TM)^2\,.
\]
\end{numba}
\begin{numba}
Then $Z:=\{\gamma\in V\colon (\rho\circ \lambda\circ p)(\gamma)
\in \bigoplus_{n\in \N} (Q\cap C^\infty_{L_n}(M,TM))^2\}$
is an open identity neighbourhood in~$V$.
For each $n\in \N$, the map
\[g_n\colon (Q\cap C^\infty_{L_n}(M,TM))^2\to
V\cap C^\infty_{L_n}(M,TM)\,, \quad g_n(\sigma,\tau):=
\sigma^{-1}*\tau
\]
is smooth.
Therefore, by \cite[Proposition~7.1]{MEA},
also the map
$g:=\bigoplus_{n\in \N}g_n\colon$
$\bigoplus_{n\in \N} (Q\cap C^\infty_{L_n}(M,TM))^2\to
\bigoplus_{n\in \N}(V\cap C^\infty_{L_n}(M,TM))$
is smooth. 
Then $g \circ \rho\circ\lambda\circ p|_Z$
is a smooth map with values in the closed
vector subspace
$\bigoplus_{n\in \N}C^\infty_{K_n}(M,TM)$ of
$\bigoplus_{n\in \N}C^\infty_{L_n}(M,TM)$,
and hence is $C^\infty$ also as a map into
this vector subspace (see Lemma~\ref{closd}).
We now consider
$g \circ \rho\circ\lambda\circ p|_Z$
as a $C^\infty$-map into the open $0$-neighbourhood
$\bigoplus_{n\in \N}(V\cap C^\infty_{K_n}(M,TM))
\sub \bigoplus_{n\in \N}C^\infty_{K_n}(M,TM)$.
Then $\Omega:=\phi^{-1}(Z)$ is an open identity neighbourhood
in $\Diff_c(M)$, and the formula
\[
\Phi|_\Omega\;=\;
\kappa^{-1}\circ g\circ \rho\circ \lambda\circ p|_Z\circ \phi|_\Omega
\]
shows that $\Phi$ is smooth.
This completes the proof of Lemma~\ref{fragment}.\Punkt
\end{numba}
\begin{rem}
Closer inspection
shows that $\Omega$ can be chosen
of the form $\Omega=\Diff_c(M)\cap \Omega_1$
for some open identity neighbourhood
$\Omega_1\sub \Diff^1_c(M)$ in the topological group
of compactly supported $C^1$-diffeomorphisms.
\end{rem}
The function $\Phi$ in the Fragmentation Lemma
is not unique, and in fact various
constructions give rise to
such functions.
The simple construction used in this section
has been adapted from~\cite{HaT}.
\section{Test function groups}\label{sectesf}
%\ma{sectesf}
%
%
In this section, we prove the results concerning
direct limit properties of test function
groups described in the Introduction.
More generally, we discuss $C^r_c(M,H)$
for $r\in \N_0\cup\{\infty\}$
and $H$ an arbitrary (not necessarily
finite-dimensional) smooth or $\K$-analytic Lie group.
\begin{numba}
We recall: If $r\in \N_0\cup\{\infty\}$
and $s\in \{\infty,\omega\}$,
$M$ is a $\sigma$-compact finite-dimensional
$C^r_\R$-manifold and $H$ a $C^s_\K$-Lie group,
then the group $C^r_c(M,H)$~of all
compactly supported \hspace*{-.3mm}$H$-valued $C^r_\R$-maps
on~$M$ is a $C^s_\K$-Lie group,
modelled on the locally convex direct limit
$C^r_c(M,L(H))=\dl\, C^r_K(M,L(H))$.\vspace{-1mm}
Its~Lie group structure
is characterized by the following property:
\emph{Let $\tilde{\phi}\colon \tilde{U}\to\tilde{V}\sub L(H)$
be a chart of~$H$ around~$1$
such that $\tilde{\phi}(1)=0$,
and $U\sub \tilde{U}$ be an open, symmetric identity neighbourhood
such that $UU\sub \tilde{U}$.
Set $V:=\tilde{\phi}(U)$
and $\phi:=\tilde{\phi}|_U^V$.
Then $C^r_c(M,U):=\{\gamma\in C^r_c(M,H)\colon
\gamma(M)\sub U\}$ is open in $C^r_c(M,H)$
and}
\[
C^r_c(M,\phi)\colon C^r_c(M,U)\to C^r_c(M,V)\,,\quad
\gamma\mto \phi\circ \gamma
\]
\emph{is chart for $C^r_c(M,H)$}
(cf.\ \cite[\S4.2]{GCX}).
\end{numba}
\begin{rem}\label{strchatf}
Let $K_1\sub K_2\sub\cdots$
be an exhaustion of~$M$ by compact
sets~$K_n$.
Since $C^r_{K_n}(M,\phi)\colon C^r_{K_n}(M,U)\to
C^r_{K_n}(M,V)\sub C^r_{K_n}(M,L(H))$ is a chart
of $C^r_{K_n}(M,H)$ for each $n\in \N$
(see \cite[\S3.2]{GCX}),
where $C^r_{K_n}(M,U)=C^r_c(M,U)\cap C^r_{K_n}(M,H)$,
we deduce that $C^r_c(M,\phi)$
is a strict direct limit chart for
$C^r_c(M,H)=\bigcup_{n\in \N}C^r_{K_n}(M,H)$.
\end{rem}
We begin our discussion of direct
limit properties with the negative results.
%
%
%\ma{negtestf}
\begin{prop}\label{negtestf}
Let $M$ be a
$\sigma$-compact, non-compact, finite-dimensional
smooth manifold of positive dimension
and $H$ be a non-discrete
Lie group whose locally convex modelling space
is smoothly regular 
$($\/for instance, a finite-dimensional
Lie group$)$.
Then there exists
a discontinuous map
\[
f\colon C^\infty_c(M,H)\to C^\infty_c(M,\R)
\]
whose restriction to $C^\infty_K(M,H)$
is~$C^\infty_\R$, for each compact set~$K\sub M$.
Hence $C^\infty_K(M,H)\not=\dl\,C^\infty_K(M,H)$\vspace{-.5mm}
as a topological space and as a $C^s_\R$-manifold,
for any $s\in \N_0\cup\{\infty\}$.
\end{prop}
The proof uses the following variant of Lemma~\ref{retrlem}.
%
%
%\ma{retract}
\begin{la}\label{retract}
Given $r \in \N_0\cup\{\infty\}$,
let $M$ be a
$\sigma$-compact finite-dimensional
$C^r$-manifold, $H$ be
a Lie group modelled on a locally convex space
which is smoothly regular,
and $P\sub H$ be an open identity neighbourhood.
Then there exists a smooth map
$\rho\colon C^r_c(M,H)\to C^r_c(M,H)$
with the following properties:
\begin{itemize}
\item[\rm (a)]
The image of $\rho$ is contained in $C^r_c(M,P)$;
\item[\rm (b)]
There exists an open identity neighbourhood
$Q\sub P$ such that
$\rho(\gamma)=\gamma$ for each $\gamma\in
C^r_c(M,Q)$; and
\item[\rm (c)]
$\rho$ restricts to a smooth map from
$C^r_K(M,H)$ to
$C^r_K(M,P)\sub C^r_K(M,H)$, for each compact subset $K\sub M$.
\end{itemize}
\end{la}
\begin{proof}
Using Lemma~\ref{retrlem}\,(e),
we find a smooth map $f\colon H\to P$
and an open identity neighbourhood
$Q\sub P$ such that $f|_Q=\id_Q$.
By \cite[Propositions~3.20 and~4.20]{GCX},
the map
\[
\rho:=C^r_c(M,f)\colon C^r_c(M,H)\to C^r_c(M,H)\,,\qquad
\gamma\mto f\circ \gamma
\]
is smooth and induces smooth self-maps
of $C^r_K(M,H)$ for each compact subset
$K\sub M$.
By construction,
it also has all other desired properties.
\end{proof}
{\bf Proof of Proposition~\ref{negtestf}.}
By \cite[Proposition~3.1]{DIS},
there exists a mapping $h\colon C^\infty(M,L(H))\to
C^\infty_c(M,\R)$ which is discontinuous at~$0$,
and such that $h|_{C^\infty_K(M,L(H))}$
is smooth for each compact subset $K\sub M$.
Choose a chart $\phi\colon U\to V\sub L(H)$
of~$H$ around~$1$
such that $\phi(1)=0$
and such that $\Psi:=C^\infty_c(M,\phi)\colon C^\infty_c(M,U)\to
C^\infty_c(M,V)$ is a chart of $C^\infty_c(M,H)$
and restricts to a chart of $C^\infty_K(M,H)$ with domain
$C^\infty_K(M,U)$, for each compact subset $K\sub M$.
Let $P\sub U$ and
$\rho\colon C^\infty_c(M,H)\to C^\infty_c(M,P)$
be as in Lemma~\ref{retract}.
Then the map
\[
f:=h\circ \Psi\circ \rho
\colon C^\infty_c(M,H)\to C^\infty_c(M,\R)
\]
is discontinuous at~$1$.
Since $\rho$ restricts to a smooth map
from $C^\infty_K(M,H)$
to $C^\infty_K(M,U)$
and $\Psi$ to a smooth map from $C^\infty_K(M,U)$
to $C^\infty_K(M,V)\sub C^\infty_K(M,L(H))$,
it follows that
$h$ is smooth on each $C^\infty_K(M,H)$.\vspace{3mm}\Punkt

\noindent
\begin{rem}
Given $r\in \N_0\cup\{\infty\}$,
let $M$ be a $\sigma$-compact, non-compact,
finite-dimensional
$C^r$-manifold of positive dimension,
and $H$ be a
Lie group modelled on a metrizable locally convex space
$\not=\{0\}$.
Then $C^r_K(M,H)$ is metrizable (cf.\
\cite[Proposition~4.19 (c) and (d)]{ZOO}).
Using Proposition~\ref{yamlocexp}\,(i),
we deduce that
$C^r_c(M,H)\not=\dl\,C^r_K(M,H)$\vspace{-.3mm}
as a topological space.
\end{rem}
Next, we establish the positive results.
%
%
%
%\ma{postestf}
\begin{prop}\label{postestf}
Let $M$ be a $\sigma$-compact,
finite-dimensional $C^r_\R$-manifold,
where $r\in \N_0\cup\{\infty\}$,
and~$H$ be a $C^\infty_\K$-Lie group modelled
on a locally convex space.
Then $C^r_c(M,H)=\dl\, C^r_K(M,H)$\vspace{-1mm}
as a topological group
and as a $C^\infty_\K$-Lie group.
\end{prop}
The next lemma (proved in Section~\ref{profrgtf})
helps us to prove Proposition~\ref{postestf}.
%
%
%\ma{fragtestf}
%
\begin{la}[Fragmentation Lemma for Test Function Groups]\label{fragtestf}
\hspace*{-.51mm}For\linebreak
$r\in\N_0\cup\{\infty\}$
and $s\in \{\infty,\omega\}$,
let $H$ be a $C^s_\K$-Lie group
modelled on a locally convex space
and $M$ be a $\sigma$-compact,
finite-dimensional $C^r_\R$-manifold.
\begin{itemize}
\item[\rm (a)]
Then
there exists a locally finite cover $(K_n)_{n\in \N}$
of~$M$ by compact sets,
an open identity neighbourhood
$\Omega\sub C_c^r(M,H)$ and a
$C^s_\K$-map
\[
\Phi :\; \Omega\to {\textstyle \prod_{n\in \N}^*}\;
C_{K_n}^r(M,H)\,,
\quad\gamma\mto \Phi(\gamma)=:(\gamma_n)_{n\in \N}
\]
such that $\Phi(1)=1$ and $\gamma=\gamma_1 \gamma_2 \,\cdots\,
\gamma_n$ for each
$\gamma\in \Omega$ and each sufficiently large~$n$.
\item[\rm (b)]
If $(U_n)_{n\in \N}$ is a locally finite
cover
of~$M$ by relatively compact, open sets,
then it can be achieved
in {\rm (a)}
that $K_n\sub U_n$ for each $n\in \N$.\Punkt
\end{itemize}
\end{la}
\begin{rem}
Considering a finite-dimensional
Lie group as a real analytic Lie group,
Lemma~\ref{fragtestf}
provides a \emph{real analytic}
fragmentation map in~this~case.
\end{rem}
{\bf Proof of Proposition~\ref{postestf}.}
By Remark~\ref{strchatf},
the hypotheses of Proposition~\ref{reduce}
are satisfied by $C_c^r(M,H)$ and any cofinal subsequence
of its Lie subgroups
$C^r_K(M,H)$.
Therefore,
we only need to show that
$C^r_c(M,H)=\dl\, C^r_K(M,H)$\vspace{-.3mm}
as a topological group.
Using Lemma~\ref{fragtestf}
instead of Lemma~\ref{fragment},
we can show this exactly
as in the proof of Proposition~\ref{posdiff}.\,\Punkt
\begin{rem}
While the case of a complex Lie group~$H$ is included
in Proposition~\ref{postestf},
we had to exclude it from Proposition~\ref{negtestf},
and the direct limit properties
of $C^r_c(M,H)$
in the category of complex manifolds
remain elusive
(because localization arguments
do not work in the complex case).
The following example makes it clear
that the direct limit
property can fail in some cases
(but there is no argument for the general case).
We consider the map
\[
f\colon C^\infty_c(M,\C)\to C^\infty_c(M\times M,\C)\, ,\quad
\gamma\mto \gamma\tensor\gamma
\]
with $(\gamma\tensor\gamma)(x,y):=\gamma(x)\gamma(y)$,
which is a homogeneous polynomial
of\linebreak degree~$2$.
It is clear from
\cite[Proposition~7.1]{ALG} and \cite[Lemma~3.7]{GCX}
that the restriction of $f$ to $C^\infty_K(M,\C)$
is a continuous homogeneous
polynomial of degree~$2$ and hence complex analytic.
However, $f$ is discontinuous
because the symmetric bilinear map
$C^\infty_c(M,\C)^2\to C^\infty_c(M\times M,\C)$
associated to~$f$ via polarization
is discontinuous (cf.\ \cite[Theorem~2.4]{TSH}).
\end{rem}
\begin{rem}
If $H$ is a real analytic Lie group,
then also $C^r_c(M,H)$ and $C^r_K(M,H)$
are real analytic Lie groups,
and one may ask whether
$C^r_c(M,H)=\dl\,C^r_K(M,H)$\vspace{-.3mm}
as a $C^\omega_\R$-Lie group
or as a $C^\omega_\R$-manifold.
However, since real analyticity is an even more
delicate property than complex analyticity,
one cannot expect
results except in special situations.
We therefore refrain from any attempt
in this direction, and merely
remind the reader of a
notorious pathology:
Already on $\R^{(\N)}=C^r_c(\N,\R)$,
non-analytic real-valued functions 
exist
which are $C^\omega_\R$ on $\R^n=C^r_{\{1,\ldots, n\}}(\N,\R)$
for each $n\in \N$ \cite[Example~10.8]{KaM}.
\end{rem}
\section{Smooth fragmentation of test functions}\label{profrgtf}
%
%\ma{profrgtf}
%
This section is devoted to the proof of
the Fragmentation Lemma for test function groups
(Lemma~\ref{fragtestf}).
We proceed in steps.
\begin{numba}
Let $(U_n)_{n\in \N}$ be a locally finite cover of~$M$
be relatively compact, open subsets $U_n\sub M$,
and $(h_n)_{n\in \N}$ be a $C^r$-partition of unity
of~$M$ such that $K_n:=\Supp(h_n)\sub U_n$
for each $n\in \N$.
For each $n\in \N_0$, we set $s_n:=\sum_{i=1}^nh_i$.
We define $U_0:=K_0:=\emptyset$.
For each $n\in \N$, we set
$W_n:=U_n\cup U_{n-1}$
and choose $\xi_n\in C^r_c(W_n,\R)$
such that $\xi_n(W_n)\sub [0,1]$
and $\xi_n|_{K_n\cup K_{n-1}}=1$.
We abbreviate $L_n:=\Supp(\xi_n)$.
\end{numba}
\begin{numba}
Pick a chart $\phi''\colon  U'' \to  V''  \sub L(H)$
of~$H$ around~$1$ such that
\mbox{$\phi''(1)=0$.}
Let $U'\sub U$ be an open, symmetric identity neighbourhood
with \mbox{$U'U'\sub U''$.}
Set $V'\!:=\phi''(U')$
and $\phi'\!:=\phi''|_{U'}^{V'}$.
Then
$C^r_c(M,\phi')\colon
\!C^r_c(M,U') \to C^r_c(M,V')\sub C^r_c(M,L(H))$
is a chart of $C^r_c(M,H)$
and
$C^r_K(M,\phi')\colon
\!C^r_K(M,U')\! \to C^r_K(M,V')$
$\sub C^r_K(M,L(H))$
is a chart of $C^r_K(M,H)$,
for each compact subset $K\sub M$.
Let $U\sub U'$ be an open, symmetric identity neighbourhood
with \mbox{$UU\sub U'$,} and set
$V:=\phi'(U)$ and $\phi:=\phi'|_U^V$.
Then the map
\begin{equation}\label{definkapp2}
\kappa:=\bigoplus_{n\in \N}\,C^r_{K_n}(M,\phi)\colon
{\textstyle\prod_{n\in \N}^*}C^r_{K_n}(M,U) \to
\bigoplus_{n\in \N}C^r_{K_n}(M,V)\,,
\end{equation}
$(\eta_n)_{n\in \N}\mto (\phi\circ \eta_n)_{n\in \N}$
is a chart of $\prod_{n\in\N}^* C^r_{K_n}(M,H)$
around~$1$
(see {\bf\ref{chawpro}}).
\end{numba}
\begin{numba}
There exists an open,
symmetric identity neighbourhood
$P\sub U$ such that $PP\sub U$.
We set $Q:=\phi(P)$,
let $S\sub Q$ be
an open $0$-neighbourhood
such that $[0,1]S=S$,
and define $R:=\phi^{-1}(S)$.
In local coordinates,
the group multiplication of $C^r_c(M,H)$
corresponds to the map
$\mu\colon Q\times Q\to V$,
$\mu(\sigma,\tau):=\sigma*\tau:=\phi(\phi^{-1}(\sigma)\phi^{-1}(\tau))$.
We write $\tau^{-1}:=\phi(\phi^{-1}(\tau)^{-1})$
for $\tau\in Q$.
\end{numba}
\begin{numba}
Given $\gamma\in C^r_c(M,R)$
and $h\in C^r_c(M,\R)$ such that $h(M)\sub [0,1]$,
we define $h\odot \gamma\in C^r_c(M,R)$
via $(h \odot \gamma)(x):=
\phi^{-1}\big(h(x)\cdot \phi(\gamma(x))\big)$.
For each $n\in \N$,
we let $\gamma_n:=(s_{n-1}\odot \gamma)^{-1} (s_n\odot \gamma)
\in C^r_c(M,U)$.
Since $(s_n\odot \gamma)(x)=(s_{n-1}\odot \gamma)(x)$
for all $x\in M\setminus K_n$,
we have $\gamma_n\in C^r_{K_n}(M,U)$.
Given $\gamma$, there is $N\in \N$
such that $U_n\cap \Supp(\gamma)=\emptyset$
for all $n>N$.
Then $s_n\odot \gamma=\gamma$ for all
$n\geq N$, whence $\gamma_1\cdots \gamma_n=\gamma$
for all $n\geq N$.
Thus
\[
\Phi\colon C^r_c(M,R)\to
{\textstyle \prod_{n\in \N}^*}C^r_{K_n}(M,H)\,,\quad
\Phi(\gamma):=(\gamma_n)_{n\in \N}
\]
will have the desired properties,
if we can show that this map is~$C^s_\K$.
\end{numba}
\begin{numba}
Since $(W_n)_{n\in \N}$
is a locally finite cover
of~$M$ by relatively compact, open
sets,
the map
\[
p\colon C^r_c(M,L(H))\to\bigoplus_{n\in \N}
C^r(W_n,L(H))\,,\quad
\gamma\mto (\gamma|_{W_n})_{n\in \N}
\]
is continuous linear
(and in fact an embedding
onto a closed vector subspace,
see \cite{SEC} or \cite[Proposition~8.13]{ZOO}).
Because $C^r(W_n,L(H))$
is a topological $C^r(W_n,\R)$-module
(see \cite{SEC} or \cite[Proposition~9.1\,(b)]{ZOO}),
the multiplication operators
$\mu_n\colon C^r(W_n,L(H))\to C^r_{L_n}(W_n,L(H))$,
$\gamma\mto \xi_n\cdot s_n|_{W_n} \cdot \gamma$
and
$\lambda_n
\colon C^r(W_n, L(H))\to C^r_{L_n}(W_n,L(H))$,
$\gamma\mto \xi_n\cdot s_{n-1}|_{W_n}\cdot\gamma$
are continuous linear.
Identifying
$C^r_{L_n} (W_n,L(H))^2$
with
$C^r_{L_n}(W_n, L(H)\times L(H))$
in the natural way (cf.\ \cite[Lemma~3.4]{GCX}),
we can consider $(\lambda_n,\mu_n)$
as a continuous linear map into
$C^r_{L_n}(W_n,L(H)^2)$.
Then
\[
\lambda:=\bigoplus_{n\in \N}(\lambda_n,\mu_n)\colon
\bigoplus_{n\in \N}
C^r(W_n,L(H))\to
\bigoplus_{n\in \N}
C^r_{L_n}(W_n,L(H)^2)
\]
is continuous linear.
Now $\lambda\circ p$
is a continuous linear map
which restricts
to a $C^s_\K$-map $f\colon C^r_c(M, S)\to
\bigoplus_{n\in \N}C^r_{L_n}(W_n,S\times S)$.
We define a mapping
$g_n\colon C^r_{L_n}(W_n,S\times S)\to
C^r_{L_n}(W_n,V)$
via
\[
g_n(\tau,\sigma):=\tau^{-1}*\sigma\,,
\]
using local inversion,
resp., the local multiplication~$*$
pointwise.
Then each $g_n$ is~$C^s_\K$,
as a consequence of \cite[Corollaries~3.11 and 3.12]{GCX},
and thus also
\[
g:=\bigoplus_{n\in \N}g_n\colon \bigoplus_{n\in \N}
C^r_{L_n}(W_n,S)\to
\bigoplus_{n\in \N}C^r_{L_n}(W_n,V)
\]
is $C^s_\K$ by \cite[Proposition~7.1 and Corollary~7.2]{MEA}.
Note that $g\circ \lambda\circ p$
has image in the closed vector subspace
$\bigoplus_{n\in \N}C^r_{K_n}(W_n,L(H))$
of $\bigoplus_{n\in \N}C^r_{L_n}(W_n,L(H))$;
we may therefore consider
$g\circ f$
as a $C^s_\K$-map into
$\bigoplus_{n\in \N}C^r_{K_n}(W_n,V)$
now (by Lemma~\ref{closd}).
For each $n\in \N$,
the map
$\rho_n\colon C^r_{K_n}(M,L(H))\to C^r_{K_n}(W_n,L(H))$,
$\gamma\mto \gamma|_{W_n}$
is an isomorphism of locally convex spaces
(see \cite[Lemma~4.24]{ZOO} or \cite{SEC}),
whence also
\[
\psi:= \bigoplus_{n\in \N} \rho_n^{-1}\colon
\bigoplus_{n\in \N}C^r_{K_n}(W_n,L(H))
\to \bigoplus_{n\in \N}C^r_{K_n}(M,L(H))
\]
is an isomorphism of locally convex spaces.
Consequently, the composition
\[
\Phi\, =\, \kappa^{-1}\circ\psi\circ
g \circ f\circ C^r_c(M,\phi|_R^S)
\]
is a $C^s_\K$-map
from $C^r_c(M,R)$ to $\prod_{n\in \N}^*C^r_{K_n}(M,H)$,
where $\kappa$ is as in~(\ref{definkapp2}).
This completes the proof of Lemma~\ref{fragtestf}.\Punkt
\end{numba}
\begin{rem}
The proof shows that $(K_n)_{n\in \N}$
can be chosen independently of $r$ and~$H$,
and that $\Omega$ can be chosen
of the form $C^r_c(M,U)$.
\end{rem}
\section{Direct limit properties of Silva-Lie groups
and Lie groups modelled on
{\boldmath $k_\omega$}-spaces}\label{silauto}
%
%
%\ma{silauto}
%
%
%
We describe conditions ensuring
that a Lie group $G=\bigcup_{n\in \N} G_n$
carries the direct limit topology
and is the direct limit in all categories
of interest.
In particular, the result applies
to many typical examples of Lie groups
modelled on Silva spaces.
We also obtain information on certain
Lie groups modelled on $k_\omega$-spaces.
%
%
%\ma{recsilva}
\begin{numba}\label{recsilva}
Recall that a locally convex space~$E$
is called a \emph{Silva space}
(or ($LS$)-space)
if it is the locally convex direct limit
$E=\bigcup_{n\in \N}E_n=\dl\,E_n$\vspace{-.3mm} of a
sequence $E_1\sub E_2\sub \cdots$
of Banach spaces such that each inclusion map $E_n\to E_{n+1}$
is a compact linear operator.
Then $E$ is Hausdorff (cf.\ \cite[\S7.3, Satz]{Flo}
and $E=\dl\,E_n$\vspace{-.6mm}
as a topological space \cite[\S7.1, Satz]{Flo}.
The ascending sequence
$(E_n)_{n\in \N}$ can always be chosen
such that, for a suitable norm
on $E_n$ defining its topology,
all closed balls $\wb{B}_r^{E_n}(x)$,
$x\in E_n$, $r>0$, are compact
in $E_{n+1}$ (cf.\ \cite[\S7.3, Proposition~1]{Hog}).
It it clear from the definition that
finite direct products of Silva spaces
are Silva spaces; this will be useful later.
\end{numba}
\begin{numba}\label{defkom}
A Hausdorff topological space~$X$ is called
a \emph{$k_\omega$-space} if there exists an ascending
sequence $K_1\sub K_2\sub \cdots$ of compact subsets
of~$X$ such that $X=\bigcup_{n\in \N}K_n$
and $U\sub X$ is open if and only if
$U\cap K_n$ is open in~$K_n$, for each $n\in \N$
(i.e., $X=\dl\,K_n$\vspace{-.6mm} as a topological space).
Then $(K_n)_{n\in \N}$ is called a
\emph{$k_\omega$-sequence} for~$X$.
For background information concerning
$k_\omega$-spaces with a view towards direct limit
constructions, see \cite{GaG}.
\end{numba}
%
%
%\ma{abund}
\begin{example}\label{abund}
The dual space $E'$
of every metrizable locally convex space~$E$
is a $k_\omega$-space when
equipped with the topology of compact
convergence (cf.\ \cite[Corollary~4.7
and Proposition~5.5]{Aus}).
We write $E'_c$ if this topology is used.
\end{example}
%
%
%\ma{silvak}
\begin{example}\label{silvak}
Every Silva space~$E$
is a $k_\omega$-space.
In fact, $E$ is reflexive
by \cite[\S9, Satz~6]{Flo}
and thus $E\isom (E'_b)'_b$,
using the topology of bounded convergence.
Since $E'_b$ is a Fr\'{e}chet-Schwartz space \cite[\S9, Satz~6]{Flo},
bounded subsets in $E'_b$ are relatively compact
and hence $(E'_b)'_b=(E'_b)'_c$.
But $(E'_b)'_c$ is $k_\omega$, by Example~\ref{abund}.
\end{example}
\begin{example}
If $E$
is an infinite-dimensional Banach space
(or, more generally, a Fr\'{e}chet space
which is not a Schwartz space),
then $E'_c$ is a locally convex space
which is a $k_\omega$-space but not a Silva-space.
In fact, if $E'_c$ was Silva,
then $(E'_c)'_b$ would be a Fr\'{e}chet-Schwartz space
\cite[\S9, Satz~6]{Flo}.
Here $(E'_c)'_b=(E'_c)'_c$
since bounded subsets of Silva spaces
are relatively compact \cite[\S7.6]{Flo}.
Since $E\isom (E'_c)'_c$ holds for every
Fr\'{e}chet space \cite[Proposition~15.2]{Bzc},
we deduce that $E$
is a Fr\'{e}chet-Schwartz space,
contradicting our hypotheses.
\end{example}
The following facts concerning $k_\omega$-spaces will be used:
%
%
%\ma{baseSilva}
\begin{la}\label{baseSilva}
\begin{itemize}
\item[\rm(a)]
If $X$ and $Y$ are $k_\omega$-spaces,
then also $X\times Y$.
\item[\rm(b)]
Let $((X_n)_{n\in \N},(i_{n,m})_{n\geq m})$
be a direct sequence of $k_\omega$-spaces
and continuous maps $i_{n,m}\colon X_m\to X_n$,
with direct limit Hausdorff
topological space~$X$.
Then $X$ is a $k_\omega$-space.
If each $i_n$ is injective,
then the direct limit topological space
$\dl\,X_n$\vspace{-.7mm} is Hausdorff.
\item[\rm(c)]
Let $((E_n)_{n\in \N},(i_{n,m})_{n\geq m})$
be a direct sequence of locally convex
spaces which are~$k_\omega$-spaces,
and continuous linear maps $i_{n,m}\colon E_m\to E_n$.
Then the Hausdorff locally
convex direct limit~$E$
coincides with the direct limit
Hausdorff topological space
$($as discussed in {\rm (b))}.
\end{itemize}
\end{la}
\begin{proof}
(a) See, e.g., \cite[Proposition~4.2\,(c)]{GaG}.

(b) The case of injective
direct sequences is covered
by \cite[Proposition~4.5]{GaG}.
In the general case,
let $i_n\colon X_n\to X$
be the limit map, and $\wb{X}_n:=i_n(X_n)$,
equipped with the quotient
topology, which is Hausdorff and hence
$k_\omega$ by \cite[Proposition~4.2]{GaG}.
Then $X=\dl\,\wb{X}_n$\vspace{-1mm}
as a topological space,
and so $X$ is $k_\omega$ by
\cite[Proposition~4.5]{GaG}.

(c) As in (b), after passing
to Hausdorff quotients we may
assume that each $i_n$ (and $i_{n,m}$)
is injective.
But this case is \cite[Proposition~7.12]{GaG}.
\end{proof}
The following situation arises
frequently:
We are given a map
\mbox{$f\colon U\to F$,}
where $E$ and $F$ are Hausdorff locally convex topological
$\K$-vector spaces
and $U$ a subset of~$E$;
we would like to show
that $U$ is open and~$f$ is $C^r_\K$ for some $r\in \N_0\cup\{\infty\}$.
We are given the
following information:
$E$ is the Hausdorff
locally convex direct limit
of a
sequence $((E_n)_{n\in\N}, (i_{n,m})_{n\geq m})$
of Hausdorff locally convex spaces~$E_n$
and continuous homomorphisms
$i_{n,m}\colon E_m\to E_n$,
with limit maps $i_n\colon E_n\to E$.
Also,
$U=\bigcup_{n\in \N}i_n(U_n)$,
where $U_n\sub E_n$ is an open
subset such that $i_{n,m}(U_m)\sub U_n$
if $n\geq m$. Finally, we assume that
$f_n:=f\circ i_n|_{U_n}\colon U_n\to F$
is $C^r_\K$ for each $n\in \N$.\vfill\pagebreak

%
%
%\ma{checksmth}
\begin{la}\label{checksmth}
In the preceding situation,
suppose that {\rm (a)} or {\rm (b)} holds:
\begin{itemize}
\item[\rm (a)]
Each $E_n$ is a Banach space
and each of the linear maps $i_{n,m}$
is compact.
\item[\rm (b)]
Each $E_n$ is a Silva space
or, more generally, a $k_\omega$-space.
\end{itemize}
Then $U$ is open in~$E$
and $f$ is $C^r_\K$.
\end{la}
\begin{proof}
Let $N_n$ be the kernel of the limit map
$E_n\to E$.
Then $\wb{E}_n:=E_n/N_n$ is
a Banach space (resp., $k_\omega$-space).
Let $q_n\colon E_n\to\wb{E}_n$ be the quotient map,
and $\wb{i}_{n,m}\colon \wb{E}_m\to\wb{E}_n$
be the continuous linear map
determined by $\wb{i}_{n,m}\circ q_m=q_n\circ i_{n,m}$
(which is again a compact operator in case~(a)).
Then $E=\dl\,\wb{E_n}$\vspace{-.4mm}
as a locally convex space,
together with the continuous linear
maps $\wb{i}_n\colon \wb{E}_n\to E$ determined
by $\wb{i}_n\circ \, q_n=i_n$.
The set $\wb{U}_n:=q_n(U_n)$ is open in
$\wb{E}_n$ and $f_n$ factors
to a map $g_n:=f\circ \wb{i}_n \colon \wb{U}_n\to F$
determined by $g_n\circ q_n|_{U_n}=f_n$,
which is $C^r_\K$ by \cite[Lemma~10.4]{Ber}.
Furthermore, $U=\bigcup_{n\in \N}\wb{i}_n(\wb{U}_n)$.
After replacing $E_n$ by $\wb{E}_n$,
we may thus assume now that each $i_{n,m}$
is injective.\\[2.7mm]
Recall that
$E=\dl\,E_n$\vspace{-.4mm}
as a topological space
(see {\bf\ref{recsilva}}
in the situation of~(a),
resp., Lemma~\ref{baseSilva} (b) and~(c)
in the situation of~(b)).
Hence $U$ is open in~$E$
and $U=\dl\,U_n$,\vspace{-.4mm}
by Lemma~\ref{baseDL}.\\[2.5mm]
To see that~$f$ is $C^r_\K$, we may assume
that $r\in \N_0$, and proceed by induction.
If $r=0$, then $f$ is continuous
since $U=\dl\,U_n$\vspace{-.4mm}
as a topological space.
Now assume that the assertion holds for~$r$
and that each $f\circ i_n$ is $C^{r+1}_\K$.
Then $f$ is continuous.
Given $x'\in U$ and $y'\in E$,
there exists $n\in \N$ and $x\in U_n$, $y\in E_n$
such that $x'=i_n(x)$, $y'=i_n(y)$.
Then the directional
derivative $d(f\circ i_n)(x,y)=\frac{d}{dt}\big|_{t=0}
f(i_n(x+ty))=\frac{d}{dt}\big|_{t=0}f(i_n(x)+ti_n(y))
=\frac{d}{dt}\big|_{t=0}f(x'+ty')=df(x',y')$
exists. The preceding calculation shows that
%
%\ma{reintp}
\begin{equation}\label{reintp}
df\circ (i_n\times i_n)|_{U_n\times E_n}\; =\;
d(f\circ i_n)
\end{equation}
for each $n\in \N$,
which is a $C^r_\K$-map.
Since $E\times E=\dl\,E_n\times E_n$\vspace{-.4mm}
is a locally convex direct limit
of the form described in (a), resp., (b)
(see {\bf\ref{recsilva}}
resp., Lemma~\ref{baseSilva}\,(a)),
the map $df$ is $C^r_\K$ by induction
and hence $f$ is $C^{r+1}_\K$.
\end{proof}
Note that $E$ in Lemma~\ref{checksmth}
is a $k_\omega$-space, by~{\bf\ref{silvak}}, resp.,
Lemma~\ref{baseSilva} (b) and~(c).
%
%
%\ma{gdDLSilva}
\begin{prop}\label{gdDLSilva}
Let $G=\bigcup_{n\in \N}G_n$
be a $C^\infty_\K$-Lie group admitting a direct
limit chart and assume that
at least one of the following conditions
is satisfied:
\begin{itemize}
\item[\rm (i)]
$G_n$ is a Banach-Lie group for each $n\in \N$,
and the inclusion map $L(G_m)\to L(G_n)$
is a compact linear operator,
for all $m<n$.
\item[\rm (ii)]
$L(G_n)$ is a $k_\omega$-space,
for each $n\in \N$.
\end{itemize}
Then $G=\dl\,G_n$\vspace{-.4mm}
as a
topological space, topological group,
$C^\infty_\K$-Lie group, and as a $C^r_\K$-manifold,
for each $r\in \N_0\cup\{\infty\}$.
\end{prop}
\begin{proof}
Let $\phi\colon U\to V\sub L(G)$ be a direct limit chart
around~$1$,
where $U=\bigcup_{n\in \N}U_n$, $V=\bigcup_{n\in \N}V_n$
and $\phi=\bigcup_{n\in \N}\phi_n$
for charts $\phi_n\colon U_n\to V_n\sub L(G_n)$
of $G_n$ around~$1$.
Suppose that $f\colon G\to X$ is a map
to a topological space (resp., $C^r_\K$-manifold)~$X$,
such that $f|_{G_n}$ is continuous
(resp., $C^r_\K$), for each $n\in \N$.
Then $(f\circ \phi)|_{V_n}$
is continuous (resp., $C^r_\K$)
for each $n\in\N$ and hence
$f\circ \phi$ is continuous (resp.,
$C^r_\K$) and thus also $f|_U$,
by Lemma~\ref{checksmth}.
Given $x\in G$, applying the same argument
to $h\colon G\to X$, $h(y):=f(xy)$,
we see that $h|_U$
(and hence also $f|_{xU}$)
is continuous, resp., $C^r_\K$.
Hence $f$ is continuous (resp., $C^r_\K$).
We have shown that $G=\dl\,G_n$\vspace{-.4mm}
as a topological space and as a $C^r_\K$-manifold.
The remaining direct limit properties follow.
\end{proof}
\section{Groups of germs of Lie group-valued maps}\label{secgems}
%
%\ma{secgems}
%
%
In this section, we begin our discussion
of the Lie group $\Gamma(K,H)$ of germs
of analytic mappings with values in a Banach-Lie group~$H$,
where $K$ is a non-empty compact subset
of a metrizable locally convex space~$X$.
For $X$ and~$H$ finite-dimensional,
$\Gamma(K,H)$ is modelled
on a Silva space, and we obtain a prime example for the
type of direct limit groups just discussed in
Proposition~\ref{gdDLSilva}\,(i).
This facilitates a complete clarification
of the direct limit properties
of $\Gamma(K,H)$ (Proposition~\ref{propsgemsilva}).
In Section~\ref{secprset},
we develop tools to tackle
$\Gamma(K,H)$ also for infinite-dimensional $X$ and~$H$
(see Section~\ref{beysilva}).
%
%
%\ma{defgem1}
\begin{numba}\label{defgem1}
Let $H$ be a Banach-Lie group over $\K\in \{\R,\C\}$
and $K\not=\emptyset$ a compact subset
of a metrizable locally convex topological
$\K$-vector space~$X$.
Then the group $G:=\Gamma(K,H)$
of germs $[\gamma]$ of $\K$-analytic
maps $\gamma\colon U\to H$ on open neighbourhoods $U\sub X$
of~$K$ is a $C^\omega_\K$-Lie group in a natural way,
with the multiplication of germs
induced by pointwise
multiplication of functions (see~\cite{HOL}).
We now recall the relevant aspects
of the construction of the Lie group
structure, starting with the case $\K=\C$.
In this case,
$\Gamma(K,H)$ is modelled on the locally
convex direct limit
$\Gamma(K,L(H))=\dl\, \Hol_b(U_n,L(H))$,\vspace{-.3mm}
where $U_1\supseteq U_2\supseteq\cdots$
is a fundamental sequence of open
neighbourhoods of~$K$.
Here $\Hol_b(U_n,L(H))=:A_n$ is
the Banach space of bounded holomorphic
functions from $U_n$ to~$L(H)$,
equipped with the supremum norm $\|.\|_{A_n}$.
The space $\Gamma(K,H)$ is Hausdorff~\cite[\S2]{HOL}.
We can (and will always) assume that
each connected component
of $U_n$ meets~$K$;
then all bonding maps
$j_{n,m}\colon \Hol_b(U_m,L(H))\to \Hol_b(U_n,L(H))$,
$\gamma\mto \gamma|_{U_n}$, $n\geq m$,
are injective and hence also all
limit maps $j_n\colon \Hol_b(U_n,L(H))\to\Gamma(K,L(H))$,
$\gamma \mto [\gamma]$.
We occasionally identify $\gamma\in \Hol_b(U_n,L(H))$
with $j_n(\gamma)=[\gamma]$.
\end{numba}
%
%
%\ma{defgem2}
\begin{numba}\label{defgem2}
If $\K=\R$, choose open
neighbourhoods $U_n$ as before
and set $\wt{U}_n:=U_n+iV_n\sub X_\C$,
where $(V_n)_{n\in \N}$
is a basis of open, balanced
$0$-neighbourhoods in~$X$.
Set $C_n:=\{\gamma\in \Hol_b(\wt{U}_n,L(H)_\C)\colon
\gamma(U_n)\sub L(H)\}$ and
\[
A_n\; :=\;
\{\gamma|_{U_n}\colon \gamma\in C_n\}\,.
\]
Then $C_n$ is a closed real vector subspace
of $\Hol_b(\wt{U}_n,L(H)_\C)$.
Because $\gamma\in \Hol_b(\wt{U}_n,L(H)_\C)$
is uniquely determined by $\gamma|_{U_n}$,
we see that $\|\gamma|_{U_n}\|_{A_n}:=\|\gamma\|_\infty$
(supremum-norm) for $\gamma$ as before
with $\gamma(U_n)\sub L(H)$
defines a norm $\|.\|_{A_n}$
on~$A_n$ making it a Banach space
isomorphic to~$C_n$.
We give $\Gamma(K,L(H))$
the vector topology making it
the locally convex direct limit\,
$\dl\,A_n$.\vspace{-.7mm}
\end{numba}
%
%
%\ma{defgem3}
\begin{numba}\label{defgem3}
We may assume that the norm $\|.\|$
on $\ch:=L(H)$ defining its topology has been
chosen such that $\|[x,y]\|\leq\|x\|\cdot\|y\|$
for all $x,y\in \ch$.
Choose $\ve\in \;]0, \frac{1}{2}\log 2]$ such that
$\exp_H|_{B_\ve^\ch(0)}$ is a diffeomorphism onto an open
identity neighbourhood in~$H$.
Then the Baker-Campbell-Hausdorff
series converges to a $\K$-analytic function
$*\colon B_\ve^{\ch}(0)\times B_\ve^{\ch}(0)\to \ch$
(see \cite[Ch.\,II, \S7, no.\,2]{Bou}).
We choose
$\delta \in \;]0,\ve]$
such that $B_\delta^\ch(0)*B_\delta^\ch(0)\sub B_\ve^\ch(0)$.
Define
\[
\Exp\colon \Gamma(K,L(H))\to \Gamma(K,H)\,,\quad
\Exp([\gamma]):=[\exp_H\circ \,\gamma]\,.
\]
Then
$\Gamma(K,H)$ can be given a $\K$-analytic
Lie group structure such that $\Psi:=\Exp|_Q$
is a $C^\omega_\K$-diffeomorphism
onto an open identity neighbourhood
in $\Gamma(K,H)$, where
$Q:=\{ [\gamma]\in \Gamma(K,L(H))\colon \gamma(K)\sub B_\delta^\ch(0)\}$
(see \cite[\S5]{HOL}).
We set $P:=\Psi(Q)$ and $\phi:=\Psi^{-1}$.
\end{numba}
%
%
%\ma{defgem4}
\begin{numba}\label{defgem4}
We now show that
$\Gamma(K,H)=\bigcup_{n\in \N}G_n$
for some Banach-Lie groups~$G_n$,
and that $\phi\colon P\to Q$
is a direct
limit chart.
To this end, let
$C^\omega(U_n,H)$ be the group
of all $\K$-analytic
$H$-valued maps on~$U_n$.
The BCH-series
defines a $\K$-analytic function
$*\colon B_\ve^{A_n}(0)\times B_\ve^{A_n}(0)\to A_n$,
such that $B_\delta^{A_n}(0)*B_\delta^{A_n}(0)\sub
B_\ve^{A_n}(0)$.
The map $\Exp_n\colon A_n \to C^\omega(U_n,H)$,
$\Exp_n(\gamma):=\exp_H\circ\,\gamma$
is injective on $B_\ve^{A_n}(0)$,
and application of point evaluations
shows that $\Exp_n(\gamma*\eta)=\Exp_n(\gamma)\Exp_n(\eta)$
for all $\gamma,\eta\in B_\ve^{A_n}(0)$.
Set $Q_n:=B_\delta^{A_n}(0)$ and $P_n:=\Exp_n(Q_n)$.
Now standard arguments show that the subgroup
$G_n^0$ of $C^\omega(U_n,H)$
generated by $\Exp_n(A_n)$
can be made a Banach-Lie group
with Lie algebra~$A_n$
and such that $\Exp_n|_{Q_n}$
is a $C^\omega_\K$-diffeomorphism
onto~$P_n$, which is open in~$G_n^0$
(cf.\ \hspace*{-.6mm}Proposition \hspace*{-.5mm}18 \hspace*{-.2mm}in
\hspace*{-.3mm}\cite[Ch.\,III, \S1, no.\,9]{Bou}).
\hspace*{-.4mm}Thus $\phi_n:=(\Exp_n|_{Q_n}^{P_n})^{-1}\colon
P_n\!\to\! Q_n$
is a chart for~$G_n^0$.
If $\K=\C$, let $G_n$
be the group of all $\gamma\in C^\omega(U_n,H)$
such that $\,\sup\,\{\|\Ad_{\gamma(x)}^H\|\colon x\in U_n\}<\infty$,
a condition which ensures that\vspace{-.2mm}
the Lie algebra homomorphism
$A_n\to A_n$, $\eta\mto (x\mto \Ad^H_{\gamma(x)}(\eta(x)))$,
is continuous linear.
If $\K=\R$, let $G_n$
be the group of all $\gamma\in C^\omega(U_n,H)$
such that $U_n\to \Aut(\ch_\C)$,\vspace{-.3mm}
$x\mto (\Ad_{\gamma(x)}^H)_\C$
has a complex analytic extension
$\wt{U}_n\to \Aut(\ch_\C)$
which is bounded.
Then $G_n^0\sub G_n$,
and standard arguments provide
a unique $\K$-analytic
manifold structure on~$G_n$
making it a Banach-Lie group with
$G_n^0$ as an open subgroup
(cf.\ \cite[Ch.\,III, \S1, no.\,9, Prop.\,18]{Bou}).
The restriction map
$i_{n,m}\colon G_m\to G_n$, $\gamma\mto\gamma|_{U_n}$
is an injective homomorphism
for $n\geq m$, which is $\K$-analytic
because $\Exp_n\circ \, j_{n,m}=i_{n,m}\circ
\Exp_m$ with $j_{n,m}\colon A_m\to A_n$
continuous linear.
Likewise,
$i_n\colon
G_n\to\Gamma(K,H)$, $\gamma\mto [\gamma]$ is an injective
homomorphism and $\K$-analytic
because $\Exp\circ \, j_n=i_n\circ \Exp_n$.
We identify $G_n$ with its image
$i_n(G_n)$ in $\Gamma(K,H)$.
Then $\Gamma(K,H)=\bigcup_{n\in \N}G_n$.
To see this, let $[\gamma]\in \Gamma(K,H)$.
If $\K=\C$, then
$U_n\to \Aut(\ch)$,
$x\mto \Ad_{\gamma(x)}^H$
is bounded for some $n\in\N$.
If $\K=\R$, then
$U_n\to \Aut(\ch_\C)$,
$x\mto (\Ad_{\gamma(x)}^H)_\C$
has a bounded complex analytic extension
to $\wt{U}_n$ for some~$n$.
Since $Q=\bigcup_{n\in \N}Q_n$,
$P=\bigcup_{n\in \N}P_n$ and
$\phi=\bigcup_{n\in \N}\phi_n$,
we see that $\phi$ is a direct limit chart.
\end{numba}
%
%
%\ma{defgem5}
\begin{numba}\label{defgem5}
If $\dim(X)<\infty$,
we assume that $U_{n+1}$ is relatively compact
in~$U_n$, for each $n\in \N$
(and $V_{n+1}$ relatively compact in $V_n$,
if $\K=\R$).
If, furthermore, $\dim(H)<\infty$,
then $L(i_{n,m})=j_{n,m} \colon A_m\to A_n$,
$\gamma\mto\gamma|_{U_n}$ is a compact operator
whenever $n>m$.
If $\K=\C$, this is a simple consequence of Montel's Theorem;
if $\K=\R$, it follows from the compactness
of the corresponding restriction map
$\Hol_b(\wt{U}_m,L(H)_\C)\to \Hol_b(\wt{U}_n,L(H)_\C)$.
\end{numba}
%
%
%\ma{propsgemsilva}
\begin{prop}\label{propsgemsilva}
For $\K\in \{\R,\C\}$,
consider a Lie group of germs
$\Gamma(K,H)$ as in {\bf \ref{defgem1}--\ref{defgem3}}
and let $G_n$ be as in {\bf\ref{defgem4}}.
If $X$ and $H$ are finite-dimensional,
then $\Gamma(K,H)=\dl\,G_n$\vspace{-.8mm} in the categories
of $C^\infty_\K$-Lie groups, topological groups,
topological spaces and $C^r_\K$-manifolds,
for each $r\in \N_0\cup\{\infty\}$.
\end{prop}
\begin{proof}
{\bf\ref{defgem4}} and
{\bf\ref{defgem5}}
guarantee the
conditions
of
Proposition~\ref{gdDLSilva}\,(i).
\end{proof}
\begin{rem}
If $\dim(H)<\infty$
and $X$ is an infinite-dimensional
Fr\'{e}chet-Schwartz space,
then $\Gamma(K,L(H))$
still is a Silva space
(cf.\ \cite[Theorem~7]{BaM}).
The preceding proposition extends to this situation.
\end{rem}
\begin{rem}
By Lemma~\ref{complf},
we have $\Gamma(K,H)_0=\bigcup_{n\in \N}G_n^0$,
with $G_n^0$ as in {\bf\ref{defgem4}}.
Replacing $\Gamma(K,H)$ by its connected
component $\Gamma(K,H)_0$
and $G_n$ by $G_n^0$,
all of the results of Proposition~\ref{propsgemsilva}
(and likewise those of Proposition~\ref{triumph}
and Corollary~\ref{gemnonsilva} below)
remain valid,
by trivial modifications of the proofs.
In many cases, $\Gamma(K,H)$
is connected (e.g.,
if~$H$ is connected and $K$ a singleton);
then simply $\Gamma(K,H)=\dl\,G_n^0$\vspace{-.4mm}
in all relevant categories.
\end{rem}
\section{\!\!Tools to identify direct limits
of topological \hspace*{-2mm}groups}\label{secprset}
%
%\ma{secprset}
%
We describe a criterion
ensuring that a topological group $G=\bigcup_{n\in \N}G_n$
is the direct limit topological group~$\dl\,G_n$.\vspace{-.3mm}
In combination with Theorem~\ref{reduce},
this facilitates to identify
Lie groups as direct limits in the category
of Lie groups, under quite weak hypotheses.
The criterion,
requiring that ``product sets are large,''
is satisfied in all
situations known to the author.
\begin{defn}
Let $G$ be a topological group
which is a union
$G=\bigcup_{n\in \N}G_n$
of an ascending sequence
$G_1\sub G_2\sub\cdots$
of topological groups such that
all of the inclusion maps
$G_m\to G_n$ and $G_n\to G$
are continuous.
We say that \emph{product sets are large in~$G$}
if the \emph{product map}
\[
\pi\colon {\textstyle \prod_{n\in \N}^*G_n}\to G\,,
\quad (g_n)_{n\in \N}\mto g_1g_2\cdots g_N
\quad
\mbox{if $\,g_n=1$ for all $\,n>N$,}
\]
takes identity neighbourhoods
in the weak direct product
to identity neighbourhoods in~$G$.
If the \emph{product map}
\[
\wt{\pi}\colon {\textstyle \prod_{n\in \N\cup(-\N)}^*G_n}\to G\,,
\quad (g_n)_{n\in \N\cup(-\N)}\mto g_{-N}\cdots g_{-1}g_1g_2\cdots g_N
\]
(with $N$ so large that
$g_n=1$ whenever $|n|>N$)
takes identity neighbourhoods to such,
then we say that \emph{two-sided product sets
are large in~$G$}.
\end{defn}
\begin{rem}
Thus product sets are large in
$G=\bigcup_{n\in \N}G_n$
if and only if
$\bigcup_{n\in \N}U_1U_2\cdots U_n$
is an identity neighbourhood in~$G$,
for each choice of identity neighbourhoods
$U_k\sub G_k$, $k\in \N$.
If product sets are large in~$G$,
then also two-sided product sets are large.
\end{rem}
The following observation provides
first examples with large product sets.
\begin{prop}\label{notintr}
Consider a topological group
$G=\bigcup_{n\in \N}G_n$
such that $G=\dl\,G_n$\vspace{-.4mm}
as a topological space.
Then product sets are large in~$G$.
\end{prop}
\begin{proof}
Consider a product set
$U=\bigcup_{n\in \N}U_1U_2\cdots U_n$
with $U_n$ an open identity neighbourhood
in~$G_n$. Then $U_1\cdots U_n$
is open in $G_n$ and thus $U$
is open in~$G$, by Lemma~\ref{baseDL}.
\end{proof}
\begin{example}
The Lie groups $G=\bigcup_{n\in \N}G_n$
considered in Proposition~\ref{gdDLSilva}
satisfy the hypothesis of Proposition~\ref{notintr},
whence product sets are large~in~$G$.
\end{example}
The following observation provides more interesting
examples:
\begin{rem}
If the product map
$\pi\colon \prod_{n\in \N}^*G_n\to G$
admits a local section
$\sigma\colon U\to \prod^*_{n\in \N}G_n$
on an identity neighbourhood $U\sub G$
which is continuous at~$1$
and takes $1$ to $1$,
then product sets are large in~$G$.
This condition is satisfied
in particular if $\pi$ admits a
continuous (or smooth) local section $\sigma$ around $1\in G$,
such that $\sigma(1)=1$.
\end{rem}
\begin{example}
Consider a test function group
$C^r_c(M,G)$, where~$M$ is
a $\sigma$-compact finite-dimensional
$C^r$-manifold
and $G$ a Lie group
modelled on a locally convex space.
Let $(A_n)_{n\in \N}$ be any exhaustion
of~$M$ by compact sets.
Then the product
map~$\pi$ of
$C^r_c(M,G)=\bigcup_{n\in \N}C^r_{A_n}(M,G)$
admits a smooth local section around~$1$
taking~$1$ to~$1$,
and hence product sets are large in
$C^r_c(M,G)=\bigcup_{n\in \N}C^r_{A_n}(M,G)$.
To see this, we use a map\linebreak
$\Phi\colon C^r_c(M,G)\subseteq \Omega\to
\prod_{n\in \N}C^r_{K_n}(M,G)$
as described in the Fragmentation\linebreak
Lemma~\ref{fragtestf}.
Pick a sequence $m_1<m_2<\cdots$
of positive integers such that
$K_n\sub A_{m_n}$ for each $n\in \N$,
and let $\psi_n\colon C^r_{K_n}(M,G)\to C^r_{A_{m_n}}(M,G)$
be the inclusion map, which is a smooth
homomorphism.
Consider the map $\psi\colon \prod_{n\in \N}^* C^r_{K_n}(M,G)\to
\prod_{n\in \N}^* C^r_{A_n}(M,G)$
sending $\gamma=(\gamma_n)_{n\in \N}$
to $(\eta_k)_{k\in \N}$,
where $\eta_k:=\gamma_n$
if $k=m_n$ for some (necessarily unique)
$n\in \N$ and $\eta_k:=1$ otherwise.
Then $\psi$ is smooth (cf.\ \cite[Proposition~7.1]{MEA})
and $\pi\circ \psi\circ \Phi=\id_\Omega$.
Thus $\sigma:=\psi\circ\Phi$ is the desired
smooth section for~$\pi$.
\end{example}
\begin{example}
Using Lemma~\ref{fragment},
the same argument
shows that the product map~$\pi$ of
$\Diff_c(M)=\bigcup_{n\in \N} \Diff_{A_n}(M)$
admits a smooth local section around~$1$
which takes~$1$ to~$1$,
for each $\sigma$-compact smooth manifold~$M$
and exhaustion $(A_n)_{n\in \N}$ of~$M$
by compact sets.
Thus product sets are large in
$\Diff_c(M)=\bigcup_{n\in \N} \Diff_{A_n}(M)$.
\end{example}
%
%
%\ma{prodDL}
\begin{prop}\label{prodDL}
If two-sided product sets
are large in
a topological group $G=\bigcup_{n\in \N}G_n$
$($in particular, if product sets are large
in~$G)$,
then $G=\dl\,G_n$\vspace{-1mm}
in the category
of topological groups.
\end{prop}
\begin{proof}
Let $f\colon G\to H$ be a homomorphism
to a topological group~$H$
such that $f_n:=f|_{G_n}$ is continuous
for all $n\in \N$.
Consider $h \colon \!\prod_{n\in \N\cup(-\N)}^*G_{|n|}\to H$,
\[
(x_n)_{n\in \N\cup (-\N)}\mto
f_N(x_{-N})\cdots f_1(x_{-1})
f_1(x_1)\cdots f_N(x_N),
\]
with $N$ so large that
$x_n=1$ for all $n\in \N\cup(-\N)$ such that $|n|>N$.
A simple modification
of the proof of Lemma~\ref{prodmaps} shows that~$h$
is continuous at~$1$.
Therefore,
for each identity
neighbourhood $V\sub H$ there exists
a family $(U_n)_{n\in \N\cup(-\N)}$ of identity
neighbourhoods $U_n\sub G_{|n|}$
such that $h(U)\sub V$,
for $U:=\prod_{n\in \N\cup (-\N)}^*U_n$.
Let $\wt{\pi}\colon \prod_{n\in \N\cup (-\N)}^*G_n\to G$ be the product map.
Then $h(U)=f(\tilde{\pi}(U))$
where $\tilde{\pi}(U)\sub G$ is an identity neighbourhood
because two-sided product sets
are large in~$G$.
As a consequence, $f$ is continuous.
\end{proof}
\begin{rem}
Following \cite[\S3.1]{Hir},
an ascending sequence $G_1\leq G_2\leq \cdots$
of topological groups with continuous inclusion maps
is said to satisfy the ``passing through assumption''
(PTA, for short),
if each $G_n$ has a basis
of symmetric identity neighbourhoods~$U$
such that, for each $m>n$
and identity neighbourhood $V\sub G_m$,
there exists an identity neighbourhood
$W\sub G_m$ such that $WU\sub UV$
(cf.\ also \cite{TSH} for a slightly different,
earlier concept).
If Condition PTA is satisfied,
then the two-sided product sets
(or ``bamboo-shoot neighbourhoods'')
$\bigcup_{n\in \N}U_{n}\cdots U_1U_1\cdots U_n$
form a basis of identity neighbourhoods
for the topology $\cO$ on $G=\bigcup_{n\in \N}G_n$
making~$G$ the direct limit topological
group (cf.\ \cite[Proposition~2.3]{TSH}).
In this case, $\cO$ is called
the ``bamboo-shoot topology'' in~\cite{TSH}.
Hence \emph{two-sided product sets
are large in the direct limit topological group~$G$
if condition~PTA is satisfied.}
\end{rem}
\begin{rem}
Let $G=\bigcup_{n\in \N}G_n$ be a topological group
such that $(G_n)_{n\in \N}$ satisfies condition~PTA.
A priori, this only provides information
concerning the direct
limit group topology~$\cO$;
it does not help us to see that the \emph{given}
topology on~$G$ coincides with~$\cO$.
\end{rem}
%
%
%\ma{banaPTA}
\begin{prop}\label{banaPTA}
Assume that $G_1\sub G_2\sub\cdots$
is an ascending sequence
of Banach-Lie groups, such that each inclusion map
$G_n\to G_{n+1}$ is a smooth homomorphism.
Then $(G_n)_{n\in \N}$ satisfies~the~PTA.
\end{prop}
\begin{proof}
For each $n\in \N$, fix a norm
$\|.\|_n$ on $L(G_n)$
defining its topology and such that
$\|[x,y]\|_n\leq \|x\|\cdot\|y\|$
for all $x,y\in L(G_n)$.
The sets $U^{(n)}_\ve:=\exp_{G_n}(B_\ve^{L(G_n)}(0))$,
$\ve>0$,
form a basis of identity neighbourhoods
in~$G_n$.
For each $m>n$,
the set $B_\ve^{L(G_n)}(0)$
is bounded in $L(G_m)$
and thus $M_m:=\sup \|B_\ve^{L(G_n)}(0)\|_m<\infty$.
Given $\delta>0$,
set $\tau:= e^{-M_m}\delta$.
Then $\|\Ad^{G_m}_u(y)\|_m=\|e^{\ad_x^{L(G_m)}}.y\|_m
\leq e^{M_m}\|y\|_m$ for each
$u\in U^{(n)}_\ve$ and $y\in L(G_m)$,
where $u=\exp_{G_n}(x)$ with
$x\in B_\ve^{L(G_n)}(0)$, say.
Thus
%
%\ma{llecal}
\begin{equation}\label{llecal}
\Ad_u^{G_m}(B_\tau^{L(G_m)}(0))\;\sub\;
B_\delta^{L(G_m)}(0)\quad
\mbox{for each $u\in U_\ve^{(n)}$.}
\end{equation}
Given $w\in U_\tau^{(m)}$
and $u\in U_\ve^{(n)}$,
say $w=\exp_{G_m}(y)$ with
$y\in B_\tau^{L(G_m)}(0)$,
we see that
$wu= uu^{-1}wu=uu^{-1}\exp_{G_m}(y)u=u\exp_{G_m}(\Ad^{G_m}_{u^{-1}}(y))
\in U^{(n)}_\ve U^{(m)}_\delta$,
using~(\ref{llecal}).
Hence $U^{(m)}_\tau U^{(n)}_\ve\sub
U^{(n)}_\ve U^{(m)}_\delta$.
We have verified the PTA.
\end{proof}
\section{Unit groups of direct limit algebras}\label{secunit}
%
%\ma{secunit}
%
The following proposition generalizes
\cite[Theorem~1]{Eda} (where all inclusion maps are
isometries)
and complements it by a Lie theoretic
perspective.
%
%\ma{Edamer}
\begin{prop}\label{Edamer}
Let $A_1\sub A_2\sub\cdots$ be an ascending sequence
of unital Banach algebras $A_n$
over~$\K$,
such that each inclusion map $A_n\to A_{n+1}$
is a continuous homomorphism of unital
algebras.
Then the following holds:
\begin{itemize}
\item[\rm (a)]
The locally convex direct limit
topology makes $A:=\bigcup_{n\in \N}A_n$
a locally m-convex topological algebra.
Its unit group $A^\times=\bigcup_{n\in \N}A_n^\times$
is open, and is a topological group
when equipped with the topology induced by~$A$.
\item[\rm (b)]
$A^\times=\dl\,A_n^\times$\vspace{-.9mm}
as a topological group.
\item[\rm (c)]
Product sets are large in
$A^\times\!=\!\bigcup_{n\in \N}A_n^\times$,
and $(A_n^\times)_{n\in \N}$
satisfies the PTA.
\end{itemize}
If $A$ is Hausdorff
$($which is automatic if the direct
sequence is strict$)$,
then $A^\times$ is a $C^\omega_\K$-Lie group
and $A^\times=\dl\,A_n^\times$\vspace{-1.1mm}
as a $C^\infty_\K$-Lie group.
\end{prop}
\begin{proof}
The PTA holds by Proposition~\ref{banaPTA}.
By \cite[Theorem~1]{DaW},
$A$ is a locally m-convex topological
algebra, i.e., the vector topology of~$A$
can be defined by a family of sub-multiplicative
seminorms (see \cite{Mca}). 
Thus $A^\times$ is
a topological group.
It is also known that
$A^\times$ is open
(Wengenroth\linebreak
communicated a proof \cite{Wen}),
but we need not use this fact here,
as an\linebreak
alternative proof is part
of the following arguments.
Consider a product set
$P:=\bigcup_{n\in \N}U_1\cdots U_n$,
with identity neighbourhoods $U_n\sub A_n^\times$.
After shrinking $U_n$, we may assume
that $U_n=\one+B^{A_n}_{\ve_n}(0)$
for some $\ve_n\in \; ]0,\frac{1}{2}]$.
Then
$\|x^{-1}\|\leq \sum_{k=0}^\infty 2^{-k}\leq 2$
for each $x\in U_n$, whence
$U_n^{-1}$ is bounded in $A_n$
and hence also in $A_m$ for each
$m\geq n$.
Thus $M_{m,n}:=\sup\{\|x^{-1}\|_m\colon x\in U_n\}<\infty$.
\emph{We claim that}
\[
\one+\bigcup_{n\in \N}\sum_{k=1}^nB_{\delta_k}^{A_k}(0)
\;\sub\; P\,,
\]
\emph{where $\delta_1:=\ve_1$
and $\delta_n:=M_{n,1}^{-1}\cdots M_{n,n-1}^{-1}\ve_n$
for integers $n\geq 2$.}
If this claim is true,
then $P$ is a neighbourhood of $\one$ in~$A$,
whence $A^\times$ an open subset
of~$A$ (cf.\ \cite[Lemma~2.6]{ALG})
and product sets are large in $A^\times$.
Therefore $A^\times=\dl\,A_n^\times$\vspace{-.5mm}
as a topological group (Proposition~\ref{prodDL}).
Since $A^\times$ is open in~$A$ and inversion
is continuous, $A^\times$ is a $C^\omega_\K$-Lie group
provided~$A$ is Hausdorff \cite[Proposition~3.2,
resp., 3.4]{ALG}.
The identity map $A^\times\to A^\times$
being a direct limit chart,
Theorem~\ref{reduce}
shows that $A^\times=\dl\,A_n^\times$\vspace{-.9mm}
also as a Lie group.\\[3mm]
\emph{Proof of the claim.}
We show
that $\one + \sum_{k=1}^nB_{\delta_k}^{A_k}(0)\sub U_1\cdots U_n\sub P$,
by induction on~$n$.
If $n=1$, then $\one+B_{\delta_1}^{A_1}(0)=U_1\sub P$.
Let $n\geq 2$ now and suppose that
$\one+\sum_{k=1}^{n-1}B_{\delta_k}^{A_k}(0)\sub U_1\cdots U_{n-1}$.
Let $y_k\in B_{\delta_k}^{A_k}(0)$
for $k\in \{1,\ldots,n\}$.
There are $x_j\in U_j$
for $j\in \{1,\ldots, n-1\}$
such that
\[
y\; :=\; \one+y_1+\cdots + y_{n-1}\; =\; x_1 \cdots x_{n-1}\,.
\]
Set $x_n:=y^{-1}(y+y_n)=\one+y^{-1}y_n$.
Then $\|y^{-1}y_n\|_n=\|x_{n-1}^{-1}\cdots x_1^{-1}y_n\|_n
< M_{n,n-1}\cdots M_{n,1}\delta_n=\ve_n$
and thus $x_n\in U_n$.
By construction, $\one + \sum_{k=1}^ny_k=y+y_n=yx_n
=x_1\cdots x_n$.
\end{proof}
If the direct sequence $A_1\sub A_2\sub\cdots$
in Proposition~\ref{Edamer} is strict, then
$A^\times\not=\dl\,A_n^\times$\vspace{-.8mm}
as a topological space
unless each $A_n$ is finite-dimensional
or the sequence $A_n$ becomes
stationary
(by Yamasaki's Theorem, see Remark~\ref{recallyam}).\\[3mm]
We also have a variant for not necessarily unital
associative algebras und algebra homomorphisms
which need not take units to units
(if units do exist).
Recall that if $A$ is an associative
$\K$-algebra (where $\K$ is $\R$ or $\C$),
then $A_e:=\K e\oplus A$ is a unital
algebra via $(re+a)(se+b)=rs e+(rb+sa+ab)$.
Then $(A,\diamond)$ with
$a\diamond b:=a+b-ab$
is a monoid with neutral element~$0$,
whose unit group is denoted $Q(A)$.
The inverse of $a\in Q(A)$ is called
the quasi-inverse of~$a$ and denoted $q(a)$.
The map $(A,\diamond)\to (A_e,\cdot)$, $a\mto e-x$
is a homomorphism of monoids
(see, e.g., \cite[\S2]{ALG}
for all of this).
%
%
%\ma{Edamer2}
\begin{prop}\label{Edamer2}
Let $A_1\sub A_2\sub\cdots$ be a sequence
of $($\/not necessarily unital\/$)$ associative Banach algebras
over~$\K$,
such that each inclusion map $A_n\to A_{n+1}$
is a continuous algebra homomorphism.
Then we have:
\begin{itemize}
\item[\rm (a)]
The locally convex direct
limit topology makes
$A:=\bigcup_{n\in \N}A_n$
a locally m-convex associative
topological algebra
with an open group
$Q(A)$ of quasi-invertible
elements and a
continuous
quasi-inversion map\linebreak
$q\colon Q(A)\to A$.
Thus $Q(A)$
is a topological group.
\item[\rm (b)]
$Q(A)=\dl\,Q(A_n)$\vspace{-.9mm}
as a topological group.
\item[\rm (c)]
Product sets are large in
$Q(A)= \bigcup_{n\in \N}Q(A_n)$,
and $(Q(A_n))_{n\in \N}$
satisfies the PTA.
\end{itemize}
If the
locally convex direct limit topology
on $A$ is Hausdorff,
then $Q(A)$ is a $C^\omega_\K$-Lie group
and $Q(A)=\dl\,Q(A_n)$\vspace{-.3mm}
in the category of $C^\infty_\K$-Lie groups.
\end{prop}
\begin{proof}
(a) The locally convex direct
limit topology on $A_e=\dl\, (\K \oplus A_n)$\vspace{-.4mm}
makes $A_e$ the direct product
$\K\times A$, where $A=\dl\, A_n$\vspace{-.3mm}
carries the locally convex direct limit topology.
Since $A_e$ is a
topological algebra with open unit
group and continuous inversion,
it follows that $Q(A)$ is open in~$A$
and~$q$ is continuous \cite[Lemma~2.8]{ALG}.
Since~$A_e$ is locally m-convex,
so is~$A$.\vspace{1mm}

(b) and (c): Consider a product set
$P:=\bigcup_{n\in \N}B^{A_1}_{\ve_1}(0)\diamond \cdots \diamond
B^{A_n}_{\ve_n}(0)$,
with $\ve_n\in \; ]0,\frac{1}{2}]$.
In~$A_e$, we then have
$e-P=\bigcup_{n\in \N}U_1\cdots U_n$
with $U_n=e-B^{A_n}_{\ve_n}(0)
=e+B^{A_n}_{\ve_n}(0)$.
Extend the norm on $A_n$ to $(A_n)_e$
via $\|re+a\|:=|r|+\|a\|$,
and define $M_{m,n}$ and $\delta_m$
as in the proof of Proposition~\ref{Edamer}.
Re-using the arguments from the proof
just cited, we see that
$e+\bigcup_{n\in \N}\sum_{k=1}^nB_{\delta_k}^{A_k}(0)
\sub e-P$.
Hence
$P$ is a neighbourhood of $0$ in~$A$
and thus product sets are large in~$Q(A)$.
Therefore
$Q(A)=\dl\,Q(A_n)$\vspace{-.5mm}
as a topological group, by Proposition~\ref{prodDL}.
The PTA holds by Proposition~\ref{banaPTA}.\\[1.7mm]
Now assume that $A$ is Hausdorff.
Then $(A_e)^\times\isom Q(A_e)$
is a $C^\omega_\K$-Lie group and $Q(A)=Q(A_e)\cap A$
(see \cite[Lemma~2.5]{ALG})
is a subgroup and submanifold of $Q(A_e)$
and therefore a $C^\omega_\K$-Lie group as well.
The identity map $Q(A)\to Q(A)$
being a direct limit chart,
Theorem~\ref{reduce}
shows that $Q(A)=\dl\,Q(A_n)$\vspace{-.7mm}
as a $C^\infty_\K$-Lie group.
\end{proof}
\section{Lie groups of germs beyond the Silva case}\label{beysilva}
%
%\ma{beysilva}
%
%
%
For all $X$, $K$, $H$
as in {\bf\ref{defgem1}} and $(G_n)_{n\in \N}$
as in {\bf\ref{defgem4}}, we show:
%
%
%\ma{triumph}
\begin{prop}\label{triumph}
Product sets are large in $\Gamma(K,H)=\bigcup_{n\in \N}G_n$.
\end{prop}
%
%
%\ma{banabana}
%
\begin{la}\label{banabana}
Let $(\ch,\|.\|)$ be a Banach-Lie algebra over~$\K$
and $R>0$ such that the BCH-series
converges to a $\K$-analytic mapping
$B_R^\ch(0)\times B_R^\ch(0)\to \ch$.
Then there exist $r\in \;]0,R]$,
a $\K$-analytic map
$F\colon B_r^\ch(0)\times B_r^\ch(0)\to B_R^\ch(0)$
and $C>0$ such that
%
%\ma{rewr}
%\ma{esti}
%
\begin{eqnarray}
x+y & = & x*F(x,y)\quad\mbox{and}\label{rewr}\\
\|F(x,y)-y\| & \leq & C\, \|x\|\,\|y\|\,,\label{esti}
\end{eqnarray}
for all $x,y\in B_r^\ch(0)$.
\end{la}
\begin{proof}
For $r\in \;]0,R]$ sufficiently small,
$F(x,y):=(-x)*(x+y)$ is defined
for all $x,y\in B_r^\ch(0)$,
$F(x,y)\in B_R^\ch(0)$ holds, and (\ref{rewr}).
Since $F(x,0)=0$ and $F(0,y)=y$
for all $x,y\in B_r^\ch(0)$,
the second order
Taylor expansion of~$F$
entails~(\ref{esti}),
after shrinking~$r$ further if necessary (see \cite[Lemma~1.7]{PAD}).
\end{proof}
Multiple products
with respect to the BCH-multiplication~$*$
are formed recursively via
$x_1*\cdots*x_n:=(x_1*\cdots*x_{n-1})*x_n$
(provided that all partial products are defined).\\[3mm]
{\bf Proof of Proposition~\ref{triumph}.}
Let $R>0$, $r$, $C$ and $*$
be as in Lemma~\ref{banabana},
applied with $\ch:=L(H)$ if $\K=\C$
(resp., $\ch:=L(H)_\C$ if $\K=\R$).
After shrinking~$r$,
we may assume that $Cr\leq\frac{1}{2}$.
Now let $(W_n)_{n\in \N}$ be a sequence
of identity neighbourhoods $W_n\sub G_n$.
After shrinking $W_n$, we may assume
that $W_n=\Exp_n(B^{A_n}_{\ve_n}(0))$
for some $\ve_n>0$, with
$A_n$ as in {\bf \ref{defgem1}} (resp.,
{\bf\ref{defgem2}}).
Let $\delta_n:= \min\{r2^{-n},\frac{\ve_n}{2}\}$
for $n\in \N$. Then
$S:= \bigcup_{n\in \N}\sum_{k=1}^n B^{A_n}_{\delta_n}(0)$
is a $0$-neighbourhood in $\Gamma(K,L(H))$.
We claim that
%
%\ma{giveslprod}
%
\begin{equation}\label{giveslprod}
\pi\Big({\textstyle \prod_{n\in \N}^* W_n}\Big)
\; \supseteq\;
\Exp(S)\,,
\end{equation}
where $\pi\colon \prod_{n\in \N}^*G_n\to\Gamma(K,H)$
is the product map.
Therefore
$\pi\big({\textstyle \prod_{n\in \N}^* W_n}\big)$
is an identity neighbourhood in $\Gamma(K,H)$,
and hence product sets are large in
$\Gamma(K,H)$.
To prove the claim,
let $z\in \Exp(S)$.
Thus $z=\Exp(\sum_{n=1}^\infty[\gamma_n])$
for some sequence
$(\gamma_n)_{n\in \N}\in \bigoplus_{n\in \N}B^{A_n}_{\delta_n}(0)$.
Choose $N\in \N$ such that $[\gamma_n]=0$
for all $n>N$.
If $\K=\C$,
we set $\eta_1:=\gamma_1$.
If $n\in \{2,\ldots, N\}$,
then $\eta_n(x):=F(\gamma_1(x)+\cdots+\gamma_{n-1}(x),\gamma_n(x))$
makes sense for each $x\in U_n$
because\linebreak
$\|\gamma_1(x)+\cdots+\gamma_{n-1}(x)\|<r$
and $\|\gamma_n(x)\|<r$,
and defines a bounded holomorphic
function $\eta_n\colon U_n\to L(H)$.
We have $\|\eta_n\|_\infty\leq (1+Cr)\|\gamma_n\|_\infty
\leq \frac{3}{2}\delta_n<\ve_n$
by (\ref{esti}),
whence $\eta_n\in B_{\ve_n}^{A_n}(0)$ and thus $\Exp_n(\eta_n)\in W_n$.
Furthermore,
$\gamma_1(x)+\cdots+\gamma_n(x)=
(\gamma_1(x)+\cdots+\gamma_{n-1}(x))*\eta_n(x)$
and hence
$\gamma_1(x)+\cdots+\gamma_n(x)=
\eta_1(x)*\eta_2(x)*\cdots*\eta_n(x)$
for each $n\in\{1,\ldots, N\}$ and $x\in U_n$,
by induction.
Therefore,
\[
\gamma_1(x)+\cdots+\gamma_N(x)\;=\;\eta_1(x)*\eta_2(x)*\cdots*\eta_N(x)
\quad \mbox{for each $x\in U_N$,}
\]
entailing that
$\sum_{n=1}^\infty[\gamma_n]
=[x\mto \eta_1(x)*\cdots*\eta_N(x)]$
and $z=\Exp(\sum_{n=1}^\infty[\gamma_n])
=[x\mto \exp_H(\eta_1(x)*\cdots*\eta_N(x))]
=[x\mto \exp_H(\eta_1(x))\cdots\exp_H(\eta_N(x))]
=\Exp([\eta_1])\cdots\Exp([\eta_N])\in
\pi\big({\textstyle \prod_{n\in \N}^* W_n}\big)$.
Thus (\ref{giveslprod}) holds.
If $\K=\R$,\vspace{-.4mm}
we use the unique bounded holomorphic
extension $\tilde{\gamma}_n\colon \wt{U}_n\to \ch$
of each $\gamma_n$
in place of $\gamma_n$
to define bounded holomorphic maps
$\tilde{\eta}_n\colon \wt{U}_n\to \ch$
along the lines of the construction of~$\eta_n$.
Set $\eta_n:=\tilde{\eta}_n|_{U_n}$.
Then $\|\eta_n\|_{A_n}=\|\tilde{\eta}_n\|_\infty<\ve_n$
for each~$n$ and we see as above
that $z=\Exp([\eta_1])\cdots\Exp([\eta_N])\in
\pi\big({\textstyle \prod_{n\in \N}^* W_n}\big)$.\vspace{2.5mm}\Punkt

\noindent
Combining Proposition~\ref{prodDL}
with Theorem~\ref{reduce}, we obtain:
%
%
%\ma{gemnonsilva}
%
\begin{cor}\label{gemnonsilva}
$\Gamma(K,H)=\dl\, G_n$\vspace{-.3mm}
holds
in the category of $C^\infty_\K$-Lie groups,
and in the category of topological groups.
\Punkt
\end{cor}
\section{Construction of Lie group structures on\\
direct limit groups}\label{secconstr}
%
%
%\ma{secconstr}
%
Consider an abstract group
$G=\bigcup_{n\in \N}G_n$ which is
the union of an ascending sequence
$G_1\sub G_2\sub \cdots$ of $C^\infty_\K$-Lie groups $G_n$,
such that the inclusion maps $i_{n,m}\colon G_m\to G_n$
(for $m\leq n$) are $C^\infty_\K$-homomorphisms
and each $G_n$ is a subgroup of~$G$.
In this section, we describe
conditions which facilitate
to construct a $C^\infty_\K$-Lie group structure on~$G$
such that $G=\dl\,G_n$\vspace{-.4mm}
as a $C^\infty_\K$-Lie group.
For finite-dimensional Lie groups
$G_n$, such a Lie group structure
has been constructed in~\cite{FUN}
(cf.\ \cite{NRW1}, \cite{NRW2},
\cite[Theorem~47.9]{KaM}
and \cite{DIR} for special cases).
The conditions formulated in this section
apply just as well to suitable infinite-dimensional
Lie groups~$G_n$.
%
%
%\ma{cacha}
\begin{numba}\label{cacha}
We shall always assume that
$G$ has a \emph{candidate for a direct limit chart},
viz.\ we assume that
there exist charts
$\phi_n\colon G_n\supseteq U_n\to
V_n\sub L(G_n)$ of $G_n$
around~$1$ for $n\in \N$
such that $U_m\sub U_n$
and $\phi_n|_{U_m}=L(i_{n,m})\circ \phi_m$
if $m\leq n$,
and $V:=\bigcup_{n\in \N}V_n$ is open in
the locally convex direct limit
$E:=\dl\, L(G_n)$,\vspace{-.4mm}
which we assume Hausdorff.
Here, we identify $L(G_m)$
with the image of $L(i_{n,m})$
in~$L(G_n)$;
this is possible because $L(i_{n,m})$ is injective
by an argument as in Remark~\ref{remdlcha}\,(a).
We define $U:=\bigcup_{n\in \N}U_n$
and $\phi:=\dl\, \phi_n\colon U\to V\sub E$.\vspace{-.4mm}
\end{numba}
It is natural to wonder whether
$\phi$ (or its restriction to a
smaller identity neighbourhood in~$U$)
can always be used as a chart around~$1$
for a Lie group structure on $G$.
Unfortunately, the answer
is {\em negative} (without extra hypotheses):
even if $\phi$ is globally defined on all of~$G$,
it need not make $G=\bigcup_n G_n$ a Lie group.
%
%
%\ma{expathoalg}
\begin{example}\label{expathoalg}
Let $A_1\subseteq A_2\subseteq \cdots$
be any ascending sequence of locally convex
unital associative topological algebras
such that 1.\,the inclusion maps
are homomorphisms of unital algebras
and topological embeddings;
2.\,the locally convex direct limit topology
renders the algebra multiplication on the union
$A:=\bigcup_{n\in \N}A_n$ discontinuous at $(1,1)$;
and 3.\,the unit group $A^\times$ is open in~$A$
(see \cite[\S10]{ALG} for such algebras).
Then $A^\times=\bigcup_{n \in \N}A_n^\times$
is a union of Lie groups
and~$A^\times$ admits the global chart $\phi:=\id_{A^\times}$,
which is a candidate for a direct limit chart
around~$1$.
However, $A^\times$ is not a Lie group because
the group multiplication is discontinuous at $(1,1)$.
\end{example}
We now describe additional
requirements ensuring that
the question just posed has an
affirmative answer. They are satisfied
in many situations.
%
%
%\ma{constrLie}
\begin{prop}\label{constrLie}
Consider an abstract group $G=\bigcup_{n\in \N}G_n$
which is the union of an ascending
sequence of $C^\infty_\K$-Lie groups.
Assume that $G$ admits a candidate
$\phi\colon U\to V\sub E:=\dl\,L(G_n)$\vspace{-.4mm}
for a direct limit chart,
and assume that condition
{\rm (i)} or {\rm (ii)}
from Proposition~{\rm \ref{gdDLSilva}}
is satisfied.
Then there exists a unique $C^\infty_\K$-Lie group
structure on $G$ making $\phi|_W$
a direct limit chart for~$G$ around~$1$,
for an open identity neighbourhood $W\sub U$.
\end{prop}
\begin{proof}
This is a special case of
Lemma~\ref{mantogp} below, applied
with $M_n:=G_n$
(we use $k:=n$, $A:=B:=G_n$
in (c) and $k:=n$, $A:=G_n$ in (d)).
\end{proof}
\begin{rem}
By Proposition~\ref{gdDLSilva},
the Lie group structure
described in Proposition~\ref{constrLie}
makes $G$ the direct limit
$\dl\,G_n$\vspace{-.4mm}
as a $C^\infty_\K$-Lie group, topological group,
topological space, and as a $C^r_\K$-manifold,
for each $r\in \N_0\cup\{\infty\}$.
\end{rem}
Consider a set $M$ which is an ascending union
$M=\bigcup_{n\in \N}M_n$ of $C^\infty_\K$-manifolds,
and $x\in M$.
Changing {\bf\ref{cacha}} in the
obvious way,\footnote{Replace
$G$ by~$M$, $L(G_n)$
by $T_xM_n$, and $L(i_{n,m})$
by $T_x(i_{n,m})$.}
we obtain the definition of a \emph{candidate for
a direct limit chart around~$x$}.
%
%
%\ma{mantogp}
\begin{la}\label{mantogp}
Let $G$ be an abstract group
which is the union $G=\bigcup_{n\in \N}M_n$
of an ascending sequence
$M_1\sub M_2\sub\cdots$
of $C^\infty_\K$-manifolds,
such that $1\in M_1$
and {\rm (a)--(d)}
hold:
\begin{itemize}
\item[\rm (a)]
The inclusion maps
$M_m\to M_n$ are~$C^\infty_\K$
for all $m\leq n$;
\item[\rm (b)]
$G$ admits a candidate $\phi\colon U\to V\sub E
:=\dl\,T_1(M_n)$\vspace{-1mm} for a direct limit chart
around~$1$;
\item[\rm (c)]
For each $n\in \N$ and $x,y\in M_n$,
there exists $k\geq n$ and open neighbourhoods
$A, B\sub M_n$ of $x$, resp., $y$
such that $AB\sub M_k$
and the group multiplication
$A\times B\to M_k$ is~$C^\infty_\K$; and
\item[\rm (d)]
For each $n\in \N$,
there exists $k\geq n$ and an open identity
neighbourhood $A \sub M_n$
such that $A^{-1} \sub M_k$
and the group inversion
$A \to M_k$, $x\mto x^{-1}$ is~$C^\infty_\K$.
\end{itemize}
Furthermore, we assume that {\rm (i)}
or {\rm (ii)} is satisfied:
\begin{itemize}
\item[\rm (i)]
$M_n$ is modelled
on a Banach space
for each $n\in \N$,
and the inclusion map
$T_1(M_m)\to T_1(M_n)$
is a compact operator for all $m<n$.
\item[\rm (ii)]
The modelling locally convex space of each $M_n$
is a $k_\omega$-space.
\end{itemize}
Then there is a unique $C^\infty_\K$-Lie group
structure on~$G$ making
$\phi$ a chart for~$G$ around~$1$.
Furthermore, $G=\dl\,M_n$\vspace{-.4mm}
as a topological
space and as a $C^r_\K$-manifold,
for each $r\in \N_0\cup\{\infty\}$.
\end{la}
\begin{proof}
Suppose that $\phi=\dl\,\phi_n$
with $\phi_n\colon U_n\to V_n$.
Equip $G$ with the topology $\cT$
turning it into the direct limit topological space
$\dl\,M_n$.\vspace{-.4mm}
Given $x\in M_n$,
consider $\lambda_x\colon G\to G$,
$\lambda_x(y):=xy$.
Hypothesis\,(c) implies
that $\lambda_x|_{M_m}$
is continuous for each $m\in \N$.
Hence $\lambda_x$ is continuous,
and hence a homeomorphism.
Likewise, all right translations
are homeomorphisms.
Let $\cS$ be the topology on
$G\times G$ making it the direct limit
$\dl\,(M_n\times M_n)$.
By Lemma~\ref{baseDL},
the topology induced by~$\cT$ on~$U$
makes~$U$ the direct limit
topological space $U=\dl\,U_n$,\vspace{-.4mm}
and $\cS$ induces on $U\times U$
the topology making it the direct
limit $\dl\,(U_n\times U_n)$.\vspace{-.4mm}
This topology is the product topology
on~$U\times U$;
this follows from
the fact that the product topology
on $E\times E$ coincides with
the locally convex direct limit
topology on $E\times E=\dl\,T_1M_n\times T_1M_n$,\vspace{-.4mm}
which makes $E\times E$
the direct limit topological
space $\dl\,T_1M_n\times T_1M_n$\vspace{-.4mm}
by {\bf\ref{recsilva}} (resp., Lemma~\ref{baseSilva} (b)
and (c)).
As a consequence of~(c),
the group multiplication
$\mu$ restricts to
a continuous map
$U\times U=\dl\,(U_n\times U_n)\to G$.\vspace{-.4mm}
Since $\mu$ is continuous
on the identity neighbourhood
$U\times U$ and all left and right
translations are homeomorphisms,
it follows that~$\mu$ is continuous.
Let $\iota\colon G\to G$ be the inversion map.
Since $y^{-1}=(x^{-1}y)^{-1}x^{-1}$
for $x\in M_n$ and $y\in M_n$,
combining (c) and (d)
we see that $\iota|_{M_n}$
is continuous on a neighbourhood
of~$x$ and hence continuous.
Hence $\iota$ is continuous and
hence~$G$ is a topological group, which is Hausdorff
because the intersection of
all identity neighbourhoods is~$\{1\}$.
For $x\in G$, define $\phi_x\colon xU\to V$,
$\phi_x(y)=\phi(x^{-1}y)$.
Given $x,y\in G$,
the map $\phi_x\circ \phi_y^{-1}$
is defined on the open set
$\phi(U\cap y^{-1}xU)$
and takes $z$
to $\phi(x^{-1}y\phi^{-1}(z))$.
In view of (c), we easily deduce from
Lemma~\ref{checksmth}
that $\phi_x\circ \phi_y^{-1}$
is~$C^\infty_\K$.
Hence the charts are compatible and
thus $G$ is a $C^\infty_\K$-manifold.
Since $\phi_{xy}\circ \lambda_x\circ \phi_y^{-1}=\id_V$
is~$C^\infty_\K$,
each left translation map $\lambda_x$
is~$C^\infty_\K$ and hence a $C^\infty_\K$-diffeomorphism.
Let $W\sub U$ be an open, symmetric identity neighbourhood
such that $WW\sub U$.
Replacing continuity by smoothness
in the above arguments, we see
(with the help of Lemma~\ref{checksmth})
that the group multiplication
$W\times W\to U$
and inversion $W\to W$ are~$C^\infty_\K$.
Similarly, (c) implies
that each inner automorphism
$c_x\colon G\to G$ takes some
identity neighbourhood
smoothly into~$U$.
Now standard arguments
provide a unique $C^\infty_\K$-Lie group structure
on $G$ making $W$
an open smooth submanifold
(see, e.g., \cite[Proposition~1.13]{GCX}).
Since $\lambda_x$ is a diffeomorphism
from $W$ onto $xW$ for each $x\in G$,
both for the manifold structure making $G$ a Lie group
and the manifold structure constructed before,
we deduce that the two manifold structures
coincide.
\end{proof}
\begin{rem}
If each $M_n$ is a finite-dimensional
$C^\infty_\K$-manifold in the situation
of Lemma~\ref{mantogp}
and each inclusion map $M_m\to M_n$,
$m\leq n$, a smooth immersion,
then a direct limit chart around~$1$
exists by \cite[Theorem~3.1]{FUN}
and thus condition~(b) of Lemma~\ref{mantogp}
is automatically satisfied.
\end{rem}
\begin{rem}
Inspecting the proof of Lemma~\ref{mantogp},
we recognize that condition (c)
can be replaced by an alternative
condition~(c)$'$:
We require that for all
$x,y\in M_n$, there exists $k\geq n$
and open neighbourhoods $A, B\sub M_n$ of $x$, resp., $y$
such that $Ay\sub M_k$, $xB\sub M_k$ and both
of the maps $\rho_y|_A\colon A\to M_k$
and $\lambda_x|_B\colon B\to M_k$
are~$C^\infty_\K$
(where $\rho_y\colon G\to G$, $g\mto gy$).
Furthermore, we
require the existence of
an open symmetric identity neighbourhood
$W\sub U$ such that $WW\sub U$
and such that (c) holds for all~$n$,
if $x,y\in M_n\cap W$.
\end{rem}
\begin{rem}
In the convenient setting of analysis,
it is easier to construct Lie group
structures on direct limit groups.
In fact, consider an abstract group
$G=\bigcup_{n\in \N}G_n$
which is the union of an ascending sequence
of convenient Lie groups.
We equip the abstract vector space $E:=\dl\,L(G_n)$\vspace{-.4mm}
with the locally convex vector
topology associated with the direct limit bornology.
We assume that the latter is Hausdorff
and require that $G$ admits
a candidate for a direct limit chart
in the convenient sense
(defined as in~{\bf\ref{cacha}}, except
that $V$ only needs to be $c^\infty$-open in~$E$,
and $E$ is topologized as just described).
Finally, we assume
that each bounded subset
in $L(G)$ is a bounded
subset of some $L(G_n)$
(regularity).
Then it is
straightforward to make
$G=\bigcup_n G_n$
a (possibly not smoothly Hausdorff)
Lie group in the sense of convenient differential
calculus,
such that $\phi$ is a chart.\footnote{Using $\phi$
and its translates $\phi_x$ as charts,
we can make $G$ a (possibly not smoothly
Hausdorff) smooth manifold
in the convenient sense.
Given a smooth curve in~$G$,
for each finite~$k$ it locally is a $Lip^k$-map into some
$G_n$ (cf.\ \cite[Corollary~1.8]{KaM}),
entailing that the group inversion
is conveniently smooth
(cf.\ \cite[Theorem~12.8 and Corollary~12.9]{KaM}).
Similarly, the group multiplication is
conveniently smooth.}
Then $G=\dl\,G_n$\vspace{-.5mm}
in the category of
(not necessarily smoothly Hausdorff)
Lie groups in the sense of convenient
differential calculus,
and in the category of
smooth manifolds in this
sense.\footnote{Using the direct limit chart
and the regularity,
we easily see that each smooth
curve in $G$ locally is a $Lip^k$-map into some
$G_n$, for each given $k\in \N_0$.
This implies that $G=\dl\,G_n$\vspace{-1.5mm}
as a $Lip^k$-manifold.
The asserted direct limit properties
follow from this.}
The group multiplication
need not be continuous
(cf.\ Example~\ref{expathoalg}).
\end{rem}
\section{Example:
Lie groups of germs of analytic diffeomorphisms}\label{diffgerm}
%
%\ma{diffgerm}
%
%
Let $K$ be a non-empty compact
subset of $X:=\K^d$, where $d\in \N$,
and
$\GermDiff(K,X)$
be the group of all germs $[\gamma]$ around~$K$
of $\K$-analytic maps $\gamma\colon U\to X$
on an open neighbourhood $U$ of~$K$
such that $\gamma|_K=\id_K$,
$\gamma(U)$ is open in~$X$
and $\gamma\colon U\to\gamma(U)$
is a $\K$-analytic diffeomorphism.
Then $\GermDiff(K,X)$
is a group in a natural way,
with group operation $[\gamma][\eta]:=
[\gamma\circ \eta|_{\eta^{-1}(U)}]$
(for $\gamma\colon U\to X$).
To illustrate the usefulness of Lemma~\ref{checksmth}
and Lemma~\ref{mantogp},
we apply them to turn $\GermDiff(K,X)$
into a $\K$-analytic Lie group,
modelled on the space $\Gamma(K,X)_K$
of germs $[\gamma]$ around~$K$ of $X$-valued
$\K$-analytic maps $\gamma\colon U\to X$
such that $\gamma|_K=0$.\vspace{-3mm}
\subsection*{{\normalsize Monoids of germs of complex analytic self-maps}}
It simplifies the construction (and provides additional
information) to consider in a first step
the monoid $\GermEnd(K,X)$
of all germs $[\gamma]$ around~$K$ of $\K$-analytic
maps $\gamma\colon U\to X$ on an open neighbourhood
of~$K$ such that $\gamma|_K=\id_K$
(with multiplication given by composition
of representatives).
We equip $\GermEnd(K,X)$ with
a $\K$-analytic manifold structure which makes the
monoid multiplication a $\K$-analytic map.
In a second step, we show that
the unit group $\GermDiff(K,X)$
of $\GermEnd(K,X)$ is open
and has a $\K$-analytic inversion map.
Until Remark~\ref{llrem},
we let $\K=\C$.
\begin{numba}
Choose a norm on~$X$.
For $n\in \N$, the sets $U_n:=K+B_{1/n}^X(0)$
form a fundamental sequence
of open neighbourhoods of~$K$ in~$X$.
The supremum norm makes
the space $\Hol_b(U_n,X)$
of bounded $X$-valued $C^\omega_\C$-maps on~$U_n$ a complex
Banach space,
and $\Hol_b(U_n,X)_K:=\{\gamma\in \Hol_b(U_n,X)\colon \gamma|_K=0\}$
is a closed vector subspace.
Let $j_{n,m}\colon
\Hol_b(U_m,X)_K\to\Hol_b(U_n,X)_\K$
be the restriction map,
for $n\geq m$.
Given an open neighbourhood
$U$ of~$K$, let
$\Hol(U,X)$
be the Fr\'{e}chet space of
all $X$-valued $C^\omega_\C$-maps on~$U$
(equipped with the compact-open topology)
and $\Hol(U,X)_K$
be its closed subspace
of functions vanishing on~$K$.
Since $\wb{U_{n+1}}=K+\wb{B}_{\frac{1}{n+1}}^X(0)$
is a compact subset of~$U_n$,
we then have continuous linear
restriction maps
%
%\ma{restrmps}
\begin{equation}\label{restrmps}
\Hol_b(U_n,X)_K\to \Hol(U_n,X)_K\to\Hol_b(U_{n+1},X)_K\,,
\end{equation}
whose composition $j_{n+1,n}$ is a compact operator
due to Montel's Theorem.
Thus, the locally convex
direct limit $\Gamma(K,X)_K=\dl\, \Hol_b(U_n,X)_K$\vspace{-.4mm}
is a Silva space.
Since each connected component
of $U_n$ meets~$K$ and hence meets~$U_{n+1}$,
the Identity Theorem implies
that each bonding map $j_{n,m}$
is injective and hence
also each limit map $j_n\colon \Hol_b(U_n,X)_K\to
\Gamma(K,X)_K$, $\gamma\mto[\gamma]$.
\end{numba}
\begin{numba}
It is useful to note that the map
\begin{equation}\label{maprho}
\rho\colon \Gamma(K,X)_K\to C(K,\cL(X))\,,\quad
[\gamma]\mto \gamma'|_K
\end{equation}
is continuous, where $\cL(X)$ is the Banach algebra
of continuous endomorphisms of~$X$
with the operator norm, and $C(K,\cL(X))$
is given the supremum norm.
Since $\Gamma(K,X)_K=\dl\,\Hol_b(U_n,X)_K$\vspace{-.4mm}
as a locally convex space, this follows
from the fact that the inclusion
maps $\Hol_b(U_n,X)\to\Hol(U_n,X)$
and $\Hol(U_n,X)\to C^\infty_\C(U_n,X)$
are continuous (see \cite[Proposition~III.15]{NeI}).
\end{numba}
\begin{numba}
It is clear from the definition that the map
\[
\Phi\colon \Gamma(K,X)_K\to \GermEnd(K,X)\,,\qquad
\Phi([\gamma]):=[\id_X+\,\gamma]
\]
is a bijection.
We use $\Phi$
to give $\GermEnd(K,X)$
a complex manifold structure
with~$\Phi^{-1}$ as a global chart.
For $\gamma,\eta\in \Gamma(K,X)_K$, we then have
%
%\ma{redcc}
\begin{equation}\label{redcc1}
\Phi^{-1}\big(\Phi([\gamma])\Phi([\eta])\big)\;=\;
[\eta]+[\gamma\circ (\id_X+\,\eta)]\,.
\end{equation}
\end{numba}
%
%
%\ma{propgemdiff}
\begin{prop}\label{propgemdiff1}
$\GermEnd(K,X)$
is a complex
analytic monoid, i.e., the multiplication map
$\GermEnd(K,X)\times \GermEnd(K,X)\to \GermEnd(K,X)$,
$([\gamma],[\eta])\mto [\gamma]\,[\eta]=[\gamma\circ\eta]$
is complex analytic.
\end{prop}
\begin{proof}
For $k\in \N$, let
$\Omega_k$ be the set of all
$[\gamma]\in \Gamma(K,X)_K$
such that $\sup\,\{\|\gamma'(x)\|\colon x\in K\}\,<\, k$.
Then $\Gamma(K,X)_K=\bigcup_{k\in \N}\Omega_k$,
and each of the sets $\Omega_k$ is an open
$0$-neighbourhood in $\Gamma(K,X)_K$,
by continuity of~$\rho$
(from (\ref{maprho}).
Hence, in view of (\ref{redcc1}),
the multiplication map will be $C^\infty_\C$ if the map
\[
f\colon \Gamma(K,X)_K\times \Omega_k\to \Gamma(K,X)_K\,,\quad
f([\gamma],[\eta]):=[\gamma\circ (\id_X+\eta)]
\]
is~$C^\infty_\C$ for each $k\in \N$.
Fix~$k$.
For $n\in \N$, set $\ell_n:=(n+1)(k+1)$,
$m_n:=\ell_n+1$,
$P_n:=\Hol_b(U_n,X)_K$,
and let $Q_n$
be the set of all $\gamma\in \Hol_b(U_n,X)_K$
such that $\,\sup\,\{\|\gamma'(x)\|\colon x\in U_{n+1}\}\,<\, k$.
Then $Q_n$ is open in $\Hol_b(U_n,X)_K$.
Hence
\[
\Gamma(K,X)_K\times \Omega_k\; =\;
\bigcup_{n\in \N}j_n(P_n)\times j_n(Q_n)\,,
\]
where $P_n\times Q_n$
is open in $\Hol_b(U_n,X)_K\times\Hol_b(U_n,X)_K$
for each $n\in \N$
and $(j_{n,m}\times j_{n,m})(P_m\times Q_m)
\sub P_n\times Q_n$ if $m\leq n$.
Using Lemma~\ref{checksmth}\,(a), we see that $f$
will be $C^\infty_\C$ if
$f|_{P_n\times Q_n}$ is $C^\infty_\C$
for each~$n\in \N$ (identifying
$\Hol_b(U_n,X)_K$ with its image in $\Gamma(K,X)_K$
using~$j_n$). 
Note that if $\eta\in Q_n$
and $y\in \wb{U_{\ell_n}}$,
say $y\in \wb{B}_{\ell_n^{-1}}^X(x)$
with $x\in K$,
then $\|(y+\eta(y))-x\|=\|y-x+\eta(y)-\eta(x)\|\leq
\|x-y\|+\|\int_0^1\eta'(x+t(y-x)).(y-x)\,dt\|
\leq (1+k)\|y-x\|\leq \frac{1+k}{\ell_n}<\frac{1}{n}$
and thus $y+\eta(y)\in B_{\frac{1}{n}}^X(x)$.
Hence $(\id_X+\,\eta)(\wb{U_{\ell_n}})\sub U_n$,
enabling us to define
\[
g_n\colon P_n\times Q_n \to \Hol_b(U_{m_n},X)_K\,,\quad
g_n(\gamma,\eta):=\gamma\circ (\id_X+\,\eta)|_{U_{m_n}}\,.
\]
Since $f|_{P_n\times Q_n}=j_{m_n}\circ g_n$,
it remains to show that~$g_n$ is $C^\infty_\C$.
The set
$S_n:=\{\eta\in \Hol(U_n,X)_K\colon (\id_X+\,\eta)(\wb{U_{\ell_n}})
\sub U_n\}$
is open in $R_n:=\Hol(U_n,X)_K$.
Using continuous linear inclusion and
restriction maps and the mapping
\[
h_n \colon R_n\times S_n \to \Hol(U_{\ell_n},X)_K\,,\quad
h_n(\gamma,\eta):=\gamma\circ (\id_X+\,\eta)|_{U_{\ell_n}},\vspace{-.8mm}
\]
we can write $g_n$ as a composition
\[
P_n\times Q_n\to R_n\times S_n\stackrel{h_n\,}{\to}
\Hol(U_{\ell_n},X)_K\to \Hol_b(U_{m_n},X)_K.
\]
Therefore, $g_n$ will be $C^\infty_\C$ if we can show
that~$h_n$ is~$C^\infty_\C$.
To this end, we exploit that $\Hol(U_n,X)=C^\infty_\C(U_n,X)$,
equipped with the compact-open $C^\infty$-topology
(cf.\ \cite[Proposition~III.15]{NeI}).
By \cite[Lemma~11.4]{ZOO}, the map
\[
C^\infty_\C(U_n,X)\times \lfloor \wb{U_{\ell_n}},U_n\rfloor_\infty
\to C^\infty_\C(U_{\ell_n},X)\,,\quad
(\gamma,\eta)\mto \gamma\circ \eta|_{U_{\ell_n}}
\]
is $C^\infty_\C$,
where
$\lfloor \wb{U_{\ell_n}},U_n\rfloor_\infty:=
\{\eta\in C^\infty_\C(U_n,X)\colon
\eta(\wb{U_{\ell_n}})\sub U_n\}$.
It easily follows that
$h_n$ is $C^\infty_\C$
as a map into $\Hol(U_{\ell_n},X)$ and
hence also as a map into
the closed vector subspace
$\Hol(U_{\ell_n},X)_K$ (by Lemma~\ref{closd}).
\end{proof}
\subsection*{{\normalsize Lie groups of germs
of complex analytic diffeomorphisms}}
\begin{la}\label{spotunitgp}
The complex analytic monoid
$\GermEnd(K,X)$ has an open unit group
$\GermEnd(K,X)^\times$.
It is given by
\begin{eqnarray}
\GermEnd(K,X)^\times &= & \GermDiff(K,X)\notag\\
&= &
\{[\gamma]\in \GermEnd(K,X)\colon \gamma'(K)\sub
\GL(X)\}\,.\label{trinity}
\end{eqnarray}
\end{la}
\begin{proof}
If $\gamma\colon U\to \gamma(U)\sub X$
is a diffeomorphism between open neighbourhoods
of~$K$, then $[\gamma]\in \GermEnd(K,X)^\times$
and thus $\GermDiff(K,X)\sub \GermEnd(K,X)^\times$.
If $[\gamma],[\eta]\in \GermEnd(K,X)^\times$
with $[\eta]=[\gamma]^{-1}$,
then $\eta\circ\gamma|_W=\id_W$
on some open neighbourhood $W$ of~$K$,
whence $\gamma|_W$ is injective
and $\gamma'(x)\in \GL(X)$ for each $x\in W$.
Hence $\gamma|_W$ is a diffeomorphism
onto an open neighbourhood of~$K$ and
thus $[\gamma]\in\GermDiff(K,X)$.\\[2.5mm]
If $[\gamma]\in \GermDiff(K,X)$, then
$\gamma'(K)\sub \GL(X)$. As the converse follows
from the next lemma, we see that (\ref{trinity})
holds. Since $C(K,\GL(X))=C(K,\cL(X))^\times$
is open in the Banach algebra $C(K,\cL(X))$
and $\rho$ from (\ref{maprho}) is continuous,
we see that $\{[\gamma]\in \GermEnd(K,X)\colon
\gamma'(K)\sub \GL(X)\}=\rho^{-1}(C(K,\GL(X))$
is open in $\GermEnd(K,X)$.
\end{proof}
\begin{la}\label{nonqu}
Let $X$ be a Banach space over
$\K\in \{\R,\C\}$, $K\sub X$ be a non-empty compact
set, $r\in \N\cup\{\infty,\omega\}$
and $f\colon U\to E$ be a $C^r_\K$-map
on an open neighbourhood~$U$ of~$K$
such that $f|_K$ is injective
and $f'(K)\sub \GL(X)$.
If $r=1$, assume $\dim(X)<\infty$.
Then there is an open neighbourhood
$V\sub U$ of~$K$ such that $f(V)\sub X$
is open and $f|_V$ is a $C^r_\K$-diffeomorphism
onto $f(V)$.
\end{la}
\begin{proof}
Define $U_n:=K+B_{1/n}^X(0)$.
There is $n_0\in \N$ such that $U_n\sub U$
and $\gamma'(U_{n_0})\sub \GL(X)$.
Then $\gamma|_{U_{n_0}}$ is a local
$C^r_\K$-diffeomorphism
by the Inverse Function Theorem
(see \cite[Theorem~5.1]{IM2} if $r\not=\omega$;
the analytic case is well known).
In particular, $\gamma(U_n)$ is open in~$X$,
for each $n\geq n_0$.
If $\gamma|_{U_n}$ fails to be injective
for each $n\geq n_0$,
we find $v_n\not=w_n\in U_n$
such that $f(v_n)=f(w_n)$.
Then $v_n\in B_{1/n}^X(x_n)$ and
$w_n\in B^X_{1/n}(y_n)$ for certain
$x_n,y_n\in K$.
There exist $k_1<k_2<\cdots$
such that $x_{n_k}\to x$
and $y_{n_k}\to y$ for certain $x,y\in K$.
Then $v_{n_k}\to x$ and $w_{n_k}\to y$,
whence $f(x)=f(y)$ by continuity
and hence $x=y$.
Let $W$ be a neighbourhood
of~$x$ on which~$f$ is injective.
Then $v_{n_k}\in W$ and $w_{n_k}\in W$
for large~$k$. Since $v_{n_k}\not=w_{n_k}$
and $f(v_{n_k})=f(w_{n_k})$,
this contradicts the injectivity of $f|_W$.
Hence, there exists $n\geq n_0$
such that $f$ is injective on $V:=U_n$.
Then $f|_V\colon V\to f(V)$
is a $C^r_\K$-diffeomorphism.
\end{proof}
We return to the default
notations of this section.
The following quantitative variant of Lemma~\ref{nonqu}
will be needed:
%
%\ma{reallylocdiff}
\begin{la}\label{reallylocdiff}
Given
$\gamma\in \Hol_b(U_n,X)_K$
with
$C\!:=\!\sup\{\|\gamma'(x)\|\colon x\in U_n\}<1$,
define $\eta:=\id_X+\,\gamma\colon U_n\to X$.
Then $\eta|_{U_{6n}}$ is injective.
\end{la}
\begin{proof}
Let
$v,w\in U_{6n}$ such that
$\eta(v)=\eta(w)$.
There are $x,y\in K$
such that $v\in B_{1/(6n)}^X(x)$
and $w\in B_{1/(6n)}^X(y)$.
Let $[x,v]$ be the line
segment joining $x$ and~$v$.
Since $[x,v]\sub B^X_{1/(6n)}(x)\sub U_n$,
the~Mean Value Theorem yields
$\|\eta(v)-\eta(x)\|\leq
\|v-x\|\cdot \max\{\|\eta'(z)\|\colon z\in [x,v]\}\leq 2\|v-x\|
<1/(3n)$.
Likewise, $\|\eta(w)-\eta(y)\|<1/(3n)$
and thus $\|y-x\|=\|\eta(y)-\eta(x)\|<2/(3n)$,
entailing that $[v,w]\sub B^X_{1/n}(x)\sub U_n$
and therefore $0=\|\eta(v)-\eta(w)\|\geq
\|v-w\|-
\|(\eta-\id_X)(v)-(\eta-\id_X)(w)\|\geq$
$\|v-w\|-C\cdot \|v-w\|$.
If $v\not=w$, we
get the contradiction
$0\geq (1-C)\|v-w\|>0$. Hence $v=w$ and thus
$\eta|_{U_{6n}}$ is injective.
\end{proof}
\begin{numba}
Lemma~\ref{spotunitgp} enables us to consider
$\GermDiff(K,X)$ as an open $C^\infty_\C$-submanifold
of $\GermEnd(K,X)$.
Then $D:=\Phi^{-1}(\GermDiff(K,X))$
is an open $0$-neighbourhood
in $\Gamma(K,X)_K$, and the restriction $\Psi$
of $\Phi^{-1}$ to a map
$\Psi\colon \GermDiff(K,X)\to D$
is a global chart for $\GermDiff(K,X)$.
We set $M_n:=\Psi^{-1}(D\cap \Hol_b(U_n,X)_K)$
and give $M_n$ the complex Banach manifold structure
with $\Psi|_{M_n}\colon M_n\to D\cap \Hol_b(U_n,X)_K$
as a global chart.
\end{numba}
%
%
%\ma{propgemdiff}
\begin{prop}\label{propgemdiff}
$\GermDiff(K,X)$ is a complex
Lie group.
Furthermore, $\GermDiff(K,X)=\dl\,M_n$\vspace{-.4mm}
both as a topological space, as a
$C^r_\C$-manifold, and as a $C^r_\R$-manifold,
for each $r\in \N_0\cup\{\infty\}$.
\end{prop}
\begin{proof}
We check the hypotheses of
Lemma~\ref{mantogp}.
Here (a) is clear and also (b),
since $\Psi$ is a candidate
for a direct limit chart.
Condition\,(i) holds
because the restriction maps
$\Hol_b(U_n,X)_K\to\Hol_b(U_{n+1},X)_K$
are compact.
The validity of
Condition (c) is clear from
the proof of Proposition~\ref{propgemdiff1}.\\[3mm]
\emph{Condition} (d) \emph{concerning the inversion map}:
It suffices to show that, for each $n\in \N$,
there exist $k> \ell> n$ and an open $0$-neighbourhood
$P\sub \Hol_b(U_n,X)_K$ such that
$\eta_\gamma:=(\id_X+\,\gamma)|_{U_\ell}$
is injective for each $\gamma\in P$,
$U_k\sub \eta_\gamma(U_\ell)$,
$g(\gamma):=\eta_\gamma^{-1}|_{U_k}-\id_X\in
\Hol_b(U_k,X)_K$, and that $g\colon P\to \Hol_b(U_k,X)_K$
is~$C^\infty_\C$. Because we can build in
continuous linear restriction maps,
it suffices to find an open $0$-neighbourhood
$Q\sub \Hol(U_n,X)_K$ and $m>\ell>n$
such that $\eta_\gamma:=(\id_X+\,\gamma)|_{U_\ell}$
is injective
for each $\gamma\in Q$,
$U_m \sub \eta_\gamma(U_\ell)$,
$h(\gamma):=(\eta_\gamma^{-1}-\id_X)|_{U_m}\in
\Hol(U_m,X)_K$, and that $h\colon P\to \Hol(U_m,X)_K$
is~$C^\infty_\C$
(then take $k:=m+1$).
We set $m:=12 (n+1)$,
$\ell:=6(n+1)$ and let
$Q$ be the set of all $\gamma\in \Hol(U_n,X)_K$
such that $\sup\{\|\gamma'(x)\|\colon
x\in \wb{U_{n+1}}\}<\frac{1}{2}$.
Then $Q$ is open in $\Hol(U_n,X)_K$
and $\eta_\gamma:=(\id_X+\,\gamma)|_{U_\ell}$
is injective for each $\gamma\in Q$,
by {\bf\ref{reallylocdiff}}.
Furthermore, $B^X_{1/(2\ell)}(x)\sub \eta_\gamma(B_{1/\ell}^X(x))$
for each $x\in K$
by \cite[Theorem~5.3\,(d)]{IM2}
(applied with $A:=\id_X$)
and thus $U_m \sub \eta_\gamma(U_\ell)$.
The map
\[
f\colon Q\times U_\ell \to X\,,\quad f(\gamma,x):=(\id_X+\,\gamma)(x)
\]
is $C^\infty_\C$ since the
evaluation map is $C^\infty_\C$
(see \cite[Proposition~11.1]{ZOO}).
Also, $f_\gamma:=f(\gamma,\sbull)\colon U_\ell\to X$
is injective by the preceding,
and is a local $C^\infty_\C$-diffeomorphism
(by the Inverse Function Theorem).
Furthermore, the map
\[
\psi\colon Q\times U_m\to X\,,\quad
\psi(\gamma,x):=f_\gamma^{-1}(x)=(\id_X+\,\gamma)^{-1}(x)=
\eta_\gamma^{-1}(x)
\]
is $C^\infty_\C$ by the Inverse Function Theorem
with Parameters (see \cite[Theorem~2.3\,(c)]{IMP}
or \cite[Theorem~5.13\,(b)]{IM2}).
Then, by the exponential
law \cite[Lemma~12.1\,(a)]{ZOO}),
also the map
\[
\psi^\vee\colon Q\to C^\infty_\C(U_m,X)=\Hol(U_m,X)\,,\quad
\gamma \mto \psi(\gamma,\sbull)=\eta_\gamma^{-1}|_{U_m}
\]
is $C^\infty_\C$,
whence also $Q\to\Hol(U_m,X)$,
$\gamma\mto \psi^\vee(\gamma)-\id_X=h(\gamma)$ is~$C^\infty_\C$.
As this map takes values in the
closed vector subspace $\Hol(U_m,X)_K$
of $\Hol(U_m,X)$,
Lemma~\ref{closd} shows that also its
co-restriction~$h$ is $C^\infty_\C$,
as required.
\end{proof}
\begin{rem}\label{llrem}
For $K$ a singleton (say, the origin),
it is well known that
$\GermDiff(\{0\},\C^d)$
is a Lie group.
The Lie group structure
has first been constructed in~\cite{Pis}, where 
this group is denoted by $\text{Gh}(d,\C)$.
\end{rem}
\subsection*{{\normalsize Lie groups of germs
of real analytic diffeomorphisms}}
Let $\K=\R$ now. It is clear that
$\Gamma(K,X)_K$ can be identified
with the set of germs $[\gamma]\in \Gamma(K,X_\C)_K$
such that $\gamma(U)\sub X$
for some neighbourhood~$U$ of~$K$ in~$X$.
In the same way, we identify
$\GermEnd(K,X)$ with a subset of $\GermEnd(K,X_\C)$.
We give $\Gamma(K,X)_K$ the topology induced
by $\Gamma(K,X_\C)_K$.
Then $\Gamma(K,X)_K$
is a closed real vector subspace
of $\Gamma(K,X_\C)$ and
$\Gamma(K,X_\C)_K=(\Gamma(K,X)_K)_\C$
as a locally convex space
(cf.\ \cite[\S4.2--4.4]{HOL}).
Then $\Phi^{-1} \colon\GermEnd(K,X_\C)\to
\Gamma(K,X_\C)_K$ (defined as before)
is a global chart of $\GermEnd(K,X_\C)$
such that $\Phi(\GermEnd(K,X))
=\Gamma(K,X)_K$, showing
that $\GermEnd(K,X)$ is a
real analytic submanifold of
$\GermEnd(K,X_\C)$.
As a consequence,
$\GermDiff(K,X)=\GermDiff(K,X_\C)\cap \GermEnd(K,X)$
is open in $\GermEnd(K,X)$.
The real analyticity of the monoid multiplication
and group inversion is inherited by
the submanifolds. Summing up:
\begin{cor}
If $\K=\R$, then
$\GermDiff(K,X)$ is a real analytic
Lie group modelled on $\Gamma(K,X)_K$.
Furthermore, $\GermEnd(K,X)$ is a real analytic
monoid with open unit group
$\GermEnd(K,X)^\times=\GermDiff(K,X)$.\Punkt
\end{cor}
\begin{rem}
It is clear that
$C_n:=\{\gamma\in \Hol_b(U_n,X_\C)\colon \gamma(U_n\cap X)\sub X\}$
is a closed vector subspace,
enabling us to make $A_n:=\{\gamma|_{U_n\cap X}\colon \gamma
\in C_n\}$ a real Banach space isometrically
isomorphic to~$C_n$.
It is easy to see that $\Gamma(K,X)_K=\dl\,A_n$\vspace{-.4mm}
as a locally convex space.
Setting $M_n:=\Phi^{-1}(A_n\cap D)$,
we easily deduce from
the proofs of Lemma~\ref{spotunitgp}
and Proposition~\ref{propgemdiff}
that the conditions of Lemma~\ref{mantogp}
are satisfied. We deduce:
\emph{$\GermDiff(K,X)=\dl\,M_n$\vspace{-.9mm}
as a topological space and as a $C^r_\R$-manifold,
for each $r\in \N_0\cup\{\infty\}$.}
\end{rem}
\begin{rem}
Note that, since $\rho$ in (\ref{maprho})
is continuous linear, the set
$E:=\{[\gamma]\in \Gamma(K,X)_K\colon \gamma'|_K=0\}$
is a closed $\K$-vector subspace of $\Gamma(K,X)_K$
(in both cases, $\K\in \{\R,\C\}$).
Then
\[
\GermDiff(K,X)^*\, :=\, \{[\gamma]\in \GermDiff(K,X)\colon \gamma'|_K=\id_X\}
\]
is a Lie subgroup of $\GermDiff(K,X)$,
since the chart $\Phi$ from above
(resp.,~its restriction the group
of germs of $C^\omega_\R$-diffeomorphisms)
takes $\GermDiff(K,X)^*$
onto $D\cap E=E$.
Because $\rho$ is continuous linear,
the homomorphism
\[
\theta\colon \GermDiff(K,X)\to C(K,\GL(X))\,,\quad
\theta([\gamma]):=\gamma'|_K
\]
is $\K$-analytic. We therefore have an exact sequence
of $C^\omega_\K$-Lie groups
\begin{equation}\label{exct}
\one\,\to\, \GermDiff(K,X)^*\,\hookrightarrow \,
\GermDiff(K,X)\,\stackrel{\theta\,}{\to}\, C(K,\GL(X))\,,
\end{equation}
where $\GermDiff(K,X)^*$ is $C^\omega_\K$-diffeomorphic
to~$E$ and hence contractible.
\end{rem}
\begin{rem}
The Lie group $\GermDiff(\{0\},\R)^*$
has also been discussed in \cite{Lei},
using different notation.
\end{rem}
In a later work, the author
hopes to discuss $\GermDiff(K,X)$
also for $X$ an
infinite-dimensional Banach space.
\section{Covering groups of direct limit groups}\label{covgp}
%
%\ma{covgp}
In this section, we discuss universal covering groups
of direct limit groups.\\[2.5mm]
First, we observe
that direct limit charts
with balanced ranges can always
be built up from chart with
balanced ranges.\\[3mm]
The setting is as follows:
Let $r\in \N_0\cup\{\infty,\omega\}$
and $M=\bigcup_{n\in \N}M_n$
be a $C^r_\K$-manifold
admitting a weak direct limit chart
$($resp., direct limit chart$)$
$\phi\colon U\to V\sub E$
around $x\in M_{n_0}$,
where $E$ is the modelling space
of~$M$.
Identifying the modelling space~$E_n$
of~$M_n$ with a subspace
of~$E$ (viz.\ each $j_n\colon E_n\to E$ and
$j_{n,m}\colon E_m\to E_n$ is the inclusion map),
we then have
$\phi=\bigcup_{n\geq n_0}\phi_n$,
$V=\bigcup_{n\geq n_0}V_n$
and $U=\bigcup_{n\geq n_0}U_n$,
with charts $\phi_n\colon U_n\to V_n\sub E_n$
such that $U_m\sub U_n$ in $m\leq n$ and $\phi_n|_{U_m}=\phi_m$.
%
%
%
%\ma{specDL}
\begin{la}\label{specDL}
If $\phi(x)=0$ and $V$ is balanced
$($i.e., $V=[-1,1]V)$,
then there are balanced $0$-neighbourhoods $V_n'\sub V_n$
such that $V_n'\sub V_{n+1}'$ for each~$n$
and $V=\bigcup_{n\geq n_0}V_n'$.
\end{la}
\begin{proof}
We define $V_n':=\{y\in V_n\colon [{-1},1]y\sub V_n\}$.
Then $0\in V_n'$,
and a simple compactness argument
shows that $V_n'$ is open.
If $y\in V_n'$ and $t\in [{-1},1]$,
then $[{-1},1]ty\sub [{-1},1]y\sub V_n$
and hence $ty\in V_n'$,
showing that $V_n'$ is balanced.
Furthermore, it is clear from the definition
that $V_n'\sub V_{n+1}'$.
To see that $V=\bigcup_{n\geq n_0}V_n'$,
take a point $y\in V$,
say $y\in V_m$.
Then $K:=[{-1},1]y\sub E_m$
is a compact subset.
Since $V$ is balanced, we have
$K\sub V$. Hence $K$ is covered
by the ascending sequence of open subsets
$K\cap V_n$, $n\geq m$,
and hence $K=K\cap V_n$
for some $n\geq m$ by compactness.
Thus $K=[{-1},1]y\sub V_n$
and therefore $y\in V_n'$. This completes the proof.
\end{proof}
\begin{rem}
In the situation of Lemma~\ref{specDL},
set $U_n':=\phi^{-1}(V_n')$
and $\phi_n':=\phi_n|_{U_n'}^{V_n'}$.
Then $\phi=\bigcup_{n\geq n_0}\phi_n'$,
where each chart $\phi_n'$ has balanced range.
\end{rem}
We now consider a connected Lie group
$G=\bigcup_{n\in \N}G_n$
admitting a weak direct limit chart,
where each $G_n$ is a connected Lie group.
We let $q_n\colon \wt{G}_n\to G_n$
and $q\colon \wt{G}\to G$ be
universal covering groups,
and $N_n$ be the kernel
of the homomorphism $j_n\colon
\wt{G}_n\to \wt{G}$ obtained by lifting
$i_n\circ q_n\colon \wt{G}_n\to G$ over~$q$,
where $i_n\colon G_n\to G$ is the inclusion map.
%
%
%\ma{univcov}
\begin{prop}\label{univcov}
In the preceding situation,
the following holds:
\begin{itemize}
\item[\rm (a)]
$\wt{G}=\dl\,\wt{G}_n$\vspace{-1mm}
and $\wt{G}=\dl\, \wt{G}_n/N_n$ as an abstract group;
\item[\rm (b)]
$\wt{G}=\bigcup_{n\in \N}\wt{G}_n/N_n$
admits a weak direct limit chart.
\end{itemize}
If $G=\bigcup_{n\in \N}G_n$
admits a direct limit chart, then
also $\wt{G}=\bigcup_{n\in \N}\wt{G}_n/N_n$
admits a direct limit chart,
and the following holds:
\begin{itemize}
\item[\rm (c)]
$G=\dl\,G_n$
as a topological space
if and only if $\wt{G}=\dl\,\wt{G}_n$\vspace{-.7mm}
as a topological space,
if and only if $\wt{G}=\dl\,\wt{G}_n/N_n$\vspace{-1mm}
as a topological space.
\item[\rm (d)]
If $L(G)$ is $C^r$-regular for
some $r\in \N_0 \cup\{\infty\}$,
then $G=\dl\,G_n$\vspace{-.9mm}
as a $C^r$-manifold
if and only if $\wt{G}=\dl\,\wt{G}_n$\vspace{-.9mm}
as a $C^r$-manifold,
if and only if $\wt{G}=\dl\,\wt{G}_n/N_n$\vspace{-.8mm}
as a $C^r$-manifold.
\end{itemize}
\end{prop}
\begin{proof}
(a) Let $i_{n,m}\colon G_m\to G_n$
be the inclusion map for $n,m\in \N$
with $n\geq m$,
and $j_{n,m}\colon \wt{G}_m\to\wt{G}_n$
be the continuous homomorphism
induced by~$i_{n,m}$,
determined by $q_n\circ j_{n,m}=i_{n,m}\circ q_m$
\cite[Proposition~A.2.32]{HaM}.
It is clear
that $j_{n,n}=\id_{\wt{G}_n}$
and $j_{n,m}\circ j_{m,\ell}=j_{n,\ell}$ if $n\geq m\geq \ell$,
whence the direct limit group
$D:=\dl\,\wt{G}_n$\vspace{-.6mm} can be formed,
with limit maps
$k_n\colon \wt{G}_n\to D$.
It is also clear that $j_n\circ j_{n,m}=j_m$,
because the calculation
$q\circ j_n\circ j_{n,m}=i_n\circ q_n\circ j_{n,m}
=i_n\circ i_{n,m}\circ q_m=i_m\circ q_m=q\circ j_m$
shows that both of the homomorphisms $j_n\circ j_{n,m}$
and $j_m$ are lifts of $q\circ j_m$
over~$q$ (and hence coincide).
Now the direct limit property
of~$D$ provides a unique homomorphism
$j\colon D \to \wt{G}$
such that $j\circ k_n=j_n$ for each $n\in \N$.\\[2.5mm]
To construct a homomorphism $h\colon \wt{G}\to D$
(which will turn out to be $j^{-1}$),
let $\phi\colon U\to V\sub L(G)$ be a direct
limit chart for~$G$ around~$1$ such that
$q(W)=U$ and $q|_W^U$ is a diffeomorphism,
for some open identity neighbourhood
$W\sub \wt{G}$.
Here $U=\bigcup_{n\in \N}U_n$,
$V=\bigcup_{n\in \N}V_n$
and $\phi=\bigcup_{n\in \N}\phi_n$
for certain charts
$\phi_n\colon U_n\to V_n\sub L(G_n)$
of $G_n$ around~$1$.
By Lemma~\ref{specDL},
after shrinking~$U$ we may assume
that~$V$ and each $V_n$
is balanced.
Let $W'\sub W$ be an open
identity neighbourhood such that $W'W'\sub W$.
Set $U':=q(W')$ and $V':=\phi(U')$;
after shrinking~$W'$ if necessary,
we may assume that $V'$
is balanced and
$V'=\bigcup_{n\in \N}V'_n$
for certain balanced open $0$-neighbourhoods
$V_n'\sub V_n$ (by Lemma~\ref{specDL}).
We let $W'':=W'\cap (W')^{-1}$ and $U'':=q(W'')$.\\[2.5mm]
Since $V_n$ is simply connected (being contractible)
and hence also~$U_n$,
the inclusion map $u_n\colon U_n\to G_n$
lifts over $q_n$ to a continuous
map $f_n\colon U_n\to\wt{G}_n$
such that $q_n\circ f_n=u_n$ and $f_n(1)=1$.
If $u_{n,m}\colon U_m\to U_n$
is the inclusion map
for $n\geq m$,
then $q_n\circ f_n\circ u_{n,m}=u_n\circ u_{n,m}=i_{n,m}\circ
u_m$ coincides with
$q_n\circ j_{n,m}\circ f_m=i_{n,m}\circ q_m\circ f_m
=i_{n,m}\circ u_m$,
whence $f_n\circ u_{n,m}=j_{n,m}\circ f_n$.
Therefore the maps $k_n\circ f_n\colon U_n\to D$
induce a map $f :=\dl\, f_n\colon U\to D$,\vspace{-1.1mm}
determined by $f \circ u_n=k_n\circ f_n$.\\[2.5mm]
We claim that
$f(xy)=f(x)f(y)$
for all $x,y\in U'$ such that $xy\in U'$.
If this is true, then
$f(q(x)q(y))=f(q(x))f(q(y))$ in particular
for all $x,y\in W''$ such that $xy\in W''$,
whence the map
$W''\to D$, $x\mto f(q(x))$ on the symmetric, open
identity neighbourhood $W''\sub \wt{G}$
extends uniquely to a homomorphism
$h\colon \wt{G}\to D$
(see \cite[Corollary~A2.26]{HaM}).
To prove the claim,
let $x,y\in U'$
and abbreviate $x_t:=\phi^{-1}(t\phi(x))$,
$y_t:=\phi^{-1}(t\phi(y))$
for $t\in [0,1]$.
There is $m\in \N$ such that $\phi(x),\phi(y)\in V'_m$
and thus $[0,1]\phi(x), [0,1]\phi(y)\sub V'_m$.
Then $g\colon [0,1]\to G_m$, 
$g(t):= x_ty_t$
is a continuous map
such that $g([0,1])\sub U$.
Since the compact set
$g([0,1])$ is the union of the ascending sequence
of open subsets $g([0,1])\cap U_n$, $n\geq m$,
there exists $n\geq m$ such that
$g([0,1])\sub U_n$.
Now both $[0,1]\to \wt{G}_n$,
$t\mto j_{n,m}(f_m(x_t)f_m(y_t))$
and $t\mto f_n(x_ty_t)$
are lifts starting in~$1$ over~$q$ of the continuous
map $[0,1]\to G_n$, $t\mto x_ty_t$.
Hence both lifts coincide,
whence $j_{n,m}(f_m(x)f_m(y))=f_n(xy)$
in particular and thus $f(x)f(y)=f(xy)$.\\[2.5mm]
We now show that $h\circ j=\id_D$.
Since $D=\bigcup_{n\in \N}k_n(\wt{G}_n)$,
we only need to show that
$(h\circ j)(x)=x$
for each $x\in k_n(\wt{G}_n)$.
Since $h\circ j\circ k_n$ is a homomorphism
and $\wt{G}_n$
is generated by $U_n'':=u_n^{-1}(W'')$,
we only need to show that $h(j(x))=x$
for each element of the form $x=k_n(f_n(z))$,
where $z\in U_n''$.
Since both $j_n\circ f_n$
and $(q|_W^U)^{-1}\circ u_n$ lift the inclusion map
$i_n\circ u_n\colon U_n\to G$,
we have $j_n\circ f_n=(q|_W^U)^{-1}\circ u_n$.
Thus, as required:
\begin{eqnarray*}
h(j(x)) &= & h(j(k_n(f_n(z))))
\;=\; h(j_n(f_n(z)))
\; =\; h((q|_{W''}^{U''})^{-1}(z))\\
&=& f(q((q|_{W''}^{U''})^{-1}(z)))
\; =\; f(z)\; =\; (f\circ u_n)(z)
\; =\; k_n(f_n(z))\; =\; x\,.
\end{eqnarray*}
Next, we show that $j\circ h=\id_{\wt{G}}$.
Since $\wt{G}$ is connected and
hence generated by $W''$,
we only need to show
that $j(h(x))=x$ for each $x\in W''$.
Set $z:=q(x)$;
then $z\in U_n'$, say.
Now $j(h(x))=j(f(q(x)))=
j(f(z))=j(f(u_n(z)))=j(k_n(f_n(z)))
=j_n(f_n(z))=(q|_W^U)^{-1}(z)=x$
indeed.\\[2.5mm]
By the preceding, $j$ is an isomorphism
with inverse~$h$.
As a consequence, $\wt{G}$, together
with the homomorphisms~$j_n$ as the limit maps,
is the direct limit $\dl\,\wt{G}_n$\vspace{-1.2mm}
as an abstract group (and hence also as a set).
Note that $N_n$ is discrete
(since $N_n\sub \ker(q_n)$);
hence $\wb{G}_n:=\wt{G}_n/N_n$
is a Lie group.
Let $p_n\colon \wt{G}_n\to \wb{G}_n$
be the quotient morphism
and $\wb{j}_{n,m}\colon \wb{G}_m \to \wb{G}_n$
be the smooth homomorphism
determined by $\wb{j}_{n,m}\circ p_m=p_n\circ j_{n,m}$.
Let $\wb{j}_n\colon \wb{G}_n\to\wt{G}$
be the smooth homomorphism obtained by factoring~$j_n$
over~$p_n$, satisfying
$\wb{j}_n\circ p_n=j_n$.
Given homomorphisms (resp., maps) $g_n\colon
\wt{G}_n\to H$
such that $g_n\circ j_{n,m}=g_m$,
each $g_n$ factors to a homomorphism (resp., map)
$\wb{g}_n\colon \wb{G}_n\to H$,
such that $\wb{g}_n\circ p_n=g_n$.
Then $\wb{g}_n\circ \wb{j}_{n,m}=\wb{g}_m$.
It is clear from this that
$\wt{G}=\dl\,\wb{G}_n$\vspace{-1mm}
as an abstract group,
together with the maps $\wb{j}_n$.
It also follows from this that
$\wt{G}=\dl\,\wt{G}_n$\vspace{-1mm}
as a Lie group,
topological group, topological space,
resp., $C^r$-manifold if and only if
$\wt{G}=\dl\,\wb{G}_n$\vspace{-.4mm}
in the respective category.\vspace{1.4mm}

(b) Let $\phi\colon U\to V$ be a weak direct limit chart
(resp., direct limit chart),
where $V=\bigcup_{n\in \N}V_n$
with each~$V_n$ balanced.
Let $\psi\colon V\to \wt{G}$ be the lift
of $\phi^{-1}$ over~$q$
with $\psi(0)=1$,
and $\psi_n\colon V_n\to \wb{G}_n$
be the lift of $\phi_n^{-1}\colon V_n\to U_n$
$\sub G_n$
over the covering morphism $\pi_n\colon \wb{G}_n\to G_n$,
$xN_n\mto q_n(x)$,
such that $\psi_n(0)=1$.
Then $\psi^{-1}$ is a chart for~$\wt{G}$,
$\psi_n^{-1}$ is a chart for $\wt{G}_n$,
and it is clear that
$\psi^{-1}=\bigcup_{n\in \N}\psi_n^{-1}$
(identifying $\wb{G}_n$ with its image in $\wt{G}_n$
under $j_n$, as before).
Hence $\psi^{-1}$ is a weak direct
limit chart (resp., direct limit chart)
for~$\wt{G}$ around~$1$.\vspace{1.4mm}

(c) and (d):
In view of\,(b),
it is immediate from
Theorem~\ref{usefcor}
that $G=\dl\,G_n$\vspace{-.8mm}
as a topological space
if and only if $\wt{G}=\dl\,\wb{G}_n$ as a topological space,\vspace{-.4mm}
and that $G=\dl\,G_n$\vspace{-.4mm}
as a $C^r$-manifold (in~(d))
if and only if $\wt{G}=\dl\,\wb{G}_n$
as a $C^r$-manifold.
By the final remarks
in the proof of~(a), the remaining
equivalences are also valid.
\end{proof}
\section{Open problems}\label{secopen}
%
%\ma{secopen}
%
%
Although the article clarifies
the direct limit properties of many
important examples,
and sheds some light on general
situations,
various related questions had to remain open.
We here compile some of these,
starting with open problems
concerning concrete examples
and then gradually turning to problems
concerning direct limits in general.
\begin{numba}
For finite $r\in \N_0$,
$s\in \N\cup\{\infty\}$
and $H$ a finite-dimensional
(or more general) Lie group,
does $C^r_c(M,H)=\dl\,C^r_K(M,H)$\vspace{-.7mm}
hold as a $C^s$-manifold?
\end{numba}
\begin{numba}
In the fully general situation of Corollary~\ref{gemnonsilva}
(beyond the Silva case),
does $\Gamma(K,H)=\dl\, G_n$\vspace{-.4mm}
hold in the category of topological
spaces\,? Does it hold in the category
of $C^r$-manifolds, for $r\in \N_0\cup\{\infty\}$\,?
\end{numba}
\begin{numba}
Does the PTA always hold for test function groups
on $\sigma$-compact manifolds
and groups of compactly
supported diffeomorphisms\,?
Do the direct limit groups $G=\dl\,G_n$\vspace{-1mm}
described in Proposition~\ref{gdDLSilva}\,(ii)
always satisfy the PTA\,?
(The author cannot see a reason why this
should be true).
\end{numba}
\begin{numba}
Suppose that a topological group (resp., Lie group)
$G$ is the direct limit
of an ascending sequence $G_1\sub G_2\sub \cdots$
of topological groups (resp., Lie groups).
Will the product map
$\pi \colon \prod^*_{n\in \N}G_n\to G$
always be continuous (resp., smooth)\,?
If $G$ is a topological group (resp., Lie group),
what can we say about the product map
$\prod^*_{n\in \N}G\to G$\,?
\end{numba}
\begin{numba}
If $X$ is a topological space
(resp., smooth manifold)
in the preceding situation and
$f\colon X\to \prod_{n\in \N}^*G_n$
a function with continuous (resp., smooth)
components, find additional
conditions ensuring that $\pi\circ f$
is continuous (resp., smooth),
e.g.\ if $G_1=G_2=\cdots=A^\times$
for a unital Banach algebra~$A$.
Functions of the form considered in the proof
of Proposition~\ref{bumponsum}
should be subsumed as a special case.
\end{numba}
\begin{numba}
Let $G$ be a Lie group which is an ascending
union of Lie groups~$G_n$,
such that $L(G)=\dl\,L(G_n)$\vspace{-.4mm}
as a locally convex space.
There does not seem to be a reason why
$G$ should admit a direct limit chart
(not even if each $G_n$ is a Lie subgroup),
but the author does not know a counterexample.
\end{numba}
\begin{numba}
Let $G$ be a Lie group which is an ascending
union of Lie groups~$G_n$,
such that $G=\dl\,G_n$\vspace{-.4mm}
as a Lie group.
There does not seem
to be a reason why this should imply
that $L(G)=\dl\, L(G_n)$\vspace{-.4mm}
as a locally convex space,
but the author does not know any examples
where this fails.
Likewise,
there does not seem to be a reason why
$G=\bigcup_{n\in \N}G_n$
should admit a direct limit chart
if $G=\dl\,G_n$\vspace{-.4mm}
as a Lie group,
but the author does not know any
counterexamples.
\end{numba}
\begin{numba}
Find an example
of a direct limit topological group
$G=\bigcup_{n\in \N}G_n$
(preferably, a strict direct limit,
and preferably a direct limit of Lie groups
with a direct limit chart)
such that product sets (or two-sided
product sets) are not large
in $G$, resp., such that
$(G_n)_{n\in \N}$ does not satisfy the PTA.
Alternatively, prove that
such counterexamples do not exist.
\end{numba}
\begin{numba}
Find a direct sequence of Lie groups which does not have a
direct limit in the category of Lie
groups.
\end{numba}
\begin{numba}
Consider a Lie group $G=\bigcup_{n\in \N}G_n$
admitting a direct limit chart.
If $G=\dl\,G_n$\vspace{-1.1mm}
as a Lie group (or as a topological group),
does it follow that
$\wt{G}=\dl\,\wt{G}_n$\vspace{-.4mm}
as a Lie group (or as a topological group)\,?
\end{numba}
\begin{numba}
We obtained some positive results
concerning direct limit properties of complex
Lie groups, but could not prove any
general criteria to rule
out direct limit properties
in categories of complex Lie groups
or manifolds (due to the lack
of localization arguments).
Direct limit properties in the categories
of real analytic Lie groups (or manifolds)
are even more inaccessible.
\end{numba}
\begin{numba}
Let $(E_i)_{i\in I}$ be an uncountable
family of $C^r$-regular locally convex spaces.
Will the locally convex direct sum
$\bigoplus_{i\in I}E_i$ be $C^r$-regular\,?
\end{numba}
\begin{numba}
Let $M$ be a paracompact finite-dimensional
$C^\infty_\R$-manifold which is not $\sigma$-compact,
and $H$ be a finite-dimensional (or more general)
Lie group.
What can be said about the direct limit properties
of $\Diff_c(M)=\bigcup_K \Diff_K(M)$
and
$C^\infty_c(M,H)=\bigcup_K C^\infty_K(M,H)$\,?
\end{numba}
\appendix
\section{Smooth regularity of direct sums}\label{appdsum}
%
%
%\ma{appdsum}
%
%\ma{bumponsum}
\begin{prop}\label{bumponsum}
If $r\in \N \cup \{\infty\}$
and $(E_n)$
is a sequence of $C^r$-regular
locally convex spaces, then
the locally convex direct sum
$\bigoplus_{n\in \N} E_n$
is $C^r$-regular.
\end{prop}
\begin{proof}
A typical $0$-neighbourhood in $E:=\bigoplus_{n\in \N}E_n$
is of the form $U:=\bigoplus_{n\in \N}U_n$
for $0$-neighbourhoods $U_n\sub E_n$.
For each $n\in \N$,
let $f_n\colon E_n\to \R$
be a $C^r$-map such that
$f_n(E_n)\sub \;]0,1]$,
$f_n(0)=1$ and $f_n(x)\leq \frac{1}{2}$
for each $x\in E_n\setminus U_n$.
Then $f\colon E\to\R$,
\[
f(x):= f(x_1)\cdots f(x_n)\quad\mbox{for $\,x\in \bigoplus_{k=1}^n U_k$}
\]
is a function such that $f(E)\sub \;]0,1]$,
$f(0)=1$ and $f(x)\leq\frac{1}{2}$
for each $x\in E\setminus U$.
Using the
exponential function
$\exp\colon \R\to \;]0,\infty[$
and logarithm $\log=\exp^{-1}$,
we can write
\[
f\; =\; \exp\, \circ \, S\, \circ \, {\textstyle \bigoplus_{n\in \N}}
(\log \circ f_n)\,,
\]
where $\bigoplus_{n\in \N}(\log \circ f_n)\colon
\bigoplus_{n\in \N}E_n\to\bigoplus_{n\in \N}\R=\R^{(\N)}$
is~$C^r$ by \cite[Proposition~7.1]{MEA}
and the map
$S\colon \R^{(\N)}\to\R$, $(r_n)_{n\in \N}\mto \sum_{n\in \N}r_n$
is continuous linear and thus smooth.
Hence $f$ is~$C^r$.
Choose a $C^r$-map $g\colon \R\to\R$
such that $g(\R)\sub [0,1]$,
$g|_{[0,\frac{1}{2}]}=0$
and $g(1)=1$.
Then $g\circ f\colon E\to \R$
is a $C^r$-map which takes $0$ to~$1$
and vanishes outside~$U$.
The assertion follows.
\end{proof}
\vfill\pagebreak

\noindent{\footnotesize{\bf Helge Gl\"{o}ckner},
TU Darmstadt, Fachbereich Mathematik AG~5,
Schlossgartenstr.\,7,\\
64289 Darmstadt, Germany. E-Mail:
{\tt gloeckner@mathematik.tu-darmstadt.de}}

{\footnotesize
\begin{thebibliography}{MM}\itemsep+.1pc
%
% 
\bibitem{Aus} Au\ss{}enhofer, L.,
{\em Contributions to the duality theory of Abelian
topological groups and to the theory
of nuclear groups}, Dissertationes Math.\ {\bf 384},
1999.
%
%
\bibitem{Bzc}
Banaszczyk, W., ``Additive Subgroups of Topological
Vector Spaces,'' Springer-Verlag, 1991.
%
%
\bibitem{Ban}
Banyaga, A., ``The Structure of Classical
Diffeomorphism Groups,''
Kluwer, 1997.
%
%
\bibitem{Ber} Bertram, W., H. Gl\"{o}ckner
and K.-H. Neeb,
\emph{Differential calculus over general base fields and rings},
Expo.\ Math.\ \textbf{22} (2004), 213--282.
%
%
\bibitem{BaM}
Bierstedt, K.-D. and R. Meise,
\emph{Nuclearity and the Schwartz property in the theory of
holomorphic functions on metrizable locally convex spaces},
pp.\,93--129 in: M\,V. Matos et al.\ (eds.),
``Infinite-Dimensional Holomorphy and Applications,''
North-Holland, 1977.
%
%
\bibitem{BaS}
Bochnak, J. and J. Siciak,
\emph{Analytic functions in topological vector spaces},
Studia Math.\ {\bf 39}, 77--112.
%
%
\bibitem{BTV}
Bourbaki, N., ``Topological Vector Spaces,
Chapters~1--5,'' Springer,
Berlin, 1987.
%
%
\bibitem{Bou}
Bourbaki, N., ``Lie Groups and Lie Algebras,
Chapters~1--3,''
Springer, Berlin, 1989.
%
%
\bibitem{DaW}
Dierolf, S. and J. Wengenroth,
\emph{Inductive limits of topological algebras},
Linear Topol.\ Spaces Complex\ Anal.\ {\bf 3} (1997), 45--49.
%
%
\bibitem{Eda}
Edamatsu, T., \emph{On the bamboo-shoot topology of certain
inductive limits of topological groups},
J. Math.\ Kyoto Univ.\ {\bf 39} (1999), 715--724.
%
%
\bibitem{Flo}
Floret, K.,
\emph{Lokalkonvexe Sequenzen mit kompakten
Abbildungen}, J. Reine Angew.\
Math.\ {\bf 247} (1971), 155--195.
%
%
\bibitem{FaK}
Fr\"{o}licher, A. and A. Kriegl,
``Linear Spaces and Differentiation Theory,''
John Wiley, 1988.
%
%
\bibitem{RES} Gl\"{o}ckner, H.,
\emph{Infinite-dimensional Lie groups without
completeness restrictions},
pp.\,43--59 in: Strasburger, A. et al.\ (Eds.),
Geometry and Analysis on Finite- and Infinite-Dimensional
Lie Groups, Banach Center Publications~\textbf{55},
Warsaw, 2002.
%
%
\bibitem{GCX}
Gl\"{o}ckner, H.,
\emph{Lie group structures on quotient groups
and universal complexifications for infinite-dimensional
Lie groups}, J. Funct.\ Anal.\ \textbf{194} (2002), 347--409.
%
%
\bibitem{ALG}
Gl\"{o}ckner, H., \emph{Algebras whose groups of units
are Lie groups}, Studia Math.\ {\bf 153} (2002),
147--177.
%
%
\bibitem{DIR} Gl\"{o}ckner, H.,
\emph{Direct limit Lie groups and manifolds},
J. Math.\ Kyoto Univ.\ \textbf{43} (2003), 1--26.
%
%
\bibitem{MEA} Gl\"{o}ckner, H., {\em Lie groups of measurable mappings},
Canad.\ J. Math.\ {\bf 55} (2003), 969--999.
%
%
\bibitem{HOL} Gl\"{o}ckner, H., {\em Lie groups of germs of analytic mappings},
pp.\ 1--16 in T. Wurzbacher (Ed.),
``Infinite-dimensional groups and manifolds,''
IRMA Lecture Notes in Math.\ and Theor.\ Physics,
de Gruyter, 2004.
%
%
\bibitem{FUN} Gl\"{o}ckner, H.,
{\em Fundamentals of direct limit Lie theory},
Compos.\ Math.\ {\bf 141} (2005), 1551--1577.
%
%
\bibitem{HOE}
{\em H\"{o}lder continuous homomorphisms between infinite-dimensional
Lie groups are smooth},
J. Funct.\ Anal.\ {\bf 228} (2005), 419--444.
%
%
\bibitem{DIS}
Gl\"{o}ckner, H.,
{\em Discontinuous non-linear mappings
on locally convex direct limits},
Publ.\ Math.\ Debrecen {\bf 68} (2006), 1--13.
%
%
\bibitem{PAD}
Gl\"{o}ckner, H.,
{\em Every smooth $p$-adic Lie group admits a compatible
analytic structure}, Forum Math.\ {\bf 18} (2006), 45--84.
%
%
\bibitem{IMP}
Gl\"{o}ckner, H.,
{\em Implicit functions from topological
vector spaces to Banach spaces},
to appear in Israel J. Math.;
cf.\ arXiv:math.GM/0303320.
%
%
\bibitem{ZOO}
Gl\"{o}ckner, H.,
\emph{Lie groups over non-discrete topological fields},
preprint,
arXiv: math.GM/0408008.
%
%
\bibitem{IM2}
Gl\"{o}ckner, H.,
\emph{Finite order differentiability properties, fixed points
and implicit functions over valued fields},
preprint,
arXiv:math.FA/0511218.
%
%
\bibitem{SMA}
Gl\"{o}ckner, H.,
{\em Direct limit groups do not have
small subgroups},
preprint, arXiv: math.GR/0602407.
%
%
\bibitem{HOM}
Gl\"{o}ckner, H.,
{\em Homotopy groups of direct limits
of infinite-dimensional Lie groups},
in preparation.
%
%
\bibitem{SEC}
Gl\"{o}ckner, H.,
\emph{Differentiable mappings between spaces
of sections},
manuscript, 2002.
%
%
\bibitem{DIF}
Gl\"{o}ckner, H.,
\emph{Patched locally convex spaces,
almost local mappings
and diffeomorphism groups of non-compact
manifolds}, manuscript, 2002.
%
%
\bibitem{GaG}
Gl\"{o}ckner, H. and R. Gramlich,
\emph{Final group topologies,
Phan systems and Pontryagin duality},
preprint, arXiv:math.GR/0603537.
%
%
\bibitem{GaN}
Gl\"{o}ckner, H. and K.-H. Neeb,
``Infinite-Dimensional Lie Groups.
Vol.\,I: Basic Theory and Main Examples,''
book in preparation.
%
%
\bibitem{HaT}
Haller, S. and J. Teichmann,
{\em Smooth perfectness through decomposition
of diffeomorphisms into fiber preserving ones},
Ann.\ Global Anal.\ Geom.\ {\bf 23} (2003),
53--63.
%
%
\bibitem{Han}
Hansen, V.\,L., \emph{Some theorems on direct limits
of expanding systems of manifolds},
Math.\ Scand.\ \textbf{29} (1971), 5--36.
%
%
\bibitem{HaR}
Hewitt, E. and K.\,A. Ross,
``Abstract Harmonic Analysis I,''
Springer-Verlag, second edition,
1979.
%
%
\bibitem{Hir}
Hirai, T., H. Shimomura, N. Tatsuuma and E. Hirai,
\emph{Inductive limits of topologies, their direct product,
and problems related to algebraic structures},
J. Math.\ Kyoto Univ.\ \textbf{41} (2001),
475--505.
%
%
\bibitem{HaM}
Hofmann, K.\,H. and S.\,A. Morris,
``The Structure of Compact Groups,''
de Gruyter, Berlin, 1998.
%
%
\bibitem{Hog}
Hogbe-Nlend. H., ``Bornologies and Functional Analysis,''
North-Holland, Amsterdam, 1977.
%
%
\bibitem{Jar}
Jarchow, H.,
``Locally Convex Spaces,'' B.\,G. Teubner, Stuttgart,
1981.
%
%
\bibitem{KaM}
Kriegl, A. and P.~W.~Michor,
``The Convenient Setting of Global Ana\-lysis,''
AMS, Providence, 1997.
%
%
\bibitem{Lei}
Leitenberger, F.,
\emph{Unitary representations and coadjoint orbits
for a group of germs of real analytic diffeomorphisms},
Math.\ Nachr.\ {\bf 169} (1994), 185--205.
%
%
\bibitem{Mca}
Michael, E.\,A.,
``Locally Multiplicatively-Convex
Topological Algebras,''
Memoirs Amer.\ Math.\ Soc. \textbf{1952},
no.\,11, 1952.
%
%
\bibitem{Mic}
Michor, P.,
``Manifolds of Differentiable
Mappings,'' Shiva Publishing,
Orpington, 1980.
%
%
\bibitem{Mil}
Milnor, J., \emph{Remarks on infinite dimensional Lie groups},
in: B. DeWitt and R. Stora (eds.),
Relativity, Groups and Topology II,
North-Holland, 1984.
%
%
\bibitem{NRW1}
Natarajan, L., E.\ Rodr\'{\i}guez-Carrington
and J.\,A. Wolf,
\emph{Differentiable structure for direct limit groups},
Letters in Math.\ Phys.\ \textbf{23} (1991), 99--109.
%
%
%
\bibitem{NRW2}
Natarajan, L., E.\ Rodr\'{\i}guez-Carrington
and J.\,A. Wolf, \emph{Locally convex Lie groups},
Nova J. Alg.\ Geom.\ \textbf{2} (1993), 59--87.
%
%
\bibitem{NeI}
Neeb, K.-H., \emph{Infinite-dimensional groups
and their representations},
pp.\,131--178 in: A. Huckleberry and T. Wurzbacher
(Eds.), ``Infinite-Dimensional K\"{a}hler
Manifolds,'' Birkh\"{a}user, Basel, 2001.
%
%
\bibitem{NeC}
Neeb, K.-H.,
{\em Central extensions of infinite-dimensional
Lie groups}, Ann.\ Inst.\ Fourier (Grenoble)
{\bf 52} (2002), 1365--1442.
%
%
\bibitem{NeA}
Neeb, K.-H.,
{\em Abelian extensions
of infinite-dimensional Lie groups},
Travaux Math.\ {\bf 15} (2004), 69--194.
%
%
\bibitem{NeN}
Neeb, K.-H.,
{\em Non-abelian extensions of infinite-dimensional
Lie groups}, erscheint in Ann.\ Inst.\ Fourier
(Grenoble); cf.\ arXiv:math.GR/0504295.
%
%
\bibitem{Pie}
Pietsch, A., ``Nuclear Locally Convex Spaces,''
Springer-Verlag, 1972.
%
%
\bibitem{Pis}
Pisanelli, D.,
\emph{An example of an infinite
Lie group}, Proc.\ Amer.\ Math.\ Soc.\
\textbf{62}, (1977), 156--160.
%
%
\bibitem{Sch}
Schaefer, H.\,H. and M\,P. Wolff,
``Topological Vector Spaces,''
Springer-Verlag, 1999.
%
%
\bibitem{TSH}
Tatsuuma, N., H. Shimomura and T. Hirai,
\emph{On group topologies and unitary representations of inductive
limits of topological groups and the case of the group of diffeomorphisms},
J. Math.\ Kyoto Univ.\ \textbf{38} (1998), 551--578.
%
%
\bibitem{Tre} Tr\`{e}ves, F.,
``Topological Vector Spaces, Distributions
and Kernels,'' Academic Press, 1967.
%
%
\bibitem{Wen}
Wengenroth, J., personal communication, 30.7.2003.
%
%
\bibitem{Yam}
Yamasaki, A., \emph{Inductive limits of general linear groups},
J. Math.\ Kyoto Univ.\ {\bf 38} (1998), 769--779.\vspace{2mm}
%
%
\end{thebibliography}}
\end{document}